\documentclass[11pt,a4paper,twoside]{article}
\usepackage{amsmath,amsfonts,amssymb,amsthm,graphics,hyperref}
\usepackage[mathscr]{euscript}
\usepackage{xy}
\xyoption{all}

\newtheorem{Theorem}[equation]{Theorem}
\newtheorem{Corollary}[equation]{Corollary}
\newtheorem{Proposition}[equation]{Proposition}
\newtheorem{Lemma}[equation]{Lemma}
\newtheorem{Remark}[equation]{Remark}

\newtheorem{Definition}[equation]{Definition}

\def\Section#1{\section{#1}\setcounter{equation}{0}}

\def\bdots{\mathinner{\mkern1mu\raise1pt\hbox{.}\mkern2mu\raise4pt\hbox{.}
           \mkern2mu\raise7pt\vbox{\kern7pt\hbox{.}}\mkern1mu}}

\def\ov#1{\overline {#1}}
\def\til#1{\widetilde {#1}}

\def\Indu#1#2#3{\Ind_{#1}^{#2}{#3}}
\def\Resu#1#2#3{\Res^{#1}_{#2}{#3}}

\def\bs{\boldsymbol}

\def\hbk{\hfill\break}

\def\today{\number\day\space
 \ifcase\month\or
	January\or February\or March\or April\or May\or June\or
	July\or August\or September\or October\or November\or December\fi
 \space\number\year}

\def\tG{\til G}
\def\tH{\til\HH}
\def\tJ{\til\JJ}
\def\tK{\til\KK}
\def\tM{\til M}
\def\tP{\til P}
\def\tU{\til U}

\def\tth{{\til\theta}}
\def\teta{\til\eta}

\def\sm{{\mathsf m}}
\def\sM{{\mathsf M}}
\def\so{{\mathsf o}}

\def\of{{\mathfrak o}_F}

\def\pf{{\mathfrak p}_F}
\def\kf{{k_F}}
\def\pif{{\varpi}_F}

\def\oe{{\mathfrak o}_E}

\def\pe{{\mathfrak p}_E}
\def\ke{{k_E}}
\def\pie{{\varpi}_E}

\def\oo{{\mathfrak o}_0}

\def\po{{\mathfrak p}_0}

\def\pio{{\varpi_0}}

\def\fM{{\mathfrak M}}

\def\fH{{\mathfrak H}}
\def\fJ{{\mathfrak J}}

\def\fa{{\mathfrak a}}
\def\fb{{\mathfrak b}}

\def\fo{{\mathfrak o}}
\def\fp{{\mathfrak p}}

\def\ko{{k_0(\beta,\Lambda)}}

\def\End{{\hbox{\rm End}}}
\def\Aut{{\hbox{\rm Aut}}}
\def\Hom{{\hbox{\rm Hom}\,}}
\def\max{\hbox{\rm max}\,}

\def\dim{\hbox{\rm dim}\,}
\def\ker{\hbox{\rm ker}\,}

\def\det{\hbox{\rm det}\,}

\def\cInd{{\hbox{\rm c-Ind}\,}}
\def\Ind{{\hbox{\rm Ind}\,}}
\def\Res{{\hbox{\rm Res}\,}}

\def\Supp{{\hbox{\rm Supp}\,}}
\def\card{\hbox{\rm card}\,}
\def\diag{\hbox{\rm diag}}

\def\lcm{\hbox{\rm lcm}\,}
\def\Lie{\hbox{\rm Lie}\,}
\DeclareMathOperator{\res}{res}

\def\GL{{GL}}

\def\tr{{\hbox{\rm tr}}}

\def\into{\hookrightarrow}
\def\ab{a_\beta}
\def\slto#1{\xrightarrow{#1}}

\def\smap{\slto{\ s\ }}
\def\abmap{\slto{\,\ab\,}}
\def\nomap{\slto{\ \ \ }}

\def\le{\leqslant}
\def\ge{\geqslant}

\def\subset{\subseteq}
\def\supset{\supseteq}

\def\CB{{\mathcal B}}
\def\CC{{\mathcal C}}

\def\CG{{\mathcal G}}
\def\CH{{\mathcal H}}
\def\CI{{\mathcal I}}

\def\CM{{\mathcal M}}

\def\CP{{\mathcal P}}
\def\CT{{\mathcal T}}
\def\CU{{\mathcal U}}

\def\BF{{\mathbb F}}

\def\BN{{\mathbb N}}

\def\BZ{{\mathbb Z}}

\def\BCx{{\mathbb C^\times}}

\def\boe{{\mathbf e}}

\def\bok{{\mathbf k}}

\def\boI{{\mathbf 1}}

\def\HH{{\mathrm H}}
\def\JJ{{\mathrm J}}
\def\KK{{\mathrm K}}
\def\PP{{\mathrm P}}

\def\a{{\alpha}}
\def\b{{\beta}}
\def\g{{\gamma}}

\def\ve{{\epsilon}}
\def\e{{\varepsilon}}
\def\k{{\kappa}}
\def\l{{\lambda}}
\def\s{{\sigma}}
\def\vth{{\vartheta}}
\def\th{{\theta}}
\def\z{{\zeta}}

\def\D{{\Delta}}
\def\Ga{{\Gamma}}
\def\L{{\Lambda}}
\def\Si{{\Sigma}}

\begin{document}

\title{The supercuspidal representations of $p$-adic classical groups}
\author{Shaun Stevens\thanks{The research for this paper was supported
    by EPSRC grant GR/T21714/01.}} 
\date{\today}
\maketitle
\begin{abstract}
Let $G$ be a unitary, symplectic or special orthogonal group over a 
locally compact non-archimedean local field of odd residual
characteristic. We construct many new supercuspidal representations 
of $G$, and Bushnell-Kutzko types for these representations. Moreover, 
we prove that every irreducible supercuspidal representation of $G$ 
arises from our constructions.
\end{abstract}
%


\section*{Introduction}

Let $G$ be a unitary, symplectic or special orthogonal group over a
locally compact non-archimedean local field of odd residual
characteristic. This article completes a series of papers (in
particular~\cite{S3,S4}) whose major goal was the construction of all
irreducible supercuspidal representations of $G$.

More precisely, in this paper, we construct pairs $(J,\l)$, where $J$
is a compact open subgroup of $G$ and $\l$ is an irreducible
representation of $J$, such that the compactly induced representation
$\cInd_J^G\l$ is irreducible and supercuspidal. Moreover, any
irreducible supercuspidal representation $\pi$ of $G$ contains some
such pair, in the sense that $\l$ is a component of the restriction of
$\pi$ to $J$, and hence $\pi\simeq\cInd_J^G\l$. In particular, this
shows that every irreducible supercuspidal representation of $G$ can
be obtained by irreducible compact induction from a compact open
subgroup. 

\medskip

The study of the representation theory of $p$-adic groups by
restriction to compact open subgroups goes back to the work of
Howe~\cite{H} and, for the general linear group in arbitrary (residual)
characteristic, was completed by Bushnell and Kutzko in~\cite{BK}. 
For classical groups $G$ in odd residual characteristic, we can go
back to the work of Moy on $U(2,1)$~\cite{Moy1} and
$GSp_4$~\cite{Moy2,Moy3}, and, for a general classical group, of
Morris in a series of papers (see, for example,~\cite{Mor2,Mor4}) which laid
much of the necessary groundwork.

More recently, results have been obtained for arbitrary connected
reductive groups. The level-zero representations -- that is, those
representations which contain the trivial character on restriction to
the pro-unipotent radical of some parahoric subgroup -- were classified
independently by Morris~\cite{Mor3} and Moy and
Prasad~\cite{MP}. For positive-level representations in the so-called
\emph{tame} case, the most general constructions of supercuspidal
representations are due to Yu~\cite{Yu}. 

Recently in~\cite{Kim}, Kim
has proved that, under some restrictive conditions on the ground field
(dependent on the group $G$), Yu's constructions do give all
supercuspidal representations. 
Although our work here is not for such a general class of group, it
has the notable advantage that the only restriction on the field is
that the residual characteristic be odd.

\bigskip

Our work for classical groups $G$ follows the spirit, and the methods,
of Bushnell and Kutzko's work~\cite{BK}. So we begin with a
non-archimedean local field $F$ of odd residual characteristic
with a (possibly trivial) galois involution with fixed field
$F_0$, and a finite-dimensional $F$-vector space $V$ equipped with a
nondegenerate $\e$-hermitian form $h$. Then $G$ is the group of
$F_0$-rational points of the connected component of the linear
algebraic group preserving the form. It will also be useful to have
the full group $G^+$ of automorphisms of $V$ preserving the form. 

The construction of supercuspidal representations
begins with certain local data, namely a \emph{skew
semisimple stratum $[\L,n,0,\b]$} (see~\S\ref{S2.1}). Here we recall 
only that $\L$ is a lattice sequence in $V$ and 
$\b$ is an endomorphism of $V$ which is skew for the adjoint
involution induced by $h$ and such that $E=F[\b]$ is a direct sum of
fields each preserved by this involution -- that is, $\b$ is an
elliptic semisimple element of the Lie algebra of $G$.
We write $G_E$ for the centralizer of $\b$ in $G$; this is a product
of classical groups over extensions of $F_0$.
The lattice sequence $\L$
defines a compact open subgroup $\PP(\L)$ in $G$, together with a
filtration, and we write $\PP(\L_{\oe})$ for the intersection
$\PP(\L)\cap G_E$. 

\medskip

In~\cite{S4}, we defined the set $\CC_-(\L,0,\b)$ of \emph{skew
semisimple characters} attached to 
a skew semisimple stratum, generalizing the notion of \emph{simple
character} from~\cite{BK}. These are certain arithmetically defined
abelian characters $\th$ of a compact open subgroup $\HH^1(\b,\L)$
with nice intertwining properties. Moreover, for fixed $\b$ but
varying lattice sequence $\L$, these semisimple characters exhibit
remarkable functorial properties, known as \emph{transfer
properties} -- again, these generalize such properties for general
linear groups. 

Attached to our stratum, there are two other compact open subgroups
$\JJ(\b,\L)\supset \JJ^1(\b,\L)$ containing $\HH^1(\b,\L)$ and
normalizing each semisimple character $\th$. Moreover, it is shown
in~\cite{S4} that there is a unique irreducible representation $\eta$
of $\JJ^1(\b,\L)$ (a \emph{Heisenberg extension}) containing a given
semisimple character $\th$, that the intertwining of $\eta$ is the
same as that of $\th$ and, moreover, that every intertwining space is
$1$-dimensional. 

\medskip

Our first task in this article is to extend the representation $\eta$
to the group $\JJ(\b,\L)$ -- this is known as the
\emph{$\b$-extension} problem. In fact, finding such an extension is
not difficult -- given the formalism and transfer properties of  
semisimple characters, the same arguments as in the case of general
linear groups in~\cite{BK} can be applied. The major difficulty is in
knowing whether we have a \emph{suitable} extension -- in general,
many of the extensions are not at all appropriate for our purposes.

For $GL_n(F)$ (or, more generally, general linear groups over a
division algebra -- see~\cite{Sec2}) the crucial property of
$\b$-extensions is that their intertwining is the same as that of the
representation $\eta$ which they extend. In our situation for
classical groups, we have not been able to prove that any extensions
have this property. Instead, we characterize $\b$-extensions by 
certain compatibility (or transfer) properties, beginning
with the case where the group $\PP(\L_{\oe})$ is a maximal compact
open subgroup of $G_E$,
analogously to the techniques in~\cite{BK}. However, even these
compatibilities have subtleties which do not appear for $GL_n(F)$.

In the case where the order $\PP(\L_{\oe})$ is maximal, we
actually extend $\eta$ to the group $\JJ^+(\b,\L)$, which is the group
analogous to $\JJ(\b,\L)$ in $G^+$ -- the representation of
$\JJ(\b,\L)$ we seek will then be given by restriction. The methods
of~\cite{BK} 
generalize to our case and we can construct an extension $\k$
of $\eta$ whose restriction to a pro-$p$ Sylow subgroup of
$\JJ^+(\b,\L)$ is determined by the transfer properties of semisimple
characters. Indeed, in most cases any extension will do -- the only
exception is in the case where one of the simple
components of the quotient $\JJ^+(\b,\L)/\JJ^1(\b,\L)$, which is a
product of reductive groups over finite fields but may be
disconnected, is isomorphic to $SL_2(\BF_3)$. This is similar to the
situation for $GL_n(F)$, where the only difficulty is with $GL_2(\BF_2)$.

\medskip

For general $\L$, we choose a lattice sequence $\L^\sM$ such that
$\PP(\L^\sM_{\oe})$ is a maximal compact subgroup of $G_E$ containing
$\PP(\L_{\oe})$ and let $\k_\sM$ be a $\b$-extension to $\JJ^+(\b,\L^\sM)$
corresponding to the semisimple character in $\CC_-(\L^\sM,0,\b)$
which is the transfer of $\th$.
We observe here that, unlike the case of $GL_n(F)$, this
does \emph{not} imply that the corresponding subgroups of $G$
satisfy $\PP(\L^\sM)\supset\PP(\L)$. This complication has already
been encountered for general linear groups over division algebras
in~\cite{Sec2}. In the case where we do have this inclusion, we may
follow~\cite{BK} and define a $\b$-extension $\k$ to $\JJ^+(\b,\L)$ to be
the (unique) 
one which satisfies an induction compatibility with respect to $\k_\sM$,
namely
$$
\Indu{\JJ^+(\b,\L)}{\PP^+(\L^\sM)}{\k}\ \simeq\ 
\Indu{\PP^+(\L_{\oe})\JJ^1(\b,\L^\sM)}{\PP^+(\L^\sM)}
   {\left(\k_\sM|_{\PP^+(\L_{\oe})\JJ^1(\b,\L^\sM)}\right)},
$$
where $\PP^+(\L)$ is the compact open subgroup of $G^+$ defined by $\L$.
When we do not have the containment of subgroups, we imitate~\cite{Sec2}
and proceed in steps from $\L$ to $\L^\sM$ by tracing a line through
the Bruhat-Tits building of $G$ and cutting it into sufficiently small
pieces. Here we are using the description of the building in terms of
\emph{lattice functions} (see~\cite{BS} or~\S\ref{S2.2}).

\medskip

There is an important difference with previous work to note here,
namely that the different possible choices for $\L^\sM$
are not conjugate in $G^+$ (and, in general, not in the group of
similitudes either). This gives the added complication that the
definition of $\b$-extension depends \emph{a priori} on the choice of
$\L^\sM$. Indeed, it is easy to see that different choices can give
different \emph{numbers} of $\b$-extensions. It is not clear whether
or not there is, in general, an extension of $\eta$ to $\JJ^+(\b,\L)$
which is a $\b$-extension for every choice of $\L^\sM$ and this is at
the heart of the problem of computing intertwining. On the other hand,
we do show (Corollary~\ref{betatwist}) that $\b$-extensions for different
choices of $\L^\sM$ cannot differ by too much.
It is essentially this same problem which prevents us from adopting
the ``bottom-up'' approach to $\b$-extensions of~\cite{Sec2},
beginning with the case when $\PP(\L_{\oe})$ is an Iwahori subgroup.

There is, nonetheless, a ``canonical'' choice for $\L^\sM$. (In fact,
there are two choices which could be described as canonical, and we
choose one of them -- see \eqref{standardmax}.) We call
$\b$-extensions defined relative to this canonical choice
\emph{standard $\b$-extensions}. The point here is that
$\PP(\L^\sM_{\oe})$ is chosen so that the Weyl group
of its reductive quotient is as large as possible -- so it is ``as
close as possible'' to being special. Using this large Weyl group, and
an Iwahori decomposition with respect to a suitable Levi subgroup of
$G$, we
are able to compute a lower bound for the intertwining of standard
$\b$-extensions (see Proposition~\ref{intkPstandard}), which is good
enough for our purposes.

\medskip

Now we add a level zero piece. Let
$\JJ^\so(\b,\L)$ denote the inverse image in $\JJ(\b,\L)$ of the
connected component of the reductive quotient
$\JJ(\b,\L)/\JJ^1(\b,\L)\cong \PP(\L_{\oe})/\PP^1(\L_{\oe})$;
likewise, we write $\PP^\so(\L_{\oe})$ for the inverse image in
$\PP(\L_{\oe})$, which is a parahoric subgroup of $G_E$. 
Let $\rho$ be the inflation to $\JJ^\so(\b,\L)$ of an irreducible
\emph{cuspidal} representation of the quotient
$\PP^\so(\L_{\oe})/\PP^1(\L_{\oe})$ and set $\l=\k\otimes\rho$. 

We can now bound above the $G$-intertwining of $\l$ (see \S~\ref{S3.2}).
The key point is that the intertwining of the
restriction of $\k$ to a pro-$p$ Sylow subgroup is the same as that
of~$\eta$ so we can first bound the intertwining in terms of the
intertwining of the restriction of $\rho$ to a pro-$p$ Iwahori
subgroup of $\PP^\so(\L_{\oe})$. Now this intertwining can be bounded
by a slight weakening of the hypotheses of a result of Morris in~\cite{M1} 
(see Proposition~\ref{Morlem}). The idea that knowing the
intertwining of the restriction to the pro-$p$ Iwahori can be enough
to bound intertwining has been used already in a special case by
Blondel (see~\cite[page 553]{Bl1}). Of course, this only works because
the level zero piece $\rho$ is \emph{cuspidal}.

\medskip

In the special case where $\PP(\L_{\oe})$ is a maximal compact
subgroup of $G_E$, the intertwining of $\l$ is found to be just
$\JJ(\b,\L)$. In particular, if $\tau$ is any irreducible
representation of $\JJ(\b,\L)$ whose restriction to $\JJ^\so(\b,\L)$
contains cuspidal $\rho$, then the compactly induced representation
$$
\pi\ = \ \cInd_{\JJ(\b,\L)}^G (\k\otimes\tau)
$$
is irreducible and hence supercuspidal. Moreover, the pair
$(\JJ(\b,\L),\k\otimes\tau)$ is a $[G,\pi]_G$-type, in the sense of
Bushnell-Kutzko~\cite{BK2}. 

\medskip

Finally, it remains to show that we have constructed every
supercuspidal representation of $G$ in this way. The ideas are similar 
to those in~\cite{BK} (see also~\cite{SS} for the case of general linear 
groups over division algebras) though the details are somewhat more 
complicated. From~\cite{S4}, we
know already that any irreducible supercuspidal representation $\pi$
contains a semisimple character $\th$. Then, possibly after changing
the lattice sequence $\L$, it is not hard to see that $\pi$ must
contain a representation $\vth=\k\otimes\rho$ of $\JJ^\so(\b,\L)$, for
$\k$ a standard $\b$-extension and 
$\rho$ a \emph{cuspidal} representation. The idea then is to show
that, if the compact subgroup $\PP(\L_{\oe})$ is not maximal then
$\pi$ has a non-zero Jacquet module, contradicting the cuspidality of
$\pi$.

In order to construct such a Jacquet module we use Bushnell-Kutzko's
method of \emph{covers}. There is a (choice of) 
Levi subgroup $M$ associated to the 
stratum such that, for any parabolic subgroup $P=MU$ with Levi factor $M$, 
the group $\JJ^\so(\b,\L)$ has an Iwahori decomposition with respect to 
$(M,P)$. We form the group 
$$
\JJ^\so_P\ =\ \HH^1(\b,\L) \left(\JJ^\so(\b,\L)\cap P\right)
$$
and let $\vth_P$ be the natural representation of $\JJ^\so_P$ on the 
$(\JJ^\so(\b,\L)\cap U)$-fixed vectors of $\vth$. Then $\pi$ certainly 
also contains $\vth_P$. Moreover, the pair 
$(\JJ^\so_P,\vth_P)$ is a \emph{decomposed pair} -- that is, $\vth_P$ 
restricts trivially to the unipotent radical of any parabolic subgroup 
with Levi component $M$, and the restriction $\vth_P|_{\JJ^\so_P\cap M}$ 
is irreducible.

In order to show that $(\JJ^\so_P,\vth_P)$ is a cover, it remains to 
show that there is an invertible element of the Hecke algebra 
$\CH(G,\vth_P)$ supported on $\JJ^\so_P\z\JJ^\so_P$, for some 
\emph{strongly $(P,\JJ^\so_P)$-positive element} $\z$ of the centre of 
$M$. Then the theory of covers shows that $\pi$ has a non-zero 
Jacquet module.

\medskip

We note that, since the compact open subgroup $\PP(\L_{\oe})$ is 
non-maximal, the reductive quotient $\PP^\so(\L_{\oe})/\PP^1(\L_{\oe})$ 
has a simple component isomorphic to some general linear group (over an 
extension of the residue field of $F_0$). There is then an involution on 
these simple components, induced by the Weyl group 
elements introduced earlier. We will need to divide into two cases:

\smallskip

(i) If the component of $\rho$ on some such general linear group is 
\emph{not} fixed by the involution, then the bound on intertwining 
from above gives the bound
$$
I_G(\vth_P)\ \subset\ \JJ^\so_P M' \JJ^\so_P,
$$
for some proper Levi subgroup $M'$ of $G$ containing $M$. Then it is 
straightforward that $(\JJ^\so_P,\vth_P)$ is a cover of 
$(\JJ^\so_P\cap M',\vth_P|_{\JJ^\so_P\cap M'})$ (see~\S\ref{5.2}).

\medskip

(ii) Otherwise every component of $\rho$ is fixed by the involution 
and we use an idea of Kutzko (which was also used in~\cite{GKS}). 
If $\L^\sM$ is a lattice sequence such that $\PP(\L^\sM_{\oe})$ is a 
maximal compact subgroup of $G_E$ containing $\PP(\L_{\oe})$ then, 
thinking of $\rho$ as a representation of $\PP(\L_{\oe})$, 
we get a support-preserving injection of Hecke algebras 
(see~\eqref{injalgmap})
$$
\CH(\PP(\L^\sM_{\oe}),\rho\otimes\chi)\ \hookrightarrow 
\CH(G,\vth_P),
$$
for some self-dual character $\chi$ (which comes from the problem of 
incompatibility of $\b$-extensions relative to different maximal 
compact subgroups of $G_E$). Moreover, the Hecke algebra on the left, 
which is just in level zero, is understood via the very general results 
of Morris in~\cite{M1}. In particular, it contains an invertible 
function supported on a certain Weyl group element (see~\S\ref{5.3} for 
details).
If we do this for two choices of $\L^\sM$, we obtain two invertible 
elements of $\CH(G,\vth_P)$ whose convolution is supported on a single 
double coset $\JJ^\so_P\xi\JJ^\so_P$, with $\xi\in M$. Raising this to 
a suitable power, we obtain the invertible function we were seeking.

\medskip

We end with some remarks on questions which are still unanswered. 
Firstly, although we present here a construction which gives all 
irreducible supercuspidal representations of $G$, we do not know 
when two apparently distinct constructions give rise to equivalent 
representations. More precisely, we do not have an ``intertwining 
implies conjugacy'' result  for the types we construct -- indeed, 
we do not even have such a result for the strata underlying the types.

Secondly, there still remains the problem of constructing types for 
the non-supercuspidal components of the Bernstein spectrum. Types have 
been constructed in certain special cases: see~\cite{Bl1} and~\cite{GKS}, 
where types are constructed for self-dual supercuspidal representations of 
the Siegel Levi subgroup; and~\cite{Bl2}, where (sometimes putative) covers 
in smaller symplectic groups are ``propagated'' to larger ones -- note that 
our work here shows that the hypothesis (H4) of~\cite{Bl2} can always be 
satisfied.
There are also some results in this direction in~\cite{Dat}, where it 
is shown that any smooth irreducible representation of $G$ contains 
some (suitably generalized) semisimple character. We remark that, 
throughout this work and its predecessor~\cite{S4}, numerous covers 
are constructed. It seems likely that a more careful examination, 
especially of the Hecke algebras in level zero used in proving 
exhaustion, should give types in all cases.

\medskip

{\bf Acknowledgements} \ I would like to thank: David Goldberg, and 
Purdue University, for their hospitality at an early stage of this work; 
the University of East Anglia, for granting me study leave, which gave 
me the chance finally to write up all the details of the arguments;
and, in no particular order, Colin Bushnell, Guy Henniart, 
Corinne Blondel, Phil Kutzko and Vincent S\'echerre, 
for their support, advice, and interest in this work at various times.
I would also like to thank Corinne Blondel, Michitaka Miyauchi and the referee for their careful reading, and for pointing out some errors and missing explanations.


\tableofcontents


\Section{A variant of a result of L.~Morris}\label{S1}

In this section, we give a mild weakening of the hypotheses in a
result of Lawrence Morris (\cite[Proposition 4.13]{M1}). Since this
result is valid in great generality, we retain all the notation of
\emph{op.\ cit.}; in particular, the notation of this section is
independent of that in the rest of the paper. We also refer the 
reader to \emph{op.\ cit.} (or to \cite{BT1}, \cite{BT2}) for details 
of the situation, which we recall briefly below.

\medskip

Let $\boldsymbol G$ be a connected reductive group over a
non-archimedean local field $F$ and let $G=\boldsymbol G(F)$ be its
group of $F$-points. Let $\boldsymbol T$ be a maximal $F$-split torus
in $\boldsymbol G$ and let $\boldsymbol N$ be its normalizer; we also 
write $T=\boldsymbol T(F)$ and $N=\boldsymbol N(F)$. Let $B$
be an Iwahori subgroup of $G$, which is the stabilizer of a chamber in
the apartment of the enlarged Bruhat-Tits building of $G$ associated
to $T$. Then $(G,B,N)$ is a generalized affine BN-pair. 

Associated to this affine BN-pair, we have the generalized affine 
Weyl group $W=N/B\cap N$; it is the semi-direct product of a 
subgroup $\Omega$ and a normal subgroup $W'$ which is the affine Weyl 
group for some split affine root system $\Sigma$. We also have a set 
$S$ of fundamental reflections in $W'$ corresponding to a simple 
system of affine roots.

Given any proper subset $J\subset S$, we let $W_J$ be the subgroup of $W$ 
generated by the reflections in $J$ and let $N_J$ be the 
corresponding subgroup of $N$. Then $P_J=BN_JB$ is a subgroup of $G$, 
which we call a \emph{standard parahoric subgroup\/}. For $J,K$ subsets 
of $S$, Morris defines in \cite[\S3.10]{M1} (see also \emph{loc.\
  cit.} Variant 3.22) a 
(non-unique, in general) set of \emph{distinguished double coset 
representatives} for $W_J\backslash W/W_K$, which we denote $D_{J,K}$. 
We can then choose $D'_{J,K}$ to be a set of inverse images in $N$ of 
$D_{J,K}$, and this is a set of \emph{distinguished double coset 
representatives} for $P_J\backslash G/P_K$.

Given $P_J$ a standard parahoric subgroup of $G$, we denote by $U_J$ 
its pro-unipotent radical, and by $M_J$ the quotient $P_J/U_J$, (the
points of) a connected
reductive group over a finite field. We write $U_B$ for the pro-radical 
$U_{\emptyset}$ of $B=P_{\emptyset}$.

\medskip

Finally, for $\rho$ a representation of a compact subgroup $K$ of $G$, 
we write $\CH(\rho|K)=\End_G\left(\cInd_K^G\rho\right)$ for the
associated Hecke algebra and 
$I_G(\rho|K)$ for the intertwining of $\rho$, which is the support of
the Hecke algebra.

\begin{Proposition}[cf. {\cite[Proposition 4.13]{M1}}] Let $J$ be a 
proper subset of $S$, let $D'_{J,J}$ be a set of distinguished double coset 
representatives for $P_J\backslash G/P_J$ and let $n\in D'_{J,J}$. 
Let $\rho$ be the inflation to $P_J$ of a cuspidal representation 
of $M_J$. If $n$ lies in the support of $\CH(\rho|U_B)$, then $wJ=J$, 
where $w\in W$ is the projection of $n$.
\label{Morlem}
\end{Proposition}

\noindent{\bf Remark} The only difference between this 
and~\cite[Proposition 4.13]{M1}
is that the hypothesis there is that $n$ lies in 
the support of $\CH(\rho|P_J)$.

\begin{proof} Write $P=P_J$ and let $P'$ denote $nPn^{-1}$, with
radical $U'$, and $U_B'=nU_Bn^{-1}$. We note that we have $P\cap~U'\subset~
U_{J\cap wJ}\subset U_B$, by~\cite[Lemma 3.21]{M1}, and $U'\subset 
U'_B$ so we have $P\cap U'\subset U_B\cap U'_B$. Then, if we replace 
every $P\cap P'$ in the proof of \emph{op.\ cit.} Proposition 4.13 by 
$U_B\cap U'_B$, the proof is identical.
\end{proof}

The interest in weakening the hypothesis to one which says only that
the restriction of $\rho$ to $U_B$ is intertwined by $n$ is suggested by
the work of Blondel in~\cite{Bl1}, in which she is able to compute
sufficient information on intertwining without calculating it
completely (\emph{op.\ cit.} Proposition IV.1 and Lemma IV.2). We will
wish to use it in two different cases:
\begin{enumerate}
\item The case when $P=P_J$ is a non-maximal parahoric subgroup of
$G$. In this case, let $^{c}\Phi_J$ denote the closure of the root
system generated by the gradients of the affine roots in $J$ and let
$\fM_J=\bs\fM_J(F)$ be the group of points of the corresponding
reductive subgroup of $\bs G$ (see~\cite[\S3.15]{M1}). 
Then $\fM_J\cap P_J$ is a parahoric subgroup of $\fM_J$ with reductive
quotient $M_J$ again.
We also put $N_J=N\cap\fM_J\cap P_J$. Then a set
of representatives for $\{n\in N:wJ=J\}/N_J$ can be taken in
the normalizer $N_N(\fM_J\cap P_J)$ of $\fM_J\cap P_J$ in $N$
(see~\cite[\S4.16]{M1}).  
\item The case where $P=P_J$ is a maximal parahoric subgroup of
$G$. In this case, the conclusion from the previous lemma is that
$n$ normalizes $P$ (see \cite[Appendix A.1]{M1}). In the classical groups
we will be looking at, the $G$-normalizer of $P$ is compact.
\end{enumerate}

We will also need to use the following
stronger result from~\cite{M1}: 

\begin{Theorem}[{\cite[Theorem 4.15]{M1}}]\label{Morlem2}
Let $J$ be a proper subset of $S$ and let $\rho$ be the inflation to
$P_J$ of a cuspidal representation  
of $M_J$. Then $I_G(\rho)\subset P_J N_N(\rho) P_J$, where
$N_N(\rho)=\{n\in N_N(\fM_J\cap P_J):{}^n\rho\simeq\rho\}$.
\end{Theorem}


\Section{Notations and preliminaries}\label{S2}

After setting up the notation for the remainder of this paper, we will
recall some of the basic objects needed in the theory. We give quite a
rapid overview, and refer the reader to~\cite{BK} and~\cite{S4} for
more details.


\subsection{Strata}\label{S2.1}

Let $F$ be a non-archimedean local field equipped with a galois
involution $\ov{\vphantom{a}\ }$ with fixed field $F_0$; we allow the
possibility $F=F_0$. Let $\of$ be the ring of integers of $F$, $\pf$
its maximal ideal and $\kf=\of/\pf$ the residue field; we assume
throughout that the residual characteristic $p:=\hbox{\rm char}\,\kf$ is not
$2$. We denote by $\oo$, $\po$, $k_0$ the same objects in $F_0$, and will
use similar notation for any non-archimedean local field. We fix a
uniformizer $\pif$ of $F$ such that $\ov\pif=-\pif$ if $F/F_0$ is
ramified, $\ov\pif=\pif$ otherwise. We put $\pio=\pif^2$ if $F/F_0$ is
ramified, $\pio=\pif$ otherwise; so $\pio$ is a uniformizer of $F_0$.
We will use analogous notation for any field extension $E$ of $F$ to which
the involution extends.

\medskip

Let $V$ be an $N$-dimensional vector space over $F$, equipped with a
nondegenerate $\epsilon$-hermitian form, with $\epsilon=\pm 1$. We put
$A=\End_F V$ and denote by $\ov{\vphantom{a}\ }$ the adjoint 
(anti-)involution on $A$ induced by $h$. Set also $\tG=\Aut_F V$ and let
$\sigma$ be the involution given by $g\mapsto\ov g^{-1}$, for
$g\in\tG$. We also have an action of $\sigma$ on the Lie algebra $A$
given by $a\mapsto -\ov a$, for $a\in A$. We put
$\Sigma=\{1,\sigma\}$, where $1$ acts as 
the identity on both $\tG$ and $A$.

We put $G^+=\tG^\Sigma=\{g\in\tG:h(gv,gw)=h(v,w)$ for all $v,w\in V\}$,
the $F_0$-points of a unitary, symplectic or orthogonal group $\bs G^+$ 
over $F_0$. Let $G$ be the $F_0$-points of the connected component 
$\bs G$ of $\bs G^+$, so that $G=G^+$ except in the orthogonal case. 

We put $A_-=A^\Sigma\simeq\Lie G$. In general, for 
$S$ a subset of $A$, we will write $S_-$ or $S^-$ for $S\cap A_-$, and, 
for $\tH$ a subgroup of $\tG$, we will write $\HH$ for $\tH\cap G$.

\medskip

Let $\psi_0$ be a character of the additive group of $F_0$, with
conductor $\po$. Then we put $\psi_F=\psi_0\circ\tr_{F/F_0}$, a
character of the additive group of $F$ with conductor $\pf$, and
$\psi_A=\psi_F\circ\tr_{A/F}$.

\medskip

An $\of$-\emph{lattice sequence\/} in $V$ is a 
function $\Lambda$ from $\BZ$ to the set of $\of$-lattices in $V$ such that
\begin{enumerate}
\item $\Lambda(k)\subset\Lambda(j)$, for $k\ge j$;
\item there exists a positive integer $e=e(\Lambda|\of)$, called the 
$\of$-period of $\Lambda$, such that $\pif\Lambda(k)=\Lambda(k+e)$, for all 
$k\in\BZ$.
\end{enumerate}

For $L$ an $\of$-lattice in $V$, we put $L^\#=\{v\in V:h(v,L)\subset\pf\}$. 
Then we call an $\of$-lattice sequence $\Lambda$ \emph{self-dual\/} if there 
exists $d\in\BZ$ such that $\Lambda(k)^\#=\Lambda(d-k)$, for all $k\in\BZ$. 
Without changing any of the objects associated to a self-dual $\of$-lattice 
sequence $\Lambda$ (except for a scale of the indices), {\bf we may (and
do) normalize all self-dual lattice sequences $\boldsymbol\Lambda$ so
that $\boldsymbol d\boldsymbol=\boldsymbol 1$ and so that all periods
are even}.
There is also a well-defined notion of the direct sum of lattice
sequences (see \cite[\S2]{BK1}). The orthogonal direct sum of self-dual lattice
sequences is itself self-dual, by the assumption $d=1$.

\medskip

Associated to an $\of$-lattice sequence $\Lambda$ in $V$, we have a
decreasing filtration $\{\fa_n(\Lambda):n\in\BZ\}$ of $A$ by
$\of$-lattices; $\fa_0$ is a hereditary $\of$-order in $A$ and $\fa_1$
is its Jacobson radical. The filtration on $A$ also gives rise to a
valuation $\nu_\Lambda$ on $A$, with $\nu_\Lambda(0)=+\infty$. 
If $\Lambda$ is self-dual, then each $\fa_n(\Lambda)$ is fixed by
$\sigma$ and $\fa_n^-=\fa_n^-(\Lambda)=\fa_n(\Lambda)\cap A_-$ gives a
filtration of $A_-$ by $\of$-lattices. Moreover,
$\nu_\Lambda$ is fixed by $\sigma$.

\medskip

Given an $\of$-lattice sequence $\Lambda$, we also put
$\til\PP=\til\PP(\Lambda)=\fa_0(\Lambda)^\times$, a compact open subgroup of
$\tG$, and $\til\PP_n=\til\PP_n(\Lambda)=1+\fa_n(\Lambda)$, for $n>0$, a
filtration of $\til\PP(\Lambda)$ by normal subgroups. 
If $\Lambda$ is self-dual, then $\til\PP$, $\til\PP_n$ are fixed by
$\sigma$ and we put $\PP^+=\PP^+(\Lambda)=\til\PP\cap G^+$, a compact
open subgroup of $G^+$, and $\PP=\PP(\L)=\PP^+\cap G$. We have a filtration
of $\PP(\Lambda)$ by normal subgroups
$\PP_n=\PP_n(\Lambda)=\til\PP_n^\Sigma=\til\PP_n\cap G$, for $n>0$.
 We also have, for $n>0$, a bijection $\fa_n^- \to 
\PP_n$ given by the Cayley map $x\mapsto C(x)=(1+\frac x2)(1-\frac
x2)^{-1}$, which is equivariant under conjugation by $\PP$.

The quotient group $\CG=\PP/\PP_1$ is (the group of rational points
of) a reductive group over the finite field $k_{F_0}$. However, it is
not, in general, connected. We denote by $\PP^\so=\PP^\so(\L)$ the
inverse image in $\PP$ of (the group of rational points
of) the connected component $\CG^\so$ of $\CG$; then $\PP^\so$ is a
parahoric subgroup of $G$.

\begin{Definition}[{\cite[\S1.5]{BK}}, {\cite[\S3.1]{BK1}}] \rm
\begin{enumerate}
\item A \emph{stratum\/} in $A$ is a 4-tuple
$[\Lambda,n,r,b]$, where $\Lambda$ is an $\of$-lattice sequence, 
$n\in\BZ$ and $r\in\BZ$ with $n\ge r\ge 0$ and $b\in\fa_{-n}(\Lambda)$.  
\item Two strata $[\Lambda,n,r,b_i]$, $i=1,2$, are called \emph{
equivalent\/} if $b_1-b_2\in\fa_{-r}(\Lambda)$.
\item A stratum $[\Lambda,n,r,b]$ is called \emph{skew\/} if
$\Lambda$ is self-dual and $b\in A_-$. 
\item A stratum $[\Lambda,n,r,b]$ is called \emph{null\/} if
$n=r$ and $b=0$.
\end{enumerate}
\label{stratumdef}
\end{Definition}

For $n\ge r\ge \frac n2>0$, an equivalence class of 
strata corresponds to a character of $\til\PP_{r+1}(\Lambda)$, by
$$
[\Lambda,n,r,b]\mapsto (\psi_b:x\mapsto \psi_A(b(x-1)),\hbox{ for
}x\in \til\PP_{r+1}),
$$
while an equivalence class of 
skew strata corresponds to a character of $\PP_{r+1}(\Lambda)$, by
$$
[\Lambda,n,r,b]\mapsto \psi_b^-=\psi_b|_{\PP_{r+1}(\Lambda)}.
$$

\begin{Definition}[{\cite[Definition 1.5.5]{BK}}, 
{\cite[Definition 5.1]{BK1}}] \rm
A stratum $[\Lambda,n,r,\beta]$ in $A$ is called \emph{simple\/} if
either it is null or the following four conditions are satisfied:
\begin{enumerate}
\item the algebra $E=F[\beta]$ is a field;
\item $\Lambda$ is an $\oe$-lattice chain (we usually write 
$\Lambda_{\oe}$ when we are thinking of it as such);
\item $\nu_{\Lambda}(\beta)=-n$;
\item $\ko<-r$.
\end{enumerate}
\label{simpledef}
\end{Definition}

Here, $\ko$ is an integer related to the action of the adjoint map
$a_\b(x)= \beta x-x \beta$ on the filtration associated to $\L$. We
refer the reader to \cite[\S1.5]{BK} or \cite[\S5]{S1} for the formal
definition and a discussion of it.

\medskip

Now let $[\Lambda,n,r,\beta]$ be a stratum in $A$ and suppose we have a
decomposition $V=\bigoplus_{i=1}^l V^i$ into $F$-subspaces. Let
$\Lambda^i$ be the 
lattice sequence in $V^i$ given by $\Lambda^i(k)=\Lambda(k)\cap V^i$
and put $\beta_i= \boI^i\beta\boI^i$, where $\boI^i$ is the projection
onto $V^i$ with kernel $\bigoplus_{j\ne i}V^j$. 

\begin{Definition}\rm\begin{enumerate}
\item We say that $V=\bigoplus_{i=1}^l V^i$ is a
\emph{splitting\/} for the stratum $[\Lambda,n,r,\beta]$ if we have
$\Lambda(k)=\bigoplus_{i=1}^l \Lambda^i(k)$, for all $k\in\BZ$, and $\beta=\sum_{i=1}^l\beta_i$. 
\item If $\CB=\{\bs v_1,...,\bs v_N\}$ is a basis for $V$ then we say that
  $\CB$ \emph{splits} the stratum $[\Lambda,n,r,\beta]$ if
  $V=\bigoplus_{i=1}^N F\bs v_i$ is a splitting for
  $[\Lambda,n,r,\beta]$.
\end{enumerate}
\end{Definition}

\begin{Definition}[{\cite[Definition 3.2]{S4}}] \rm A stratum
$[\Lambda,n,r,\beta]$ in $A$ is
called \emph{semisimple\/} if either it is null or
$\nu_\Lambda(\beta)=-n$ and there exists a splitting 
$V=\bigoplus_{i=1}^l V^i$ for the stratum such that 
\begin{enumerate}
\item for $1\le i\le l$, $[\Lambda^i,q_i,r,\beta_i]$ is a
simple or null stratum in $A^{ii}$, where $q_i=r$ if $\beta_i=0$,
$q_i=-\nu_{\Lambda^i}(\beta_i)$ otherwise; and 
\item for $1\le i,j\le l$, $i\ne j$, the stratum 
$[\Lambda^i\oplus\Lambda^j,q,r, \beta_i+\beta_j]$ is not equivalent
to a simple or null stratum, with $q=\max \{q_i,q_j\}$.
\end{enumerate}
In this case, the splitting is uniquely determined (upto ordering) by
the stratum. We will use the block notation
$A^{ij}=\Hom_F(V^j,V^i)$ and put $\CM_\b=\bigoplus_{i=1}^l A^{ii}$.
\label{semisimpledef}
\end{Definition}

Let $[\Lambda,n,r,\beta]$ be a semisimple stratum as above and put
$E=F[\b]=\bigoplus_{i=1}^l E_i$, where $E_i=F[\b_i]$. 
We will sometimes say that $\L$ is an $\oe$-lattice sequence to mean
that $\L=\bigoplus_{i=1}^l \L^i$ and each $\L^i$ is an
$\fo_{E_i}$-lattice sequence in $V^i$.

We let $B=B_\beta$ denote the $A$-centralizer of $\beta$, so that
$B=\bigoplus_{i=1}^l B_i$, where $B_i$
is the centralizer of $\beta_i$ in $A^{ii}$. 
We write $\tG_E=B^\times$, $\tG_i=\Aut_F(V^i)$ and
$\tG_{E_i}=B_i^\times=\tG_i\cap\tG_E$, so that $\tM_\b=\CM_\b^\times$ is
a Levi subgroup of $\tG$ and $\tG_E=\prod_{i=1}^l\tG_{E_i}\subset\tM_\b$. 
Each $\tG_{E_i}$ is (the group of $F_0$-points of) the restriction of
scalars to $F_0$ of a general linear group over $E_i$. 
We also write $\fb_n(\L)=\fa_n(\L)\cap B$, for $n\in\BZ$, which gives
the filtration induced on $B$ by thinking of $\L$ as an $\oe$-lattice
sequence.

\medskip

\begin{Definition}\rm
We say that a semisimple stratum $[\L,n,r,\b]$ is \emph{skew} if $\L$ is self-dual, 
$\b\in A_-$ and the associated splitting $V=\bigoplus_{i=1}^l V^i$ is
orthogonal with respect to the form $h$.
\end{Definition}

If $[\L,n,r,\b]$ is skew semisimple with splitting
$V=\bigoplus_{i=1}^l V^i$, then the strata
$[\Lambda^i,q_i,r,\beta_i]$ are all skew simple strata, for the form
$h_i$ induced by $h$ on $V^i$. The
involution on $F$ extends to each $E_i$, as the adjoint involution, and we write $E_{i,0}$ for
the subfield of fixed points; it is a subfield of index $2$ except in
the case $E_i=F=F_0$ (so that $\b_i=0$).

For $[\Lambda,n,r,\beta]$ a skew semisimple stratum, we put $G_E^+=B\cap
G^+$ and $G_E=B\cap G$; then $G_E^+$ is (the group of $F_0$-points of) a
reductive group over $F_0$ and $G_E$ is (the group of $F_0$-points of) its
connected component. We write $G^+_E=\prod_{i=1}^l G^+_{E_i}$ and
$G_E=\prod_{i=1}^l G_{E_i}$; then each $G^+_{E_i}$ is (the
group of $F_0$-points of) the restriction of scalars to $F_0$ of a unitary,
symplectic or orthogonal group over $E_{i,0}$; moreover, at
most one of the $G^+_{E_i}$ is symplectic or orthogonal -- this happens in 
the case $\beta_i=0$ and $F=F_0$. If $G^+_{E_i}$ is not orthogonal,
then we have $G_{E_i}=G_{E_i}^+$; if it is orthogonal, then $G_{E_i}$
is special orthogonal.

Note that, by~\cite[Lemmas 5.2--5.5]{BS}, for each $i$ there is a
nondegenerate $E_i/E_{i,0}$ $\ve$-hermitian form $f_i$ on $V^i$ such
that the two notions of lattice duality for $\fo_{E_i}$-lattices in $V^i$
given by $h_i$ and by $f_i$ coincide. Then $G^+_{E_i}$ is the group of
fixed points in $\tG_{E_i}$ of the involution determined by $f_i$.

\medskip

We write $\til\PP_n(\L_{\oe})=\til\PP_n(\L)\cap\tG_E=1+\fb_n(\L)$ and
$\PP_n(\L_{\oe})=\PP_n(\L)\cap G_E$, for $n> 0$. Similarly, we write 
$\til\PP(\L_{\oe})=\til\PP(\L)\cap\tG_E$ and
$\PP^+(\L_{\oe})=\PP^+(\L)\cap G_E^+$, $\PP(\L_{\oe})=\PP(\L)\cap
G_E$. Let $\PP^\so(\L_{\oe})$ denote the inverse image in
$\PP(\L_{\oe})$ of the connected component of the quotient
$\PP(\L_{\oe})/\PP_1(\L_{\oe})$;  note that this can be smaller than
$\PP^\so(\L)\cap G_E$.

\begin{Lemma}[cf.\ {\cite[Theorem 1.6.1]{BK}}] Let $[\L,n,0,\b]$ be a skew
semisimple stratum. Then
\begin{eqnarray*}
&\PP_1(\L)x\PP_1(\L) \cap G^+_E\ =\ \PP_1(\L_{\oe}) x \PP_1(\L_{\oe}),&
\hbox{ for }x\in G^+_E; \\
&\PP_1(\L)x\PP_1(\L) \cap G_E\ =\ \PP_1(\L_{\oe}) x \PP_1(\L_{\oe}),&
\hbox{ for }x\in G_E.
\end{eqnarray*}
\label{inter}
\end{Lemma}

\begin{proof} 
We begin by proving the corresponding statement for $\tG$: for
$x\in\tG_E$,
\begin{equation}\label{intertG}
\til\PP_1(\L)x\til\PP_1(\L) \cap\tG_E = 
\til\PP_1(\L_{\oe}) x \til\PP_1(\L_{\oe}).
\end{equation}
The simple case for this is given by~\cite[Theorem 1.6.1]{BK}. For the
general case, we have $\tG_E\subset\tM_\b$ and $\tM_\b$ is the fixed-point
subgroup in $\tG$ of a $2$-group of automorphisms (namely the group of
conjugations by 
$\sum_{i=1}^l \e_i\boI^i$, for $\e_i=\pm 1$). Since $p\ne 2$,
$x\in\tM_\b$ and $\til\PP_1(\L)$ is a pro-$p$ group,~\cite[Lemma 2.1]{S1}
implies that $\til\PP_1(\L)x\til\PP_1(\L) \cap\tM_\b =
(\til\PP_1(\L)\cap\tM_\b) x 
(\til\PP_1(\L)\cap\tM_\b)$. The result now follows from the simple case, by
intersecting with $\tG_E$.

Now to prove the first statement of the Lemma, we take the $\s$-fixed points 
of~\eqref{intertG} and apply~\cite[Lemma 2.1]{S1} again. The second
statement follows by intersection with $G$, since $\PP_1(\L)\subset G$.
\end{proof}


\subsection{Buildings}\label{S2.2}

We will need a few results on inclusions of the hereditary orders given
by lattice sequences, which are most easily proved using the
Bruhat-Tits affine
building. We recall briefly the description of the building in terms of
lattice sequences -- see~\cite{BL} and~\cite{BS} for more details and
proofs.

There is a (unique upto translation) $\tG$-set isomorphism between the
(extended) affine building $\CI(\tG,F_0)$ and the set of lattice
\emph{functions} on $V$, of which the lattice sequences are the subset
of ``rational points'' (i.e.\ barycentres of vertices with rational
weights). The involution $\s$ on $\tG$ induces an involution on
$\CI(\tG,F_0)$ which, for a unique choice of the isomorphism above,
coincides with the map $\L\mapsto\L^\#$, where
$\L^\#(k)=\left(\L(1-k)\right)^\#$, on lattice sequences. Moreover, we
can identify the set of fixed points with the affine building
$\CI(G^+,F_0)$ of $G^+$. Thus the affine building is identified with the
set of self-dual lattice functions on $V$. In particular, we will
identify a self-dual lattice function with the corresponding points in
$\CI(G^+,F_0)$ and $\CI(\tG,F_0)$; thus we can, for example, talk about
the line segment joining two lattice sequences.

\begin{Lemma}[cf.\ {\cite[Lemme 1.7]{Sec2}}]\label{lattcont}
Given (self-dual) $\oe$-lattice sequences $\L$, $\L'$ with $\fb_0(\L)\supset
\fb_0(\L')$, there exists a (self-dual) $\oe$-lattice sequence $\L''$
such that
$$
\fb_0(\L'')=\fb_0(\L')\qquad\hbox{and}\qquad
\fa_0(\L)\supset \fa_0(\L'').
$$
\end{Lemma}

Note that we are identifying an $\oe$-lattice sequence both with a point of $\CI(\tG_E,F_0)$ and a point of $\CI(\tG,F_0)$, which gives a continuous affine injection of buildings $\CI(\tG_E,F_0)\into\CI(\tG,F_0)$ (see~\cite{BL}). 

\begin{proof} 
We use the following facts, translated into the language of lattice
sequences from~\cite[Lemme 1.5]{Sec2}):
\begin{enumerate}
\item There is a neighbourhood of $\L$ in $\CI(\tG,F_0)$ such that,
for all $\L''$ in this neighbourhood, $\fa_0(\L)\supset\fa_0(\L'')$.
\item There is a sequence of $\oe$-lattice sequences $\L_j$ such that
$\fb_0(\L_j)=\fb_0(\L')$ and the limit of $(\L_j)$ is $\L$.
\end{enumerate}
Note that (ii) is given by~\cite[Lemme 1.5(ii)]{Sec2} only in the
simple case; in the general case, we use the facts that
$\CI(\tG_E,F_0)\simeq \prod_{i=1}^l \CI(\tG_{E_i},F_0)$ and
$\CI(\tG_{E_i},F_0)\simeq \CI(B_i^\times,E_i)$, so we can
apply the simple case in each $\CI(B_i^\times,E_i)$.

The result is now immediate: for large enough $j$, the sequence $\L_j$
of (ii) is inside the neighbourhood given by (i) so $\L''=\L_j$ has
the required properties.

For the self-dual case, note that we may assume that the sequence in
(ii) consists of self-dual lattice sequences, by replacing $\L_j$ by
the barycentre of $\L_j$ and $\L_j^\#$. The result now follows in the
same way.
\end{proof}

\begin{Corollary}\label{contmin}
Given any (self-dual) $\oe$-lattice sequence $\L$, there exists a
(self-dual) $\oe$-lattice sequence $\L^{\sm}$ such that
$\fb_0(\L^\sm)$ is a minimal (self-dual) $\oe$-order in $B$ and
$\fa_0(\L^\sm)\subset \fa_0(\L)$.
\end{Corollary}

{\bf Remark} The same result is not true for maximal
(self-dual) orders: given a (self-dual) $\oe$-lattice sequence $\L$,
there is \emph{not} necessarily a (self-dual) $\oe$-lattice sequence
$\L^{\sM}$ such that $\fb_0(\L^\sM)$ is a maximal (self-dual)
$\oe$-order in $B$ and $\fa_0(\L^\sM)\supset \fa_0(\L)$.

\begin{Lemma}\label{intermediate}
Let $\L,\L'$ be (self-dual) $\oe$-lattice sequences in $A$ with
$\fb_0(\L)=\fb_0(\L')$. Then there is a sequence
$\L=\L_0,\L_1,...,\L_t=\L'$ of (self-dual) $\oe$-lattice sequences on the line
segment $[\L,\L']$ such that, for
$1\le i\le t$:
\begin{enumerate}
\item $\fb_0(\L_i)=\fb_0(\L)$;
\item either $\fa_0(\L_i)\subset\fa_0(\L_{i-1})$ or
  $\fa_0(\L_i)\supset\fa_0(\L_{i-1})$.
\end{enumerate}
\end{Lemma}

\begin{proof}
Note that any lattice sequence $\L''$ on the line segment $[\L,\L']$ is
certainly a (self-dual) $\oe$-lattice sequence and we have
$\fb_0(\L'')=\fb_0(\L)$. On the other hand, the map $\L''\mapsto
\fa_0(\L'')$ is locally constant on this line segment, except at a
discrete (so, in particular, finite) set of points -- namely, the
points where the line segment meets some wall in the building. This
divides the line segment into a finite set of subsegments, on the
interior of which $\L''\mapsto\fa_0(\L'')$ is constant. Taking one
point from the interior of each subsegment, together with the
endpoints of the subsegments, gives the required sequence.
\end{proof}

\begin{Remark}\label{rmk:closure}\rm
Suppose $\L_1,\L_2,\L$ are $\oe$-lattice sequences such that
$\fa_0(\L_i)\supset\fa_0(\L)$, for $i=1,2$. Then, for any $\L'$ on the
line segment $[\L_1,\L_2]$, we also have $\fa_0(\L')\supset\fa_0(\L)$
-- this is clear since $\fa_0(\L')\supset\fa_0(\L)$ if and only if
$\L'$ is contained in the closure of the facet containing $\L$. 

Likewise, if $\fa_0(\L_i)\subset\fa_0(\L)$, for $i=1,2$, then, for any
$\L'$ on the line segment $[\L_1,\L_2]$, we also have
$\fa_0(\L')\subset\fa_0(\L)$. 
One can check thus using the explicit description of the building in
terms of lattice sequences as follows. We choose an apartment
containing both $\L_1$ and $\L_2$; then also $\L$ is in the same
apartment, since it is in the closure of each $\L_i$. One can then use
the explicit description of $\fa_0(\L_i)$
from~\cite[Corollary~4.5]{BL}. The details are left to the reader.
\end{Remark}

If $\L_1,\ldots,\L_t$ are (self-dual) $\oe$-lattice sequences then we write $[\L_1,\ldots,\L_t]$ for their (closed) convex hull. If this is actually a polygon then we assume that the $\L_i$ are distinct and that we have ordered the lattice sequences so that the boundary is $\cup_{i=1}^t [\L_i,\L_{i+1}]$, with the understanding that $\L_{t+1}=\L_0$.

\begin{Lemma}\label{lem:triangle}
Let $\L,\L',\L''$ be (self-dual) $\oe$-lattice sequences with $\fb_0(\L)$ and $\fb_0(\L')$ both contained in $\fb_0(\L'')$. Then there is a decomposition of the triangle $[\L,\L',\L'']$ into a union of triangles $[\L_i,\L'_i,\L''_i]$, where $\L_i,\L'_i,\L''_i$ are (self-dual) $\oe$-lattice sequences with the property that $\fa_0(\L_i)\subset\fa_(\L'_i)\subset\fa_0(\L''_i)$.
\end{Lemma}

Note that the hypotheses imply that $\L,\L',\L''$ lie in a common apartment of $\CI(\tG_E,F_0)$; since the embedding of buildings $\CI(\tG_E,F_0)\into\CI(\tG,F_0)$ respects apartments, they also lie in a common apartment of $\CI(\tG,F_0)$, so that the convex hull $[\L,\L',\L'']$ is genuinely a triangle. The same is true in the self-dual case.

\begin{proof}
The idea is essentially the same as that of Lemma~\ref{intermediate}; indeed, it already gives the result if $\L,\L',\L''$ are collinear, so we assume this is not the case.

Removing all walls of the building from $[\L,\L',\L'']$ splits it into a disjoint union of convex polygons. We split the closure $[\Ga_1,...,\Ga_t]$ of each such polygon into triangles as follows: choose a point $\Ga_0$ in the interior of the polygon and points $\Ga'_i$ on each (open) edge $(\Ga_i,\Ga_{i+1})$, for $1\le i\le t$; then the triangles are
$$
[\Ga_0,\Ga'_i,\Ga_i]\quad\hbox{and}\quad [\Ga_0,\Ga'_i,\Ga_{i+1}],\quad
\hbox{for }1\le i\le t.
$$
Taking all these triangles for all the polygons gives a decomposition as required.
\end{proof}

We note that, in the situation of Lemma~\ref{lem:triangle}, if we actually have $\fb_0(\L)=\fb_0(\L')=\fb_0(\L'')$, then also $\fb_0(\L_i)=\fb_0(\L'_i)=\fb_0(\L''_i)=\fb_0(\L)$.


\Section{Semisimple characters and Heisenberg extensions}\label{SSSH}

We will now recall results from \cite{S4} on the definition and
properties of semisimple characters, in particular their ``transfer
property''. Then
we look at the same transfer properties for the Heisenberg extensions
which are induced by these characters, as in \cite[\S5.1]{BK} and
\cite[\S2.2]{Sec2}, and deduce the existence of certain extensions to some
larger pro-$p$ groups. Finally, we characterize these extensions in
terms of their intertwining.


\subsection{Semisimple characters}

We will not repeat in detail the definitions of the objects used here
but refer the reader to~\cite[\S3]{BK} and~\cite[\S3]{S4} for more
details, in the simple and semisimple cases respectively.

Let $[\L,n,0,\b]$ be a semisimple stratum with associated splitting
$V=\bigoplus_{i=1}^l V^i$. We use the same notation as in
\S\ref{S2.1}, so $\L=\bigoplus_{i=1}^l \L^i$, $\b=\sum_{i=1}^l \b_i$,
$B=C_A(E)=\bigoplus_{i=1}^l B_i$, $\tG_E=B^\times$,
$\til\PP(\L_{\oe})=\til\PP(\L)\cap\tG_E$, etc.

In~\cite[\S3.2]{S4} we define, inductively along $\ko$, a pair of
$\of$-orders
$$
\fH(\b,\L)\subset\fJ(\b,\L)\subset\fa_0(\L).
$$
This gives us two compact open subgroups
$\tH(\b,\L)=\fH(\b,\L)^\times$ and $\tJ(\b,\L)=\fJ(\b,\L)^\times$ of
$\tG$, filtered by $\tH^{m+1}(\b,\L)=\tH(\b,\L)\cap\til\PP_{m+1}(\L)$ and  
$\tJ^{m+1}(\b,\L)=\tJ(\b,\L)\cap\til\PP_{m+1}(\L)$, for $m\ge 0$.
Moreover, we have $\tJ(\b,\L)=\til\PP(\L_{\oe})\tJ^1(\b,\L)$.

When the stratum $[\L,n,0,\b]$ is skew, $\fH$ and $\fJ$ are stable
under the involution so, as usual, we put $\JJ^+(\b,\L)=\tJ(\b,\L)\cap
G^+$ and $\JJ(\b,\L)=\tJ(\b,\L)\cap G$ etc. We also put
$\JJ^\so(\b,\L)=\PP^\so(\L_{\oe})\JJ^1(\b,\L)$, the inverse image in
$\JJ(\b,\L)$ of the connected component of $\JJ(\b,\L)/\JJ^1(\b,\L)$.

\medskip

Returning to the case of a general semisimple stratum,
in~\cite[Definition 3.13]{S4} we
define, again inductively along $\ko$, a set $\CC(\L,0,\b)$ of abelian
characters $\tth$ of $\tH^1(\b,\L)$, called \emph{semisimple characters},  
by the following two properties:
\begin{enumerate}
\item we have $\tth|_{\tH^1\cap\tM_\b}=\bigotimes_{i=1}^l \tth_i$, for some 
\emph{simple} characters $\tth_i\in\CC(\L^i,0,\b_i)$ (in the sense
of~\cite[\S3.3]{BK});
\item writing $r=-\ko$ there is, by~\cite[Proposition 3.4]{S4}, 
a semisimple stratum $[\L,n,r,\g]$ equivalent to $[\L,n,r,\b]$, with 
$\g\in\CM_\b$; then $\tth|_{\tH^{[\frac r2]+1}(\b,\L)} 
= \tth_0|_{\tH^{[\frac r2]+1}(\b,\L)} \psi_{\b-\g}$, for some semisimple 
character $\tth_0\in\CC(\L,0,\g)$.
\end{enumerate}
Moreover, by~\cite[Lemma 3.15]{S4}, $\tth$ is trivial on the unipotent 
radical $\tU$of any parabolic subgroup of $\tG$ containing $\tM_\b$. 

When the stratum is skew, this set $\CC(\L,0,\b)$ is stable under the
involution and the set $\CC_-(\L,0,\b)$ of restrictions to
$\HH^1(\b,\L)$ of semisimple characters is called the set of
\emph{skew semisimple characters}. This can also be described in
slightly different terms via the Glauberman correspondence
(see~\cite[\S3.6]{S4} and also~\S\ref{S4.glaub}). 

\medskip

Recall that, given a representation $\rho$ of a subgroup $K$ of $G$
and $g\in G^+$, the \emph{$g$-intertwining space} of $\rho$ is
$$
I_g(\rho)=I_g(\rho|K)=\Hom_{K\cap {}^gK}(\rho,{}^g\rho),
$$
where ${}^g\rho$ is the representation of ${}^gK=gKg^{-1}$ given by 
${}^g\rho(gkg^{-1})=\rho(k)$. The $G$-\emph{intertwining} of $\rho$ is
$$
I_G(\rho)=I_G(\rho|K)=\{g\in G: I_g(\rho)\ne\{0\}\},
$$
and, similarly, we have the $G^+$-intertwining $I_{G^+}(\rho)$.
We recall the following intertwining result:

\begin{Proposition}[{\cite[Proposition 3.27]{S4}}]
\label{thetaint}
Let $[\L,n,0,\b]$ be a skew semisimple stratum and
$\th\in\CC_-(\L,0,\b)$ a skew semisimple character. 
Then $\th$ is normalized by $\JJ(\b,\L)$ and
$$
I_{G^+}(\th) = \JJ^1(\b,\L) G^+_E \JJ^1(\b,\L),\qquad
I_{G}(\th) = \JJ^1(\b,\L) G_E \JJ^1(\b,\L).
$$
\end{Proposition}

Of equally great importance is the so-called ``transfer property'' of
semisimple characters:

\begin{Proposition}[{\cite[Propositions 3.26, 3.32]{S4}}] 
Let $[\L,n,0,\b]$ and $[\L',n',0,\b]$ be semisimple strata in $A$.
Then there exists a canonical bijection 
$$
\tau_{\L,\L',\b}:\CC(\L,0,\b)\to\CC(\L',0,\b)
$$ 
such that, for $\tth\in\CC(\L,0,\b)$, 
$\tth':=\tau_{\L,\L',\b}(\tth)$ is the unique semisimple character 
in $\CC(\L',0,\b)$ such that $\tG_E \cap 
I_{\tG}(\tth,\tth')\ne\emptyset$. Moreover $\tG_E\subset 
I_{\tG}(\tth,\tth')$.

If the strata are skew then the bijection $\tau_{\L,\L',\b}$ commutes
with the involution $\s$; in particular, $\tau_{\L,\L',\b}$ restricts
to a bijection $\tau_{\L,\L',\b}:\CC_-(\L,0,\b)\to\CC_-(\L',0,\b)$.
\end{Proposition}

\begin{Remark}\rm In~\cite[\S3.5]{S4}, the transfer is given in a more
general situation between $\CC(\L,m,\b)$ and $\CC(\L',m',\b)$, where
$\CC(\L,m,\b)$ denotes the set of restrictions to $\tH^{m+1}(\b,\L)$
and $m,m'$ are related by $\lfloor m/e(\L|\of)\rfloor = \lfloor
m'/e(\L'|\of)\rfloor$. This is not quite the correct condition on the
integers $m,m'$. It should be replaced as follows: put 
$$
e=\frac{e(\L|\of)}{\lcm_{1\le i\le l}\{e(E_i/F)\}}, 
$$
and likewise $e'$; then we require that 
\begin{equation}\label{correctm/e}
\lfloor m/e\rfloor = \lfloor m'/e'\rfloor. 
\end{equation}
In particular, this ensures that $\lfloor
m/e(\L^i|\fo_{E_i})\rfloor = \lfloor
m'/e(\L^{\prime i}|\fo_{E_i})\rfloor$, for each $i$, so the transfer
from the simple case can be used in each block. Notice that the
integers called $m_0,m_0'$ in~\cite[\S3.5]{S4} also satisfy this
condition since, by~\cite[\S3.1]{S4} we have
$k_0(\b,\L)/e(\L|\of)=k_0(\b,\L')/e(\L'|\of)$ so surely 
$k_0(\b,\L)/e=k_0(\b,\L')/e'$. In particular, the proof
of~\cite[Proposition 3.26]{S4} is valid under the
condition~\eqref{correctm/e}. 
\end{Remark}


\subsection{Heisenberg extensions}

Let $[\L,n,0,\b]$ be a skew semisimple stratum and
$\th\in\CC_-(\L,0,\b)$ a skew semisimple character. 

\begin{Proposition}[{\cite[Corollary 3.29, Proposition 3.31]{S4}}] There
is a unique
irreducible representation $\eta$ of $\JJ^1$ such that $\eta$
contains $\th$. Moreover, $\dim\eta=(\JJ^1:\HH^1)^{\frac 12}$ and
the intertwining of $\eta$ is given by
$$
\dim I_g(\eta)=\begin{cases} 1&\hbox{ if }g\in \JJ^1 G^+_E
\JJ^1;\\
0&\hbox{ otherwise}. \end{cases}
$$
\label{eta}
\end{Proposition}

We also have a relationship between the representations $\eta$ for
different $\oe$-lattice sequences $\L$.

\begin{Lemma}[cf.\ {\cite[Proposition 5.1.2]{BK}}] 
Let $[\L^i,n_i,0,\b]$ be skew
semisimple strata, for $i=1,2$, and let $\th_i\in\CC(\L^i,0,\b)$; let
$\eta_i$ be the representation of $\JJ^1_i=\JJ^1(\b,\L^i)$ given by
Proposition~\ref{eta}. Then 
$$
\frac{\dim\eta_1}{\dim\eta_2}=\frac{\left(\JJ^1_1:\JJ^1_2\right)} 
{\left(\PP_1(\L^1_{\oe}):\PP_1(\L^2_{\oe})\right)}.
$$
\label{etacomp}
\end{Lemma}

\begin{proof} From~\cite[Lemma 3.17]{S4}, we have exact sequences in $A_-$
$$
0\ \nomap\ \fb_1(\L^i)_-\ \nomap\ \fJ^1(\beta,\L^i)_-\ \abmap\
(\fH^1(\b,\L^i))^*_-\ \smap\ \fb_0(\L^i)_-\ \nomap\ 0,
$$
where $s:A\to B$ is a \emph{tame corestriction} relative to $E/F$ which 
commutes with the involution
(see~\cite[\S1.3]{BK} and~\cite[Lemma 4.4]{S2} for the definition).
Since the Cayley map is a bijection $\fJ^1(\b,\L^i)_-\to
\JJ^1_i$ (and similarly for $\HH$), the result follows from
Proposition~\ref{eta} exactly as in the proof of~\cite[Proposition
  5.1.2]{BK}.
\end{proof}

Now suppose we are given:
\begin{itemize}
\item skew semisimple strata $[\L,n,0,\b]$, $[\L^\sm,n_\sm,0,\b]$ 
and $[\L^\sM,n_\sM,0,\b]$ with (necessarily) the same splitting, such that
$\fa_0(\L^\sm)\subset\fa_0(\L)\subset \fa_0(\L^\sM)$; by duality we then
also have that $\fa_1(\L^\sm)\supset\fa_1(\L) \supset\fa_1(\L^\sM)$. 
\item semisimple characters $\th\in\CC_-(\L,0,\b)$,
$\th_\sm\in\CC_-(\L^\sm,0,\b)$ and $\th_\sM\in\CC_-(\L^\sM,0,\b)$, related by
the correspondences: 
$$
\th_\sm=\tau_{\L,\L^\sm,\b}(\th), \qquad
\th_\sM=\tau_{\L,\L^\sM,\b}(\th)
$$ 
\item $\eta$ (respectively $\eta_\sm$, $\eta_\sM$) the unique irreducible
representation of $\JJ^1=\JJ^1(\b,\L)$ (respectively
$\JJ^1_\sm=~\JJ^1(\b,\L^\sm)$, $\JJ^1_\sM=\JJ^1(\b,\L^\sM)$)
containing $\th$ (respectively $\th_\sm$, $\th_\sM$), given by
Proposition~\ref{eta}. 
\end{itemize}

We form the group $\JJ^1_{\sm,\sM}=\PP_1(\L^\sm_{\oe})\JJ^1_\sM$; note
that this is indeed a group,
since $\PP(\L^\sM_{\oe})$ normalizes $\JJ^1_\sM$ and contains
$\PP_1(\L^\sm_{\oe})$.

\begin{Proposition}[cf.\ {\cite[Propositions 5.1.4, 5.1.19]{BK}}, 
{\cite[Proposition 2.12]{Sec2}}]
There exists a unique
irreducible representation $\eta_{\sm,\sM}$ of $\JJ^1_{\sm,\sM}$ such that
\begin{enumerate}
\item $\eta_{\sm,\sM}|_{\JJ^1_\sM}=\eta_\sM$;
\item $\eta_{\sm,\sM}$ and $\eta_\sm$ induce equivalent irreducible
representations of $\PP_1(\L^\sm)$.
\end{enumerate}
Moreover, the intertwining of $\eta_{\sm,\sM}$ is given by
$$
\dim I_g(\eta_{\sm,\sM})=\begin{cases} 1&\hbox{ if }g\in \JJ^1_{\sm,\sM} G_E^+
\JJ^1_{\sm,\sM};\\
0&\hbox{ otherwise}. \end{cases}
$$
\label{etatilde}
\end{Proposition}

\begin{proof} Given Proposition \ref{eta}, Lemma \ref{etacomp} and
Lemma \ref{inter}, the proof of the existence and uniqueness of
$\eta_{\sm,\sM}$ is the same as that of \cite[Proposition 5.1.14]{BK} 
while the proof
of the intertwining is the same as that of \emph{op.\ cit.\/}
Proposition 5.1.19.
\end{proof}

{\bf Remark} The representation $\eta_{\sm,\sM}$ depends only on
$\fb_0(\L^\sm)$, not on the lattice sequence $\L^\sm$. For suppose
$\L^{\sm'}$ is another $\oe$-lattice sequence with
$\fb_0(\L^\sm)=\fb_0(\L^{\sm'})$ and
$\fa_0(\L^{\sm'})\subset\fa_0(\L^\sM)$, and let $\eta_{\sm',\sM}$ be
the representation of $\JJ^1_{\sm,\sM}$ given by
Proposition~\ref{etatilde} applied with $\L^{\sm'}$ in place of $\L^\sm$.

By Lemma~\ref{intermediate}, there is a sequence
$\L^\sm=\L_0,\L_1,...,\L_t=\L^{\sm'}$ of $\oe$-lattice sequences on
the line segment $[\L^\sm,\L^{\sm'}]$ such that, for $1\le i\le t$:
\begin{enumerate}
\item $\fb_0(\L_i)=\fb_0(\L^\sm)$;
\item either $\fa_0(\L_i)\subset\fa_0(\L_{i-1})$ or
  $\fa_0(\L_i)\supset\fa_0(\L_{i-1})$.
\end{enumerate}
Let $\th_i=\tau_{\L,\L^i,\b}(\th)$ and let $\eta_i$ be the unique
irreducible representation of $\JJ^1_i=\JJ^1(\b,\L_i)$ containing $\th_i$.
Since each $\L_i$ is on the line segment
$[\L^\sm,\L^{\sm'}]$, we have $\fa_0(\L_i)\subset\fa_0(\L^\sM)$ (see Remark~\ref{rmk:closure}).
Let $\eta_{i,\sM}$ denote the representation of $\JJ^1_{\sm,\sM}$ given by
Proposition~\ref{etatilde} applied with $\L_i$ in place of $\L^\sm$.

Fix $i$ and suppose $\fa_0(\L_i)\subset\fa_0(\L_{i-1})$ (the opposite
case is symmetrically similar). Applying Proposition~\ref{etatilde}
with the pair $(\L_i,\L_{i-1})$ in place of $(\L^\sm,\L^\sM)$, we see
that $\eta_i,\eta_{i-1}$ induce equivalent irreducible representations
of $\PP_1(\L_i)$. But then
$$
\Indu{\JJ^1_{\sm,\sM}}{\PP_1(\L_i)}{\eta_{i-1,\sM}} \simeq 
\Indu{\PP_1(\L_{i-1})}{\PP_1(\L_i)} 
{\Indu{\JJ^1_{\sm,\sM}}{\PP_1(\L_{i-1})}{\eta_{i-1,\sM}}} \simeq 
\Indu{\PP_1(\L_{i-1})}{\PP_1(\L_i)} 
{\Indu{\JJ^1_{i-1}}{\PP_1(\L_{i-1})}{\eta_{i-1}}} \simeq 
\Indu{\JJ^1_i}{\PP_1(\L_i)}{\eta_i}.
$$
Hence, by the uniqueness in Proposition~\ref{etatilde}, we have
$\eta_{i,\sM}=\eta_{i-1,\sM}$. Applying this with all pairs $(i,i-1)$,
we see that $\eta_{\sm',\sM}=\eta_{\sm,\sM}$, as required.

Let $\eta_{\L,\sM}$ denote the representation of
$\JJ^1_{\L,\sM}=\PP_1(\L_{\oe})\JJ^1_\sM$ given by 
applying Proposition \ref{etatilde} with $\L$ in place of $\L^\sm$. As 
in \cite[Proposition 5.1.18]{BK} (and with the same proof), we have
the following compatibility property:

\begin{Proposition}\label{etatildecompat}
With notation as above, we have
$$
\eta_{\sm,\sM}|_{\JJ^1_{\L,\sM}} = \eta_{\L,\sM}.
$$
\end{Proposition}


\subsection{Characters of Sylow subgroups}

We continue with the notation of the previous section but suppose, in
addition, that $\fb_0(\L^\sm)$ is a minimal self-dual $\oe$-order in $B$.
It will be useful, in this situation, to have another characterization of the
extension $\eta_{\sm,\sM}$ of $\eta_\sM$. For this we examine the
characters of the quotient group $\JJ^1_{\sm,\sM}/\JJ^1_\sM \cong 
\PP_1(\L^\sm_{\oe})/\PP_1(\L^\sM_{\oe})$. We will think of this as a
subgroup of $\JJ^+_\sM/\JJ^1_\sM \cong 
\PP^+(\L^\sM_{\oe})/\PP_1(\L^\sM_{\oe})$, where $\JJ^+_\sM=\JJ^+(\b,\L^\sM)$.
Since $\fb_0(\L^\sm)$ is a minimal self-dual $\oe$-order, the group
$\JJ^1_{\sm,\sM}$ is in fact a Sylow pro-$p$ subgroup of $\JJ^+_\sM$.

We identify
$\PP^+(\L^\sM_{\oe})/\PP_1(\L^\sM_{\oe})$ with the product of (the
points of) reductive groups $\CG_i$ over $k_{E_{i,0}}$, which are not
necessarily connected; the
possibilities here are as follows:
\begin{itemize}
\item if $E_i/E_{i,0}$ is quadratic unramified then $\CG_i$ is the
  product of (the restriction of scalars to $k_{E_{i,0}}$ of) some general
  linear groups over $k_{E_i}$ and at most two
  unitary groups over $k_{E_{i,0}}$;
\item if $E_i/E_{i,0}$ is quadratic ramified then $\CG_i$ is the
  product of some general linear groups, 
  at most one symplectic group, and at most one orthogonal group, 
  over $k_{E_i}=k_{E_{i,0}}$.
\item if $E_i=F=F_0$ then $\CG_i$ is a product of some general
  linear groups and at most two symplectic
  (if $\ve=-1$) or orthogonal (if $\ve=+1$) groups over $\kf$.
\end{itemize}
The image of $\PP_1(\L^\sm_{\oe})$ identifies with $\CU$, the product of
the unipotent radicals $\CU_i$ of Borel subgroups of $\CG_i$.

We also remark, for future reference, that $\PP(\L^\sM_{\oe})$ is a
\emph{maximal} compact subgroup of $G_E$ if and only if there are no factors
isomorphic to a general linear group in any of the three cases
above. 

\begin{Lemma}\label{charext}
Let $\phi$ be a character of $\JJ^1_{\sm,\sM}/\JJ^1_\sM \cong \CU$ 
which is intertwined by all of $\JJ^+_\sM$. Then $\phi$ extends to a
character of $\JJ^+_\sM$. 
\end{Lemma}

\begin{proof}
We think of $\phi=\otimes\phi_i$ as a character of
$\CU=\prod_{i=1}^l\CU_i$, where $\phi_i$ is a character of $\CU_i$
intertwined by all of $\CG_i$.

We claim that the characters $\phi_i$ extend to characters, also
denoted $\phi_i$, of $\CG_i$; indeed, in all but one case, $\phi_i$
will be trivial. Inflating $\phi=\otimes \phi_i$ to
$\JJ_{\sM}$ gives the required extension.

\medskip

We fix $i$ and now omit it from our notation. Let $\CT$ be a maximal
torus in $\CG$, and $\CB$ a Borel subgroup with unipotent radical $\CU$
containing $\CT$. Let $\Si$ be the root system of $\CG$ relative to
$\CT$, and let $\Si^+$ be the positive roots, and $\D$ the simple roots,
determined by $\CB$. For $\a\in\Si$, write $\CU_\a$ for the
corresponding root subgroup, so $\CU=\prod_{\a\in\Si^+}\CU_\a$. By
\cite[page 129]{DM}, the derived subgroup of $\CU$ is
$\prod_{\a\in\Si^+\setminus\D}\CU_\a$ (note that the only
(quasi-simple) exceptions to
this are groups of type $B_2$ or $F_4$ over the field of $2$ elements
and $G_2$ over the field of $3$ elements, none of which can occur in
our situation). Hence $\phi$ must take the form
$$
\phi=\prod_{\a\in\D}\phi_\a, \qquad\hbox{for }\phi_\a\hbox{ a
character of }\CU_\a.
$$
Now write $\D=\cup\D_j$, with each $\D_j$ irreducible.
If
$\card\D_j\ge 2$ then, for all $\a\in\D_j$, there is a Weyl-group
element $w$ such that $w(\a)\in\Si^+\setminus\D$. If $n$ in the
$\CG$-normalizer of $\CT$ reduces
to $w$ then, since $n$ intertwines $\phi$, we have
$$
\phi(u_\a)=\phi(nu_\a n^{-1})=1, \qquad\hbox{for }u_\a\in \CU_\a,
$$
since $nu_\a n^{-1}\in\CU_{w(\a)}\subset\ker\phi$. Hence $\phi_\a$ is
trivial. In particular, if $\card\D_j\ge 2$ for all $j$, then $\phi$ is
trivial and extends to the trivial representation of $\CG$.

This leaves the cases when $\card\D_j=1$, for some $j$; in this case
the corresponding factor of $\CG$ is one of: 
\begin{enumerate}
\item a general linear group of dimension $2$, that is $GL_2(\ke)$;
\item a symplectic group of dimension $2$, that is $SL_2(\ke)$;
\item a unitary group of dimension $2$ or $3$;
\item an orthogonal group of dimension $3$ or $4$.
\end{enumerate} 
These cases need to be treated individually. We omit the subscript $j$.

\medskip

(i) The proof in this case is contained in the proof
of~\cite[Proposition 5.2.4]{BK} (see bottom of page 168). Indeed,
since $p\ne 2$, we have $|\ke|>2$ so the proof there shows that the
character $\phi$ is in fact trivial and extends to the trivial
character of $\CG$.

\medskip

(ii) We write the elements of $\CG\simeq SL_2(\ke)$ as matrices and
$\CU$ as the group of upper triangular unipotent matrices. The
character $\phi$ is given by
$$
\phi\begin{pmatrix} 1&n\\ 0&1\end{pmatrix} = \chi(n),\quad\hbox{ for
}n\in \ke,
$$
where $\chi$ is an additive character of $\ke$. Since $\phi$ is
intertwined by all of $\CG$, it is normalized by $\diag(a,a^{-1})$,
for $a\in \ke^\times$. Hence
$$
\chi\left((a^2-1)n\right) = 1,\quad\hbox{ for all }a\in \ke^\times,
n\in \ke. 
$$
Providing $|\ke|>3$, there is an $a\in \ke^\times$ such that $a^2\ne
1$ and hence $\chi$ is trivial and $\phi$ extends to the trivial
representation of $\CG$. 

If $|\ke|=3$ then the two non-trivial characters of $\CU$ are also
intertwined by all of $\CG=SL_2(\BF_3)$. However, in this case $\CU$
is a quotient of $\CG$ (by its $2$-Sylow subgroup) 
so $\phi$ extends in any case. (In the case of the
non-trivial characters of $\CU$, it extends to a $1$-dimensional
cuspidal representation of $\CG$.)

\medskip

(iii) The proof here is almost identical, though there is no exceptional
case as in (ii). In the $2$-dimensional case we conjugate by
$\diag(a,\ov a^{-1})$ and use the facts that the norm map $\ke\to
k_{E,0}$ is surjective and $|k_{E,0}|>2$; in the $3$-dimensional case,
we conjugate by $\diag(a,1,\ov a^{-1})$.

\medskip

(iv) Again, the proof here is the same, with no exceptional case. In
the $3$-dimensional case we conjugate by $\diag(a,1,a^{-1})$ and in
the $4$-dimensional case by $\diag(a,1,1,a^{-1})$.
\end{proof}

\begin{Lemma}\label{inttriv}
Let $\phi$ be a character of $\JJ^1_{\sm,\sM}/\JJ^1_\sM \cong \CU$ 
which is intertwined by all of $G_E$. Then $\phi$ is trivial.
\end{Lemma}

\begin{proof} 
Using the notation of the proof of Lemma~\ref{charext}, we have
$\phi=\otimes\phi_i$, and each character $\phi_i$ is intertwined by
$G_{E_i}$. It is enough to show that each $\phi_i$ is trivial so it is
enough to check the simple case. Moreover, we have seen in the proof of
Lemma~\ref{charext} that $\phi$ is trivial except in the case where
$\CG$ has a factor isomorphic to $SL(2,\BF_3)$. Hence we may restrict
to the case (in the notation of the proof of Lemma~\ref{charext}) when
$|\D|=1$ and $\CG\simeq SL_2(\BF_3)$.

We choose a self-dual $E$-basis which splits the
lattice sequence $\L^\sm_{\oe}$ (so also $\L^\sM_{\oe}$). Then, after
scaling and permuting the basis if necessary, this identifies 
\begin{eqnarray*}
\fb_0(\L^\sm_{\oe}) 
&\hbox{ with }&
\begin{pmatrix} \oe & \oe \\ \pe & \oe \end{pmatrix};
\qquad\hbox{ and } \\
\fb_0(\L^\sM_{\oe})
&\hbox{ with }&
\begin{pmatrix} \oe & \oe \\ \oe & \oe \end{pmatrix}
\quad\hbox{ or }\quad
\begin{pmatrix} \oe & \pe^{-1} \\ \pe & \oe \end{pmatrix}.
\end{eqnarray*}
(The second case occurs when $E/E_0$ is ramified and the form $h$ is
hermitian.) We consider 
$$
x\ =\ \begin{pmatrix} 1&1\\0&1\end{pmatrix} 
\quad\hbox{ or }\quad
\begin{pmatrix} 1&0\\ \pie&1\end{pmatrix}
$$ 
in the two cases, so that
$x\PP_1(\L^\sM_{\oe})$ generates $\CU$. We put
$g=\diag(\pie,\pie^{-1})$ or $\diag(\pie^{-1},\pie)$ in the two
cases. Then $gxg^{-1}\in\PP_1(\L^\sM_{\oe})$. Since $g\in G_E$
intertwines $\phi$ and $\phi$ is trivial on $\PP_1(\L^\sM_{\oe})$, we
see that
$$
\phi(x)=\phi(gxg^{-1}) = 1.
$$
Since $x\PP_1(\L^\sM_{\oe})$ generates $\CU$, the character $\phi$ is
indeed trivial.
\end{proof}

\begin{Corollary}\label{uniqueetatilde}
If $\fb_0(\L^\sm)$ is a minimal self-dual $\oe$-order contained in $\fb_0(\L^\sM)$, 
then the representation $\eta_{\sm,\sM}$ is the unique extension of
$\eta_\sM$ to $\JJ^1_{\sm,\sM}$ which is intertwined by all of $G_E$.
\end{Corollary}

\begin{proof}
Suppose $\eta'$ is another such extension. Then
$\eta'=\eta_{\sm,\sM}\otimes\phi$, for some character $\phi$ of
$\JJ^1_{\sm,\sM}/\JJ^1_\sM$ which is intertwined by all
of $G_E$. By Lemma~\ref{inttriv}, such a character is trivial so
$\eta'=\eta_{\sm,\sM}$. 
\end{proof}


\subsection{Analogous results for $\tG$}

The results of the previous two sections all have analogues for
$\tG$. The proofs are almost identical, though a little
simpler. Indeed, since so far we have only been looking at
representations of pro-$p$ subgroups, we could have obtained the
results above by first proving the analogous results for $\tG$ and
then using the Glauberman transfer, as in previous work
(see~\cite{S2},~\cite{S4}). We state those results we will need for
$\tG$ here, leaving the modifications of the proofs of the previous
sections as an exercise for the reader.

\begin{Proposition}\label{etatildetG}
Let $[\L^\sm,n_\sm,0,\b]$ and $[\L^\sM,n_\sM,0,\b]$ be semisimple
strata with $\fa_0(\L^\sm)\subset\fa_0(\L^\sM)$. Let
$\tth_\sm\in\CC(\L^\sm,0,\b)$ be a semisimple character and
$\tth_\sM=\tau_{\L^\sm,\L^\sM,\b}(\tth_\sm)$. Then there is a unique
irreducible representation $\teta_\sm$ (respectively $\teta_\sM$) of
$\tJ^1(\b,\L^\sm)$ (respectively $\tJ^1(\b,\L^\sM)$) which contains
$\tth_\sm$ (respectively $\tth_\sM$).

Put $\tJ^1_{\sm,\sM}=\til\PP_1(\L^\sm_{\oe})\tJ^1(\b,\L^\sM)$. There
is a unique irreducible representation $\teta_{\sm,\sM}$ of
$\tJ^1_{\sm,\sM}$ such that
\begin{enumerate}
\item $\teta_{\sm,\sM}|_{\tJ^1(\b,\L^\sM)}=\teta_\sM$;
\item $\teta_{\sm,\sM}$ and $\teta_\sm$ induce equivalent irreducible
representations of $\til\PP_1(\L^\sm)$.
\end{enumerate}
Moreover, if $\fb_0(\L^\sm)$ is a minimal $\oe$-order in $B$ then
$\teta_{\sm,\sM}$ is the unique extension of $\teta_\sM$ to
$\tJ^1_{\sm,\sM}$ which is intertwined by $\tG_E$. Finally, the intertwining of
$\teta_{\sm,\sM}$ is given by
$$
\dim I_g(\teta_{\sm,\sM})=\begin{cases} 1&\hbox{ if }g\in \tJ^1_{\sm,\sM} \tG_E
\tJ^1_{\sm,\sM};\\
0&\hbox{ otherwise}. \end{cases}
$$
\end{Proposition}


\Section{Beta extensions}\label{S3}

This section is devoted to the study of \emph{$\b$-extensions}, which
are certain rather special extensions of the Heisenberg
representations $\eta$ to $\JJ^+(\b,\L)$. In the case of $GL_N(F)$ (or
$GL_m(D)$ for $D$ an $F$-central  division algebra) these extensions
are characterized by their intertwining, which should contain the
whole centralizer of $\b$ (see~\cite[\S5.2]{BK},~\cite[\S2.4]{Sec2}
respectively). In our situation we have been unable to do this, but
instead give a different characterization (which, from the results
of~\cite[\S5.2]{BK}, is also valid for $GL_N(F)$) in terms of the
representation $\eta_{\sm,\sM}$ of the previous section and by
transfer properties.


\subsection{Beta-extensions in the maximal case}\label{S3.1}

Let $[\L^\sM,n,0,\b]$ be a skew semisimple stratum such that
$\fb_0(\L^\sM)$ is a maximal self-dual $\oe$-order in $B$. Let
$\th_\sM\in\CC_-(\L^\sM,0,\b)$ be a skew semisimple character and let
$\eta_\sM$ be the unique irreducible representation of
$\JJ^1_\sM=\JJ^1(\b,\L^\sM)$ containing $\th_\sM$. By
Corollary~\ref{contmin}, there exists a self-dual $\oe$-lattice
sequence $\L^\sm$ in $V$ such that $\fb_0(\L^\sm)$ is a minimal
self-dual $\oe$-order in $B$ and $\fa_0(\L^\sm)\subset\fa_0(\L^\sM)$. 
Let $\eta_{\sm,\sM}$ be the representation of
$\JJ^1_{\sm,\sM}=\PP_1(\L^\sm_{\oe})\JJ^1_\sM$ given by
Proposition~\ref{etatilde}. Recall that this depends only on the group
$\JJ^1_{\sm,\sM}$, not on the choice of $\L^\sm$.

Following the ideas of~\cite[\S5.2]{BK}, we show that there is a
representation of $\JJ^+_\sM=\JJ^+_\sM(\b,\L)$ which 
extends~$\eta_{\sm,\sM}$.

\begin{Theorem}[cf.\ {\cite[Proposition 5.2.4]{BK}}] With notation as
above, there exists a representation $\k_\sM$ of
$\JJ^+_\sM$ which extends $\eta_{\sm,\sM}$. Moreover, if
$\k_\sM'$ is another such representation then
$\k_\sM'=\k_\sM\otimes\chi$, for some character $\chi$ of
$\PP^+(\Lambda^\sM_{\oe})/\PP_1(\Lambda^\sM_{\oe})$ which is trivial on
the subgroup generated by all its unipotent subgroups.
\label{betamax}
\end{Theorem}

We call such a representation $\k_\sM$ a \emph{$\beta$-extension} of
$\eta_\sM$. Note that the definition of $\k_\sM$ is independent of the
choice of $\L^\sm$, since $\PP^+(\L^\sM_{\oe})$ acts transitively by
conjugation on the
minimal self-dual $\oe$-orders contained in $\fb_0(\L^\sM)$ and the
representation $\eta_{\sm,\sM}$ depends only on $\fb_0(\L^\sm)$, not
on the lattice sequence $\L^\sm$.

\begin{proof} The proof is essentially the same as that 
of~\cite[Proposition 5.2.4]{BK} (see also~\cite[Lemme 6.4]{Bla}) so
we will omit some of the details and refer the reader to \emph{loc.\ cit.\/}. 

Since $\JJ^+_\sM$ normalizes $\eta_\sM$, we can extend $\eta_\sM$ to a
projective representation of $\JJ^+_\sM$. This gives rise to a $2$-cocycle
$\boldsymbol\a$ of $\JJ_\sM$ in $\BCx$. Since the dimension of $\eta_\sM$ is
a power of $p$, the order of the cohomology class of this cocycle is a
power of $p$. But $\eta_\sM$ admits an extension to $\JJ^1_{\sm,\sM}$ so the
restriction of the cocycle to $\JJ^1_{\sm,\sM}$ is cohomologous to
zero. However, $\JJ^1_{\sm,\sM}$ is a Sylow pro-$p$ subgroup of
$\JJ^+_\sM$ so the original cohomology class of $\boldsymbol\a$ is
trivial and hence $\eta_\sM$ extends to a linear representation $\l$
of $\JJ^+_\sM$.

Now we compare $\l|_{\JJ^1_{\sm,\sM}}$ to $\eta_{\sm,\sM}$; it is of the form
$\eta_{\sm,\sM}\otimes \phi$, for some abelian character $\phi$ of
$\PP_1(\L^\sm_{\oe})/\PP_1(\L^\sM_{\oe})$ which is intertwined by all of
$\PP^+(\L^\sM_{\oe})$. By Lemma~\ref{charext}, any such character
extends to a character, also denoted $\phi$, of
$\PP^+(\L^\sM_{\oe})$. Then $\k_{\sM}:=\l\otimes{\phi}^{-1}$ extends
$\eta_{\sm,\sM}$, as required. Note also that the final assertion of
the Theorem is clear.
\end{proof}

\begin{Remark}\rm
The conclusions of Theorem~\ref{betamax} are also true without the assumption that $\fb_0(\L^\sM)$ is \emph{maximal}, as is the case for $\tG$ (see~\cite[Remark,~p.170]{BK}). However, as for $\tG$, this gives many more extensions than are useful for our purposes -- instead, we will define $\b$-extensions in the non-maximal case by a compatibility property with the maximal case.
\end{Remark}


\subsection{Beta-extensions in the general case}\label{S3.3}

Let $[\L,n,0,\b]$ be any skew semisimple stratum, let $\th\in\CC_-(\L,0,\b)$ 
be a skew semisimple character and let $\eta$ be the representation of
$\JJ^1=\JJ^1(\b,\L)$ given by Proposition~\ref{eta}. 
We are going to define $\b$-extensions of $\eta$ in the case where 
$\fb_0(\L)$ is not necessarily a maximal self-dual $\oe$-order, by 
transfer from the maximal case of Theorem~\ref{betamax}. In order to
obtain compatibility bijections, we will see that we do not want
to extend our representations all the way to the group $\JJ^+(\b,\L)$,
but rather to the group   
$\JJ^\so=\JJ^\so(\b,\L)=\PP^\so(\L_{\oe})\JJ^1$, the inverse image in
$\JJ(\b,\L)$ of the connected component of
$\JJ(\b,\L)/\JJ^1(\b,\L)\cong
\PP(\L_{\oe})/\PP_1(\L_{\oe})$. Nevertheless, we initially define
$\b$-extensions to the group $\JJ^+(\b,\L)$ and then restrict back to
$\JJ^\so(\b,\L)$; this is because we need to be sure that our $\b$-extensions to $\JJ^\so$ extend to $\JJ^+(\b,\L)$.

\begin{Lemma}[cf.\ {\cite[Proposition 5.2.5]{BK}, \cite[Lemme 2.23, 
Proposition 2.9]{Sec2}}]
\label{exttransfer}
Using the notation above, 
let $[\L^\sm,n_\sm,0,\b]$ be a skew semisimple stratum such that
$\fb_0(\L)\supset\fb_0(\L^\sm)$. Let
$\th_\sm=\tau_{\L,\L^\sm,\b}(\th)$ and let $\eta_\sm$ be the unique
irreducible representation
of $\JJ^1_\sm=\JJ^1(\b,\L^\sm)$ which contains $\th_\sm$. Then
there is a canonical bijection between the set of extensions $\k_\sm$ of
$\eta_\sm$ to $\JJ^+_\sm=\JJ^+(\b,\L^\sm)$ and the set of extensions
$\k$ of $\eta$ to $\JJ^+_{\sm,\L}=\PP^+(\L^\sm_{\oe})\JJ^1$. 
If $\fa_0(\L)\supset\fa_0(\L^\sm)$ then it is given as follows: 
given $\k_\sm$ (respectively $\k$) there is a unique $\k$
(respectively $\k_\sm$) such that $\k_\sm|_{\JJ^+_\sm}$ and
$\k|_{\JJ^+_{\sm,\L}}$ induce equivalent irreducible representations of
$\PP^+(\L^\sm_{\oe})\PP_1(\L^\sm)$.

Moreover, if $g\in G_E^+$ then $g$ intertwines $\k_\sm|_{\JJ^+_\sm}$ if
and only if it intertwines the corresponding $\k|_{\JJ^+_{\sm,\L}}$.
\end{Lemma} 

Note that, in the special case of Lemma~\ref{exttransfer} when
$\fb_0(\L)=\fb_0(\L^\sm)$, we have $\JJ^+_{\sm,\L}=\JJ^+(\b,\L)$. 

\begin{proof}
\underline{\emph{Case 1}} We assume first that $\fa_0(\L)\supset\fa_0(\L^\sm)$.
The proof of the first statement is identical to that
of~\cite[Proposition 5.2.5]{BK}. Given the simple intersection
property Lemma~\ref{inter}, the proof of the second statement is
identical to that of \cite[Proposition 2.9]{Sec2}.

For the general case, we first need a compatibility condition:

\begin{Lemma}[cf.\ {\cite[(5.2.14)]{BK}}]\label{lem:firstcompat}
Suppose $[\L,n,0,\b]$,
$[\L',n',0,\b]$ and $[\L^\sM,n_\sM,0,\b]$ are semisimple strata such
that $\fa_0(\L)\subset\fa_0(\L')\subset\fa_0(\L^\sM)$. Let
$\th\in\CC_-(\L,0,\b)$ and $\th'=\tau_{\L,\L',\b}$,
$\th_\sM=\tau_{\L,\L^\sM,\b}$, and let $\eta$ (respectively $\eta'$,
$\eta_\sM$) be the irreducible representation of $\JJ^1(\b,\L)$
(respectively, $\JJ^1(\b,\L')$, $\JJ^1(\b,\L^\sM)$) given by
Proposition~\ref{eta}. Put $\JJ^+_{\L,\sM}=\PP^+(\L_{\oe})\JJ^1(\b,\L^\sM)$,
$\JJ^+_{\L',\sM}=\PP^+(\L'_{\oe})\JJ^1(\b,\L^\sM)$ and 
$\JJ^+_{\L,\L'}=\PP^+(\L_{\oe})\JJ^1(\b,\L')$ so
$\JJ^+_{\L,\sM}\subset\JJ^+_{\L',\sM}$. Then the following diagram
commutes, where the maps are given by (Case 1 of) Lemma~\ref{exttransfer}, or by restriction: 
$$
\xymatrix{
\left\{\hbox{extensions of }\eta_\sM\hbox{ to }\JJ^+_{\L',\sM}\right\} 
\ar[r]^{\res}\ar[d] & 
\left\{\hbox{extensions of }\eta_\sM\hbox{ to }\JJ^+_{\L,\sM}\right\} 
\ar[dr] &  \\
\left\{\hbox{extensions of }\eta'\hbox{ to }\JJ^+(\b,\L')\right\} 
\ar[r]^{\res} & 
\left\{\hbox{extensions of }\eta'\hbox{ to }\JJ^+_{\L,\L'}\right\} 
\ar[r] & 
\left\{\hbox{extensions of }\eta\hbox{ to }\JJ^+(\b,\L)\right\} 
}
$$
\end{Lemma}

\begin{proof}
Suppose $\k_\sM$ is an extension of $\eta_\sM$ to $\JJ^+_{\L',\sM}$,
let $\k'$ be the corresponding extension of $\eta'$ to $\JJ^+(\b,\L')$,
and let $\k$ be the extension of $\eta$ to $\JJ^+(\b,\L)$ which
corresponds to $\k'|_{\JJ^+_{\L,\L'}}$. We form the group
$\PP^+(\L_{\oe})\PP_1(\L')$ and note that we have
$$
\Indu{\JJ^+_{\L,\L'}}{\PP^+(\L_{\oe})\PP_1(\L')}
{\Resu{\JJ^+(\b,\L')}{\JJ^+_{\L,\L'}}{\k'}}\ \simeq\ 
\Resu{\PP^+(\L_{\oe})\PP_1(\L)}{\PP^+(\L_{\oe})\PP_1(\L')}
{\Indu{\JJ^+(\b,\L')}{\PP^+(\L_{\oe})\PP_1(\L)}{\k'}},
$$
since there is only one double coset in the Mackey restriction formula
here. Similarly, we have
\begin{eqnarray*}
\Indu{\JJ^+_{\L,\sM}}{\PP^+(\L_{\oe})\PP_1(\L')}
{\Resu{\JJ^+_{\L',\sM}}{\JJ^+_{\L,\sM}}{\k_\sM}} &\simeq&
\Resu{\PP^+(\L'_{\oe})\PP_1(\L')}{\PP^+(\L_{\oe})\PP_1(\L')}
{\Indu{\JJ^+_{\L',\sM}}{\PP^+(\L'_{\oe})\PP_1(\L')}{\k_{\sM}}} \\ &\simeq&
\Resu{\PP^+(\L'_{\oe})\PP_1(\L')}{\PP^+(\L_{\oe})\PP_1(\L')}
{\Indu{\JJ^+(\b,\L')}{\PP^+(\L'_{\oe})\PP_1(\L')}{\k'}},
\end{eqnarray*}
where the second equivalence is from the compatibility of $\k'$ and
$\k_\sM$. Hence 
\begin{eqnarray*}
\Indu{\JJ^+_{\L,\sM}}{\PP^+(\L_{\oe})\PP_1(\L)}
{\Resu{\JJ^+_{\L',\sM}}{\JJ^+_{\L,\sM}}{\k_\sM}} &\simeq&
\Indu{\PP^+(\L_{\oe})\PP_1(\L')}{\PP^+(\L_{\oe})\PP_1(\L)}
{\Resu{\PP^+(\L'_{\oe})\PP_1(\L')}{\PP^+(\L_{\oe})\PP_1(\L')}
{\Indu{\JJ^+(\b,\L')}{\PP^+(\L'_{\oe})\PP_1(\L')}{\k'}}} \\ &\simeq&
\Indu{\JJ^+_{\L,\L'}}{\PP^+(\L_{\oe})\PP_1(\L)}
{\Resu{\JJ^+(\b,\L')}{\JJ^+_{\L,\L'}}{\k'}} \simeq
\Indu{\JJ^+(\b,\L)}{\PP^+(\L_{\oe})\PP_1(\L)}{\k}. 
\end{eqnarray*}
Hence $\k$ does indeed correspond to $\k_\sM|_{\JJ^+_{\L,\sM}}$ via
(Case 1 of) Lemma~\ref{exttransfer}.
\end{proof}

\underline{\emph{Case 2}} We return to the proof of Lemma~\ref{exttransfer} and consider now the case where $\fb_0(\L)=\fb_0(\L^\sm)$. Let
$\L=\L_0,\L_1,...,\L_t=\L^\sm$ be a sequence of self-dual $\oe$-lattice sequences given
by Lemma~\ref{intermediate}, so that $\fb_0(\L_i)=\fb_0(\L)$ and either
$\fa_0(\L_i)\subset\fa_0(\L_{i-1})$ or
$\fa_0(\L_i)\supset\fa_0(\L_{i-1})$, for $1\le i\le t$.
For each $i$, let $\th_i=\tau_{\L,\L_i,\b}(\th)$
and let $\eta_i$ be the representation of 
$\JJ^1(\b,\L_i)$ given by Proposition~\ref{eta}.
Then, using the first case repeatedly (with $\{\L_i,\L_{i+1}\}$ in place of $\{\L,\L^\sm\}$) and composing the bijections it gives, we get the required bijection between the set of extensions $\k_\sm$ of
$\eta_\sm$ to $\JJ^+_\sm$ and the set of extensions
$\k$ of $\eta$ to $\JJ^+_{\sm,\L}=\JJ^+(\b,\L)$. The assertion concerning intertwining also follows by repeated use of the first case.

To prove that this bijection is independent of the choices of the $\L_i$, we must use Lemma~\ref{lem:firstcompat}. So suppose that $\L=\L'_0,\L'_1,...,\L'_{t'}=\L^\sm$ is another sequence of self-dual $\oe$-lattice sequences such that $\fb_0(\L_i)=\fb_0(\L)$ and either
$\fa_0(\L'_i)\subset\fa_0(\L'_{i-1})$ or
$\fa_0(\L'_i)\supset\fa_0(\L'_{i-1})$, for $1\le i\le t'$.

We consider the triangles
$$
[\L,\L_i,\L_{i+1}],\hbox{ for }1\le i\le t-1,\quad\hbox{and}\quad
[\L'_{i-1},\L'_i,\L^\sm],\hbox{ for }1\le i\le t'-1.
$$
Using Lemma~\ref{lem:triangle}, each can be decomposed into smaller triangles $[\Ga_j,\Ga'_j,\Ga''_j]$ such that $\fa_0(\Ga_j)\subset\fa_0(\Ga'_j)\subset\fa_0(\Ga''_j)$. We write $\th_{\Ga_j},\th_{\Ga'_j},\th_{\Ga''_j}$ for the transfers of $\th$, and $\eta_{\Ga_j},\eta_{\Ga'_j},\eta_{\Ga''_j}$ for the representations given by Proposition~\ref{eta}. Then Lemma~\ref{lem:firstcompat} gives us commutative diagrams
$$
\xymatrix{
\left\{\hbox{extensions of }\eta_{\Ga''_j}\hbox{ to }\JJ^+(\b,\Ga''_j)\right\} 
\ar[d]\ar[dr] &  \\
\left\{\hbox{extensions of }\eta_{\Ga'_j}\hbox{ to }\JJ^+(\b,\Ga'_j)\right\} 
\ar[r] & 
\left\{\hbox{extensions of }\eta_{\Ga_j}\hbox{ to }\JJ^+(\b,\Ga_j)\right\} 
}
$$
Putting together all these commutative diagrams, we see that the bijection
$$
\left\{\hbox{extensions of }\eta\hbox{ to }\JJ^+(\b,\L)\right\} 
\ \longleftrightarrow\  
\left\{\hbox{extensions of }\eta_\sm\hbox{ to }\JJ^+_\sm\right\} 
$$
is indeed independent of the choices.

\underline{\emph{Case 3}} Finally, consider the general case $\fb_0(\L)\supset\fb_0(\L^\sm)$. By Lemma~\ref{lattcont}, there is a self-dual $\oe$-lattice sequence $\L'$ such that $\fb_0(\L')=\fb_0(\L^\sm)$ and $\fa_0(\L)\supset \fa_0(\L')$. Put $\th'=\tau_{\L,\L',\b}$ and let $\eta'$ be the irreducible representation of $\JJ^1(\b,\L')$ given by Proposition~\ref{eta}.

Applying Case 2 (with $\L'$ in place of $\L$), we get a canonical bijection between the set of extensions $\k_\sm$ of $\eta_\sm$ to $\JJ^+_\sm$ and the set of extensions $\k'$ of $\eta'$ to $\JJ^+(\b,\L')$. Now applying Case 1 (with $\L'$ in place of $\L^\sm$), we get a canonical bijection between the set of extensions $\k'$ of
$\eta'$ to $\JJ^+(\b,\L')$ and the set of extensions
$\k$ of $\eta$ to $\JJ^+_{\sm,\L}=\PP^+(\L'_{\oe})\JJ^1$. The required bijection is the composition of these two and the assertion on intertwining is immediate. 

It remains to check that the bijection is independent of the choice of $\L'$. Suppose $\L''$ is another $\oe$-lattice sequence such that $\fb_0(\L'')=\fb_0(\L^\sm)$ and $\fa_0(\L)\supset \fa_0(\L'')$, put $\th''=\tau_{\L,\L'',\b}$ and let $\eta''$ be the irreducible representation of $\JJ^1(\b,\L'')$ given by Proposition~\ref{eta}. 
Then we have a diagram in which the arrows are given by Cases 1 and 2 of Lemma~\ref{exttransfer}
$$
\xymatrix{
\left\{\hbox{extensions of }\eta\hbox{ to }\JJ^+_{\sm,\L}\right\} 
\ar[r]\ar[d] & 
\left\{\hbox{extensions of }\eta''\hbox{ to }\JJ^+(\b,\L'')\right\} 
\ar[d]\ar[dl] \\
\left\{\hbox{extensions of }\eta'\hbox{ to }\JJ^+(\b,\L')\right\} 
\ar[r] & 
\left\{\hbox{extensions of }\eta_\sm\hbox{ to }\JJ^+_\sm\right\} 
}
$$
The lower right triangle commutes by Case 2, while we can see that the upper left triangle commutes by splitting $[\L',\L'',\L]$ into triangles (using Lemma~\ref{lem:triangle}) and arguing as in Case 2. Hence the whole diagram commutes, as required.
\end{proof}

\medskip

We can now define $\b$-extensions for an
arbitrary skew semisimple stratum $[\L,n,0,\b]$ and character
$\th\in\CC_-(\L,0,\b)$. 

\begin{Definition}\rm\label{kappadef}
Let $[\L^\sM,n_\sM,0,\b]$ be a skew semisimple stratum with
$\fb_0(\L^\sM)$ a maximal self-dual $\oe$-order in $B$ and 
$\fb_0(\L^\sM)\supset\fb_0(\L)$. Let $\th_\sM=\tau_{\L,\L^\sM,\b}$,
let $\eta_\sM$ be the unique irreducible representation of
$\JJ^1_\sM=\JJ^1(\b,\L^\sM)$ which contains $\th_\sM$, and let $\k_\sM$ be a 
$\b$-extension of $\eta_\sM$ to $\JJ^+_\sM=\JJ^+(\b,\L^\sM)$, given by 
Theorem~\ref{betamax}.
Then the \emph{$\b$-extension of $\eta$ to $\JJ^+$ 
relative to $\L^\sM$, compatible with $\k_\sM$} is the unique representation $\k$ of $\JJ^+$ which corresponds to $\k_\sM|_{\PP^+(\L_{\oe})\JJ^1_\sM}$ under the canonical bijection of Lemma~\ref{exttransfer}.
\end{Definition}

The definition of $\b$-extension, together with
Lemma~\ref{exttransfer}, immediately gives us a lower bound for their
intertwining

\begin{Corollary}\label{intkappa}
Let $\k$ be a $\beta$-extension of $\eta$ relative to
$\L^\sM$. Then $\PP^+(\L^\sM_{\oe})\subset I_{G^+}(\k)$.
\end{Corollary}

\begin{proof}
Note first that
certainly $\PP^+(\L^\sM_{\oe})\subset I_{G^+}(\k_\sM)$, since
$\PP^+(\L^\sM_{\oe})\subset \JJ^+_\sM$. In particular, $\PP^+(\L^\sM_{\oe})$
intertwines $\k_\sM|_{\PP^+(\L_{\oe})\JJ^1_\sM}$ so, by
Lemma~\ref{exttransfer}, it also intertwines $\k$.
\end{proof}

This is, of course, much less than the situation for simple characters
for $\tG$, where
$\b$-extensions are intertwined by the whole $\tG$-centralizer of
$\b$. When $\PP(\L^\sM_{\oe})$ is a good maximal compact subgroup of $G_E$ it
should also be possible to prove such a result here, using the Cartan
decomposition and imitating the proofs in~\cite[\S5.2]{BK}. 
Alternatively, if $G$ is an unramified unitary group
then the techniques of \S\ref{S4} (see also~\cite{Sec2}) may provide a
proof of this for any $\L^\sM$. In general, the situation is much less
clear. Although, in~\S\ref{S4.beta} we do improve our lower bound for
$\b$-extensions, we do not pursue these strong intertwining results
here since we will not need them. 

\medskip

It is useful also to note that we can describe the restriction of a 
$\b$-extension to the pro-$p$ Sylow subgroup of $\JJ^+(\b,\L)$ in terms
of the representation given by Proposition~\ref{etatilde}. The proof is
identical to that of~\cite[Proposition 5.2.6]{BK}.

\begin{Proposition}[cf.\ {\cite[Proposition 5.2.6]{BK}}]
\label{kapparesetatilde} 
Let $[\L,n,0,\b]$, $[\L^\sm,n_\sm,0,\b]$ be skew semisimple strata
with $\fb_0(\L^\sm)$ a minimal self-dual $\oe$-order in $B$ and
$\fa_0(\L^\sm)\subset\fa_0(\L)$. Let $\th\in\CC_-(\L,0,\b)$, let
$\th_{\sm}=\tau_{\L^\sm,\L,\b}(\th)$, and let $\eta$ (respectively 
$\eta_\sm$) be the irreducible representation of $\JJ^1(\b,\L)$
(respectively, $\JJ^1(\b,\L^\sm)$) given by
Proposition~\ref{eta}. Put
$\JJ^1_{\sm,\L}=\PP_1(\L^\sm_{\oe})\JJ^1(\b,\L)$, a pro-$p$ Sylow
subgroup of $\JJ^+(\b,\L)$, and let $\eta_{\sm,\L}$ be the irreducible
representation of $\JJ^1_{\sm,\L}$ given by Proposition~\ref{etatilde}
applied with $\L$ in place of $\L^\sM$. Then, for any $\b$-extension $\k$ 
of $\eta$ to $\JJ^+(\b,\L)$, 
$$
\k|_{\JJ^1_{\sm,\L}} = \eta_{\sm,\L}.
$$
\end{Proposition}

\medskip

If $\k$ is a $\b$-extension of $\eta$ to $\JJ^+(\b,\L)$ then we say
that the restriction $\k|_{\JJ^\so(\b,\L)}$ is a $\b$-extension of
$\eta$ to $\JJ^\so(\b,\L)$. We have a compatibility property analogous
to Lemma~\ref{exttransfer} for extensions to $\JJ^\so(\b,\L)$.
Indeed, the following Proposition shows that this compatibility is a
bijection for $\b$-extensions to $\JJ^\so(\b,\L)$, as in~\cite[(5.2.14)]{BK}.

\begin{Proposition}[cf.\ {\cite[(5.2.14)]{BK}}]\label{compat}
Let $[\L,n,0,\b]$ and $[\L',n',0,\b]$ be skew semisimple strata such 
that there exists $[\L^\sM,n_\sM,0,\b]$ a skew semisimple stratum with 
$\fb_0(\L^\sM)$ a maximal self-dual $\oe$-order containing both $\fb_0(\L)$ 
and $\fb_0(\L')$. Let $\th\in\CC_-(\L,0,\b)$ be a skew semisimple character 
and let $\th'=\tau_{\L,\L',\b}(\th)$. Let $\eta$, $\eta'$ be the 
representations of $\JJ^1(\b,\L)$, $\JJ^1(\b,\L')$ respectively, given by
Proposition~\ref{eta}. Then there is a canonical bijection
between the set of $\b$-extensions $\k'$ of $\eta'$ to $\JJ^\so(\b,\L')$
relative to $\L^\sM$ and the set of $\b$-extensions $\k$ of $\eta$ to
$\JJ^\so(\b,\L)$ relative to $\L^\sM$. If 
$\fa_0(\L)\subset\fa_0(\L')$ then it is given as follows: given $\k$
(respectively $\k'$) there is a unique $\k'$ (respectively $\k$) 
such that $\k'|_{\PP^\so(\L_{\oe})\JJ^1(\b,\L')}$ and $\k$ induce 
equivalent irreducible representations of $\PP^\so(\L_{\oe})\PP_1(\L)$.
\end{Proposition}

We say that $\k'$ is \emph{compatible with $\k$} if they correspond under 
the bijection of Proposition~\ref{compat}.

\begin{proof}
We observe first that, using Lemma~\ref{exttransfer} (or rather, its analogue for $\JJ^\so$), we can reduce to the case where $\fa_0(\L^\sM)$ contains both $\fa_0(\L)$ and $\fa_0(\L')$: Lemma~\ref{lattcont} shows that there is a self-dual $\oe$-lattice sequence $\L''$ with $\fb_0(\L'')=\fb_0(\L)$ and $\fa_0(\L'')\subset\fa_0(\L^\sM)$ and, by Lemma~\ref{exttransfer}, there is a canonical bijection between the extensions of $\eta$ to $\JJ^\so(\b,\L)$ and the extensions of $\eta''$ to $\JJ^\so(\b,\L'')$. Then we replace $\L$ by $\L''$ (and similarly for $\L'$).

When $\fa_0(\L)\subset\fa_0(\L')$, 
from the commutativity of the diagram in Lemma~\ref{lem:firstcompat}, and by
restriction, we get a well-defined 
\emph{surjection}
\begin{equation}\label{betabij}
\left\{\b\hbox{-extensions of }\eta'\hbox{ to }\JJ^\so(\b,\L')\hbox{ relative to }\L^\sM\right\} 
\ \to\ 
\left\{\b\hbox{-extensions of }\eta\hbox{ to }\JJ^\so(\b,\L)\hbox{ relative to }\L^\sM\right\}. 
\end{equation}
To see that we have a bijection, we need the following Lemma:

\begin{Lemma}\label{betamaxres}
Let $[\L^\sM,n_\sM,0,\b]$ be a skew semisimple stratum with
$\fb_0(\L^\sM)$ a maximal self-dual $\oe$-order in $B$ containing $\fb_0(\L)$. 
Let $\th_\sM\in\CC_-(\L^\sM,0,\b)$, let $\eta_\sM$ be the representation
of $\JJ^1_\sM=\JJ^1(\b,\L^\sM)$ given by Proposition~\ref{eta}, and let 
$\k_\sM,\k'_\sM$ be $\b$-extensions of $\eta_\sM$ to 
$\JJ^+_\sM=\JJ^+(\b,\L^\sM)$, given by Theorem~\ref{betamax}. 
Put $\JJ^\so_\sM=\JJ^\so(\b,\L^\sM)$ and 
$\JJ^\so_{\L,\sM}=\PP^\so(\L_{\oe})\JJ^1_\sM$. Then
$$
\k_\sM|_{\JJ^\so_{\L,\sM}} \simeq \k'_\sM|_{\JJ^\so_{\L,\sM}}
\quad\iff\quad
\k_\sM|_{\JJ^\so_\sM} \simeq \k'_{\sM}|_{\JJ^\so_\sM}.
$$
\end{Lemma}

\begin{proof} By Theorem~\ref{betamax}, we have
$\k'_\sM=\k_\sM\otimes\chi$, for some character $\chi$ of
$\CG^+=\PP^+(\L^\sM_{\oe})/\PP_1(\L^\sM_{\oe})$ which is trivial on the
subgroup generated by all its unipotent subgroups. If 
$\k_\sM|_{\JJ^\so_{\L,\sM}} \simeq \k'_\sM|_{\JJ^\so_{\L,\sM}}$ then
$\chi$ is also trivial on $\CP^\so=\PP^\so(\L_{\oe})/\PP_1(\L^\sM_{\oe})$,
which is a parabolic subgroup of $\CG^\so$. Hence $\chi$ is trivial on all
of $\CG^\so$ and $\k_\sM|_{\JJ^\so_\sM} \simeq \k'_{\sM}|_{\JJ^\so_\sM}$. 
The converse is trivial.
\end{proof}

Returning to Proposition~\ref{compat}, we see that Lemmas~\ref{betamaxres} 
and~\ref{exttransfer} (together with the analogous Lemma for
$\JJ^\so(\b,\L)$) imply that both sides of \eqref{betabij} have the 
same number of elements, which completes the proof when 
$\fa_0(\L)\subset\fa_0(\L')$. The general case follows immediately.
\end{proof}

{\bf Remark} In the simple case for $\tG=GL_N(F)$, all the maximal
$\oe$-orders in $B$ containing $\fb_0(\L)$ are conjugate by $\tG_E$ so
the set of $\b$-extensions relative to $\L^\sM$ is independent of the 
choice of $\L^\sM$. This is no longer true in our situation. However, we 
will see in Corollary~\ref{betatwist} that there is some sort of compatibility between 
the $\b$-extensions defined relative to different maximal self-dual
$\oe$-orders in $B$. On the other hand, it will be useful in many 
situations to have a ``canonical'' choice for $\L^\sM$, as follows:

Given a skew semisimple stratum $[\L,n,0,\b]$, write $\L=\oplus_i
\L^i$, with $\L^i$ an $\fo_{E_i}$-lattice sequence. Define the $\fo_{E_i}$
lattice sequence $\fM_\L^i$ by 
\begin{equation}\label{standardmax}
\fM_\L^i(2r+s) = \fp_{E_i}^r\L^i(s),\qquad\hbox{for }r\in\BZ,s=0,1,
\end{equation}
and put $\fM_\L=\oplus_{i}\fM_\L^i$. Then $[\fM_\L,n_\fM,0,\b]$ is a skew
semisimple stratum, for some integer $n_\fM$, and
$\fb_0(\fM_\L)$ is a maximal self-dual $\oe$-order containing
$\fb_0(\L)$. 

Let $\th\in\CC_-(\L,0,\b)$ and let $\eta$ be the unique irreducible
representation of 
$\JJ^1(\b,\L)$ containing $\th$. We call a $\b$-extension
of $\eta$ relative to $\fM_\L$ a \emph{standard $\b$-extension}. 

\medskip

In particular, if $[\L,n,0,\b]$ and $[\L',n',0,\b]$ are skew semisimple
strata with $\fM_\L=\fM_{\L'}$ then, by Proposition~\ref{compat} there
is a bijection between the standard $\b$-extensions of $\eta$ to
$\JJ^\so(\b,\L)$ and those of $\eta'$ to $\JJ^\so(\b,\L')$, where
$\eta'$ is the unique irreducible representation of $\JJ^1(\b,\L')$
containing $\th'=\tau_{\L,\L',\b}(\th)$.


\subsection{Analogous results for $\tG$}

If $[\L,n,0,\b]$ is a semisimple stratum and $\tth\in\CC(\L,0,\b)$ is a
semisimple character then the results of the previous sections can be
imitated to construct $\b$-extensions $\til\k$ of $\teta$ to
$\tJ(\b,\L)$. Given Proposition~\ref{etatildetG}, the proofs are identical 
to those for $G$ and are left as an exercise for the reader.
In this situation the notion of $\b$-extension is
independent of the choice of $\L^\sM$ such that $\fb_0(\L^\sM)$ is a
maximal $\oe$-order containing $\fb_0(\L)$, the reductive quotient
$\tJ(\b,\L)/\tJ^1(\b,\L)$ is always connected and there is a canonical
bijection between $\b$-extensions, as in
Proposition~\ref{compat}. 

We will not, in fact, use these results for semisimple characters, but
only for simple characters, when they all come
from~\cite[\S5.2]{BK}.


\Section{Iwahori factorizations}\label{S4}

In this section we imitate the constructions of~\cite[\S7]{BK} to get 
Iwahori factorizations of semisimple characters. 
We do this first for $\tG$ and then, by Glauberman transfer, for $G^+$
and $G$. We end with some similar results for $\b$-extensions.


\subsection{Iwahori decompositions}\label{S4.1}

Let $[\L,n,0,\b]$ be a semisimple stratum with associated splitting
$V=\bigoplus_{i=1}^l V^i$ and all the other usual notation. 
Let $V=\bigoplus_{j=1}^{m} W^{(j)}$ be a
decomposition of $V$ into subspaces such that
\begin{enumerate}
\item $W^{(j)}=\bigoplus_{i=1}^l \left(W^{(j)}\cap V^i\right)$, for
$1\le j\le m$, so also $V^i=\bigoplus_{j=1}^m \left(W^{(j)}\cap
  V^i\right)$, for $1\le i\le l$;
\item $W^{(j)}\cap V^i$ is an $E_i$-subspace of $V^i$, for $1\le j\le m$
and $1\le i\le l$.
\end{enumerate}

\begin{Definition}\label{defsub}\rm We say that the decomposition 
$V=\bigoplus_{j=1}^m W^{(j)}$ is:
\begin{enumerate}
\item \emph{subordinate to $[\L,n,0.\b]$} if, for all $r\in\BZ$,
  $\L(r)=\bigoplus_{j=1}^m\left(\L(r)\cap W^{(j)}\right)$;
\item \emph{properly subordinate to $[\L,n,0,\b]$} if it is subordinate 
to $[\L,n,0,\b]$ and, for all $r\in\BZ$ and for $1\le i\le l$, there is 
at most one $j$, $1\le j\le m$, such that 
$$
\left(\L(r)\cap W^{(j)}\cap V^i\right) 
\supsetneq 
\left(\L(r+1)\cap W^{(j)}\cap V^i\right).
$$
\end{enumerate}
\end{Definition}

An explanation of how these decompositions arise is in order at this
point. We consider first the case when the stratum $[\L,n,0,\b]$ is
simple. Let $\CB$ be an $E$-basis for $V$ which splits the lattice
sequence $\L$. If we decompose $\CB$ as a disjoint union of subsets
$\CB_j$, $1\le j\le m$, and define $W^{(j)}$ to be the $E$-linear span
of $\CB_j$, then we obtain a decomposition of $V$ subordinate to the
stratum. Moreover, any decomposition subordinate to the stratum arises
in this way.

For a properly subordinate decomposition, we need to be a little more
careful. Writing $e=e(\L|\oe)$ for the $\oe$-period of $\L$, we
decompose $\CB$ into a disjoint union of subsets $\CB_j$, $1\le j\le
e$, by saying that a vector $\bs w\in\CB$ lies in $\CB_j$ if and only
if there exists $r\in\BZ$ such that $\bs w\in\L(j+re) \setminus
\L(j+re+1)$. Let $W^{(j)}$ be the $E$-linear space of $\CB_j$. Then
any decomposition of $V$ which coarsens the decomposition
$V=\bigoplus_{j=1}^e W^{(j)}$ (that is, each subspace is a sum of
some $W^{(j)}$) is properly subordinate to the stratum. Again, any
decomposition properly subordinate to the stratum arises in this
way. 

\medskip

For the semisimple case we take, for $1\le i\le l$, a decomposition
$V^i=\bigoplus_{j=1}^m W^{(i,j)}$ which is (properly) subordinate to
the simple stratum $[\L^i,n_i,0,\b_i]$ (and where we allow some of the
subspaces $W^{(i,j)}$ to be trivial). Then we get a decomposition of
$V$ which is (properly) subordinate to $[\L,n,0,\b]$ by setting
$W^{(j)}=\bigoplus_{i=1}^l W^{(i,j)}$.

\begin{Proposition}[cf.\ {\cite[Theorem 7.1.14]{BK}}]\label{IwHJ}
Let $[\L,n,0,\b]$ be a semisimple stratum and let 
$V=\bigoplus_{j=1}^m W^{(j)}$ be a decomposition as above. Let $\tM$ be
the Levi subgroup of $\tG$ which is the stabilizer of this
decomposition and let $\tP$ be any parabolic subgroup with Levi
component $\tM$ and unipotent radical $\tU$.
\begin{enumerate}
\item If the decomposition is subordinate to $[\L,n,0,\b]$, then the 
groups $\tH^1(\b,\L)$ and $\tJ^1(\b,\L)$ have Iwahori decompositions with 
respect to $(\tM,\tP)$.
\item If the decomposition is properly subordinate to $[\L,n,0,\b]$, 
then the group $\tJ(\b,\L)$ also has an Iwahori decomposition with 
respect to $(\tM,\tP)$ and, moreover, $\tJ(\b,\L)\cap\tU=\tJ^1(\b,\L)\cap\tU$. 
\end{enumerate}
\end{Proposition}

\begin{proof} (i) For $1\le j\le m$, let $\boe^{(j)}$ denote the projection
onto $W^{(j)}$ with kernel $\bigoplus_{k\ne j} W^{(k)}$. The
definition of subordinate (Definition~\ref{defsub}(i)) ensures that
$\boe^{(j)}\in\fa_0(\L)$. Moreover, 
the property that each $W^{(j)}\cap V^i$ is an $E_i$-subspace of $V^i$,
for $1\le i\le l$, ensures that $\boe^{(j)}\in B$. Hence
$\boe^{(j)}\in\fb_0(\L)$. 

Now $\tJ^1(\b,\L)=1+\fJ^1(\b,\L)$ and, by~\cite[Lemma 3.10]{S4},
$\fJ^1(\b,\L)$ is a $\fb_0(\L)$-bimodule. In particular,
$\boe^{(j)}\fJ^1(\b,\L)\boe^{(k)}\subset\fJ^1(\b,\L)$ so, 
by~\cite[Proposition 10.4]{BH1}, $\tJ^1(\b,\L)$ has an Iwahori
decomposition with respect to $(\tM,\tP)$. The same proof works for
$\tH^1(\b,\L)$. 

(ii) Suppose now that the decomposition is properly subordinate to
$[\L,n,0,\b]$. We start with a lemma. 

\begin{Lemma}[cf.\ {\cite[Lemma 7.1.15]{BK}}]\label{properPE}
In this situation we
have $\til\PP(\L_{\oe})\cap \tU=\til\PP_1(\L_{\oe})\cap\tU$.
\end{Lemma}

\begin{proof} We prove the corresponding additive statement, which is
$\boe^{(j)}\fb_0(\L)\boe^{(k)}\subset\fb_1(\L)$, for $j\ne k$. So suppose
$x\in\boe^{(j)}\fb_0(\L)\boe^{(k)}$. For any $r\in\BZ$ and for $s\ne k$ we
have $x(\L(r)\cap W^{(s)})=0$. On the other hand, for any $r\in\BZ$ and
$1\le i\le l$,
either $\L(r)\cap W^{(k)}\cap V^i=\L(r+1)\cap W^{(k)}\cap V^i$ or
$\L(r)\cap W^{(j)}\cap V^i=\L(r+1)\cap
W^{(j)}\cap V^i$, since the decomposition is properly subordinate. In
either case, we get $x(\L(r)\cap W^{(k)}\cap V^i)\subset \L(r+1)\cap
W^{(j)}\cap V^i$. Hence
$$
x\L(r) = x\left(\L(r)\cap W^{(k)}\right) 
= \bigoplus_{i=1}^l x\left(\L(r)\cap W^{(k)}\cap V^i\right)
\subset \L(r+1),
$$ 
and $x\in\fa_1(\L)$ as required.
\end{proof}

Now the same proof as that of~\cite[Theorem 7.1.14]{BK} shows that
$\til\PP(\L_{\oe})$ has an Iwahori decomposition with respect to
$(\tM,\tP)$ and that $\til\PP(\L_{\oe})=\left(\til\PP(\L_{\oe})\cap\tM\right)
\til\PP_1(\L_{\oe})$. Hence $\JJ(\b,\L)=\left(\til\PP(\L_{\oe})\cap\tM\right) 
\JJ^1(\b,\L)$ has an Iwahori decomposition with respect to $(\tM,\tP)$,
and the final assertion follows also.
\end{proof}

Let $[\L,n,0,\b]$ be a semisimple
stratum and let $V=\bigoplus_{j=1}^m W^{(j)}$ be a decomposition subordinate
to the stratum. We write $\L^{(j)}$ for the lattice sequence in
$W^{(j)}$ given by
$$
\L^{(j)}(r)=\L(r)\cap W^{(j)},\quad\hbox{for }r\in\BZ,
$$
and put $\b^{(j)}=\boe^{(j)}\b\boe^{(j)}$, where $\boe^{(j)}$ denotes
the projection onto $W^{(j)}$ with kernel $\bigoplus_{k\ne j} W^{(k)}$. 
Then, since each $V^i\cap W^{(j)}$ is an $E_i$-subspace of $V^i$, the
stratum $[\L^{(j)},n^{(j)},0,\b^{(j)}]$ is a semisimple
stratum in $A^{(j)}=\End_F(W^{(j)})$, for some integer $n^{(j)}$, with
splitting $W^{(j)}=\bigoplus_{i=1}^l \left(V^i\cap
W^{(j)}\right)$. Note that some of the pieces of this splitting may be
trivial so, strictly speaking, the splitting consists only of the
non-zero pieces. We also write $B^{(j)}$ for the centralizer of
$\b^{(j)}$ in $A^{(j)}$. 

\begin{Proposition}\label{HintM}
 Let $[\L,n,0,\b]$ be a semisimple stratum and let 
$V=\bigoplus_{j=1}^m W^{(j)}$ be a decomposition subordinate to the
stratum. Let $\tM$ be the Levi subgroup of $\tG$ which is the
stabilizer of this decomposition. Then
$$
\tH^1(\b,\L)\cap\tM = \prod_{j=1}^m \tH^1(\b^{(j)},\L^{(j)}),
$$
and likewise for $\tJ^1(\b,\L)$. If the decomposition is properly
subordinate to the stratum then we also have
$$
\tJ(\b,\L)\cap\tM = \prod_{j=1}^m \tJ(\b^{(j)},\L^{(j)}).
$$
\end{Proposition}

\begin{proof} 
We will prove the additive statement
$$
\boe^{(j)}\fH(\b,\L)\boe^{(j)} = \fH(\b^{(j)},\L^{(j)}).
$$
The same proof will give the corresponding statement for $\fJ(\b,\L)$ and the 
result follows.

Let $\CM$ denote the stabilizer in $A$ of the decomposition, so that $\tM=\CM^\times$. We proceed by induction on $\ko$, noting that the base 
case is a null stratum, where the result is clear since 
$\fH(\b,\L)=\fa_0(\L)$. So we turn to the inductive step and put $r=-\ko$.

Put $W^{(i,j)}=V^i\cap W^{(j)}$ so that $V^i=\bigoplus_{j=1}^m W^{(i,j)}$ 
is (properly) subordinate to the simple stratum $[\L^i,n_i,0,\b_i]$. 
We choose an $E_i$-basis $\CB_{i,j}$ for each $W^{(i,j)}$ which splits 
the lattice sequence $\L^{(i,j)}$ given by 
$\L^{(i,j)}(r)=\L(r)\cap W^{(i,j)}$, and let $Y^{(i,j)}$ be the 
$F$-linear span of $\CB_{i,j}$. Putting $Y^i=\bigoplus_{j=1}^m Y^{(i,j)}$, 
we get a ``generalized $(W,E)$-decomposition'' $V^i=Y^i\otimes_F E_i$
(see~\cite[\S5.3]{BK1}), and 
hence an embedding $\iota_{Y^i}:\End_F(E_i)\into A^{ii}=\End_F(V^i)$ with 
image in $A^{ii}\cap\CM$.

According to~\cite[Proposition 3.4]{S4}, we can choose a semisimple stratum 
$[\L,n,r,\g]$ equivalent to $[\L,n,r,\b]$ with $\g$ in the image of $\prod_i \iota_{Y^i}$. In particular, we have $\g\in\CM$ so the decomposition is subordinate to $[\L,n,r,\g]$ also. Hence, by the inductive hypothesis, we 
have
$$
\boe^{(j)}\fH(\g,\L)\boe^{(j)} = \fH(\g^{(j)},\L^{(j)}),
$$
where $\g^{(j)}=\boe^{(j)}\g\boe^{(j)}$. 
The semisimple stratum $[\L^{(j)},n^{(j)},r,\g^{(j)}]$ 
is equivalent to $[\L^{(j)},n^{(j)},r,\b^{(j)}]$ so, since 
$\boe^{(j)}\in\fb_0(\L)$, we get
\begin{eqnarray*}
\boe^{(j)}\fH(\b,\L)\boe^{(j)} &=& \boe^{(j)}\fb_0(\L)\boe^{(j)} + 
\boe^{(j)}\fH^{\left[\frac r2\right]+1}(\g,\L)\boe^{(j)} \\
&=& \fb_0(\L^{(j)}) + \fH^{\left[\frac r2\right]+1}(\g^{(j)},\L^{(j)}) \\
&=& \fH(\b^{(j)},\L^{(j)}),
\end{eqnarray*}
as required.
\end{proof}

Continuing in the same situation, let $\tP$ be any parabolic subgroup
of $\tG$ with Levi subgroup $\tM$, and let $\tU$ denote its unipotent
radical. We define the groups
$$
\tH^1_{\tP}=\tH^1(\b,\L)\left(\tJ^1(\b,\L)\cap\tU\right),
\qquad\hbox{and}\qquad
\tJ^1_{\tP}=\tH^1(\b,\L)\left(\tJ^1(\b,\L)\cap\tP\right).
$$
If the decomposition is properly subordinate to the stratum, we also put
$$
\tJ_{\tP}=\tH^1(\b,\L)\left(\tJ(\b,\L)\cap\tP\right).
$$
All these groups have Iwahori decompositions with respect to any
parabolic subgroup with Levi component $\tM$.


\subsection{Iwahori factorization of semisimple characters}\label{S4.2}

We continue in the same situation, so
$[\L,n,0,\b]$ is a semisimple stratum and  
$V=\bigoplus_{j=1}^m W^{(j)}$ is a decomposition subordinate to the
stratum. Let $\tM$ be
the Levi subgroup of $\tG$ which is the stabilizer of this
decomposition and let $\tP$ be any parabolic subgroup with Levi
component $\tM$. Write $\tU$ for the unipotent radical of $\tP$.

\begin{Proposition}[cf.\ {\cite[Proposition 7.1.19]{BK}}]\label{Iwtheta}
In this situation, let $\tth\in\CC(\L,0,\b)$. Then the character
$\tth|_{\tH^1(\b,\L)\cap\tU}$ is trivial. After identifying
$\tH^1(\b,\L)\cap\tM$ with 
$\prod_{j=1}^m\tH^1(\b^{(j)},\L^{(j)})$, we have
$$
\tth|_{\tH^1(\b,\L)\cap\tM} = \bigotimes_{j=1}^m \tth^{(j)},
$$
where $\tth^{(j)}
\in \CC(\L^{(j)},0,\b^{(j)})$.
\end{Proposition}

In this situation, we shall write $\tth^{(j)}=\tau_{\L,\L^{(j)},\b}(\tth)$. 
By the definition of simple character in~\cite[\S5]{BK1}, this is
consistent with the notion of transfer between 
different algebras for simple characters. We 
do not pursue the functorial properties of this transfer map but note that it 
gives a map $\tau_{\L,\L^{(j)},\b}:\CC(\L,0,\b) \to \CC(\L^{(j)},0,\b^{(j)})$ 
which is \emph{not} in general injective; this happens if one of the subspaces 
$W^{(j)}\cap V^i=\{0\}$ but $\b_i\ne 0$.

\begin{proof} 
We proceed by induction on $\ko$, the base case again being 
the null case, when the only semisimple character is the trivial one and 
the result is clear. So we consider the induction step and put $r=-\ko$. 
We will use again the notation of the proof of Proposition~\ref{HintM}, 
so that $[\L,n,r,\g]$ is the carefully chosen semisimple stratum equivalent 
to $[\L,n,r,\b]$ from there. 

Recall that $\tth$ has the properties:
\begin{enumerate}
\item we have $\tth|_{\tH^1\cap\tG_i}=\tth_i$, for some simple character 
$\tth_i\in\CC(\L^i,0,\b_i)$, where $\tG_i=\Aut_F(V^i)$;
\item $\tth|_{\tH^{\left[\frac r2\right]+1}(\b,\L)} 
= \tth_0|_{\tH^{\left[\frac r2\right]+1}(\b,\L)} \psi_{\b-\g}$, for
some semisimple 
character $\tth_0\in\CC(\L,0,\g)$.
\end{enumerate}
Moreover, these properties characterize the semisimple characters in 
$\CC(\L,0,\b)$. 

For the first assertion of the Lemma, let $\tM_\b$ denote the
stabilizer of the splitting $V=\bigoplus_{i=1}^l V^i$ and let $\tU_\b$
be the unipotent radical of some parabolic subgroup $\tP_\b$ with Levi
component 
$\tM_\b$. According to~\cite[Lemma 3.15(i)]{S4}, the restriction
$\tth|_{\tH^1\cap\tU_\b}$ is trivial. Since $\tH^1\cap\tU$ has an
Iwahori decomposition with respect to $(\tM_\b,\tP_\b)$, we are
reduced to showing that $\tth$ restricts trivially to
$\tH^1\cap\tU\cap\tM_\b$. But here we may work block-by-block so the
result follows from the simple case in~\cite[Proposition 5.2]{BK1}.

Now we turn to the second assertion.
From Proposition~\ref{HintM} we have that $\tH^1(\b,\L)\cap\Aut_F(W^{(j)}) 
=\tH^1(\b^{(j)},\L^{(j)})$ so we may consider $\tth^{(j)}
=\tth|_{\tH^1(\b^{(j)},\L^{(j)})}$. To show that it is a 
semisimple character in $\CC(\L^{(j)},0,\b^{(j)})$ we need 
only check that it satisfies the two properties (i), (ii) above.

The splitting of  
$[\L^{(j)},n^{(j)},0,\b^{(j)}]$ is given by 
$W^{(j)}=\bigoplus_{i=1}^l W^{(i,j)}$, where $W^{(i,j)}=V^i\cap W^{(j)}$ 
and we ignore any trivial pieces of the splitting. We put 
$\tG^{(i,j)}=\Aut_F(W^{(i,j)})$, let $\boI^{(i,j)}$ be the projection 
from $W^{(j)}$ to $W^{(i,j)}$ with kernel $\bigoplus_{k\ne i} W^{(k,j)}$, 
and let $\boe^{(i,j)}$ be the projection from $V^i$ to $W^{(i,j)}$ with 
kernel $\bigoplus_{k\ne j} W^{(i,k)}$. 
We put $\b^{(j)}_i=\boe^{(i,j)}\b_i\boe^{(i,j)}=
\boI^{(i,j)}\b^{(j)}\boI^{(i,j)}$.

(i) We have $\tH^1(\b^{(j)},\L^{(j)})\cap\tG^{(i,j)} = 
\tH^1(\b^{(j)}_i,\L^{(i,j)})$ and 
$$
\tth^{(j)}|_{\tH^1(\b^{(j)},\L^{(j)})\cap\tG^{(i,j)}} 
\ =\ \left.\left(\tth|_{\tH^1(\b,\L)\cap\tG^i}\right)
\right|_{\tH^1(\b^{(j)}_i,\L^{(i,j)})}
\ =\ \tth_i|_{\tH^1(\b^{(j)}_i,\L^{(i,j)})},
$$
which, by~\cite[\S5]{BK1}, is the simple character 
$\tau_{\L^i,\L^{(i,j)},\b_i}(\tth_i)$ in $\CC(\L^{(i,j)},0,\b^{(j)}_i)$.

(ii) The stratum $[\L^{(j)},n^{(j)},0,\g^{(j)}]$ is a semisimple stratum 
equivalent to $[\L^{(j)},n^{(j)},0,\b^{(j)}]$ and it is straightforward 
that $\psi_{\b-\g}|_{\tH^{[\frac r2]+1}(\b^{(j)},\L^{(j)})} 
= \psi_{\b^{(j)}-\g^{(j)}}$. Thus we have
\begin{eqnarray*}
\tth^{(j)}|_{\tH^{[\frac r2]+1}(\b^{(j)},\L^{(j)})}& =& 
\left.\left(\tth|_{\tH^{[\frac r2]+1}(\b,\L)}\right)
\right|_{\tH^{[\frac r2]+1}(\b^{(j)},\L^{(j)})} \\ &=& 
\left.\left(\tth_0|_{\tH^{[\frac r2]+1}(\b,\L)} \psi_{\b-\g}\right)
\right|_{\tH^{[\frac r2]+1}(\b^{(j)},\L^{(j)})} \ =\ 
\tth_0^{(j)}|_{\tH^{[\frac r2]+1}(\b^{(j)},\L^{(j)})} 
\psi_{\b^{(j)}-\g^{(j)}},
\end{eqnarray*}
for some $\tth_0^{(j)}\in\CC(\L^{(j)},0,\g^{(j)})$, by the inductive 
hypothesis.
\end{proof}

We continue with the same situation. From~\cite[Proposition 3.24]{S4}, 
the pairing
$$
\bok_\tth(x,y)\ =\ \tth[x,y],\qquad x,y\in\tJ^1(\b,\L),
$$
induces a nondegenerate alternating form on $\tJ^1(\b,\L)/\tH^1(\b,\L)$. 
Likewise, for $j=1,...,m$ we have a form $\bok_{\tth^{(j)}}$ on 
$\tJ^1(\b^{(j)},\L^{(j)})/\tH^1(\b^{(j)},\L^{(j)})$.

Let $\tU_l$ denote the
unipotent radical of the parabolic subgroup opposite $\tP$. As an
immediate consequence of Proposition~\ref{Iwtheta} and the Iwahori
decompositions for $\tH^1(\b,\L)$ and $\tJ^1(\b,\L)$, we get:

\begin{Lemma}[cf.\ {\cite[Proposition 7.2.3]{BK}}]\label{orththeta}
\begin{enumerate} 
\item The subspaces $\tJ^1\cap \tU_l/\tH^1\cap \tU_l$ and 
$\tJ^1\cap \tU/\tH^1\cap \tU$ of $\tJ^1/\tH^1$ are both totally
isotropic for the form $\bok_\tth$ and orthogonal to 
$\tJ^1\cap\tM/\tH^1\cap\tM$.
\item The restriction of $\bok_\tth$ to the group
$$
\tJ^1\cap\tM/\tH^1\cap\tM = \prod_{j=1}^m 
\tJ^1(\b^{(j)},\L^{(j)}) /
\tH^1(\b^{(j)},\L^{(j)})
$$
is the orthogonal sum of the pairings $\bok_{\tth^{(j)}}$, where 
$\tth^{(j)}=\tau_{\L,\L^{(j)},\b}(\tth)$.
\item We have an orthogonal sum decomposition
$$
\frac{\tJ^1}{\tH^1} = \frac{\tJ^1\cap\tM}{\tH^1\cap\tM} \perp \left(
\frac{\tJ^1\cap\tU_l}{\tH^1\cap\tU_l}\times 
\frac{\tJ^1\cap\tU}{\tH^1\cap\tU}\right).
$$
\end{enumerate}
\end{Lemma}

Recall that we have defined the group
$\tH^1_{\tP}=\tH^1(\b,\L)\left(\tJ^1(\b,\L)\cap\tU\right)$, and, similarly, 
$\tJ^1_{\tP}$. We define the
character $\tth_{\tP}$ of $\tH^1_{\tP}$ by
$$
\tth_{\tP}(hj)=\tth(h),\quad\hbox{ for }h\in\tH^1(\b,\L),\ 
j\in\tJ^1(\b,\L)\cap\tU.
$$
This is well-defined by Lemma~\ref{orththeta}. As in~\cite[\S7.2]{BK},
we immediately get

\begin{Corollary}[cf.\ {\cite[Propositions 7.2.4, 7.2.9]{BK}}]
\label{etaP}
There is a unique irreducible representation $\teta_{\tP}$ of
$\tJ^1_{\tP}$ such that $\teta_{\tP}|_{\tH^1_{\tP}}$ contains
$\tth_{\tP}$. Moreover, $\teta\simeq \Indu{\tJ^1_{\tP}}{\tJ^1(\b,\L)}
{\teta_{\tP}}$ and, for each $y\in \tG_E$, there is a unique 
$(\tJ^1_{\tP},\tJ^1_{\tP})$-double coset in
$\tJ^1(\b,\L)y\tJ^1(\b,\L)$ which intertwines $\teta_{\tP}$.
\end{Corollary}

\begin{Lemma}[cf.\ {\cite[Proposition 2.4]{GKS}}]\label{intetaP}
We have $I_{\tG}(\tth_{\tP})=I_{\tG}(\teta_{\tP})=\tJ^1_{\tP} \tG_E
\tJ^1_{\tP}$.
\end{Lemma}

\begin{proof} We have $I_{\tG}(\tth_{\tP})=I_{\tG}(\teta_{\tP})$ so, by
Corollary~\ref{etaP}, we need only check that all of $\tG_E$
intertwines $\tth_{\tP}$.

Let $\CB_{i,j}$ be an $E_i$-basis for $W^{(i,j)}=V^i\cap W^{(j)}$ which 
splits the lattice sequence $\L^{(i,j)}$ and consider the decomposition 
\begin{equation}\label{mindecomp}
V=\bigoplus_{i,j}\bigoplus_{\bs v\in\CB_{i,j}} E_i\bs v.
\end{equation}
This is subordinate to the stratum $[\L,n,0,\b]$ and refines both the
decomposition $V=\bigoplus_{j=1}^m W^{(j)}$ and the splitting
$V=\bigoplus_{i=1}^l V^i$ of the stratum. 
Let $\tM_0$ denote the stabilizer of the decomposition~\eqref{mindecomp} 
and let
$\tP_0\subset\tP$ be a parabolic subgroup with Levi component
$\tM_0$. We write $\tU_0$ for the unipotent radical of $\tP_0$, so
$\tU_0\supset\tU$.

Let $\CP_0\subset A$ denote the Lie algebra of $\tP_0$. Then
$\CP_0\cap\fb_0(\L) + \fb_1(\L)$ is a minimal $\oe$-order in $B$,
which must be of the form $\fb_0(\L^\sm)$, for some
$\oe$-lattice sequence $\L^\sm$ in $V$. Then, for a suitable integer
$n_\sm$, the stratum $[\L^\sm,n_\sm,0,\b]$ is semisimple and the
decomposition \eqref{mindecomp} is subordinate to it. Moreover,
$\til\PP(\L^\sm_{\oe})$ is an Iwahori subgroup of $\tG_E$ contained in
$\til\PP(\L_{\oe})$. 

\medskip

Now, by construction, $\til\PP(\L^\sm_{\oe}) =
(\til\PP(\L_{\oe})\cap\tP_0)\til\PP^1(\L_{\oe})\subset
\tJ_{\tP}$ and, since $\tJ_{\tP}$ normalizes $\tth_{\tP}$, we need
only show that some set of double coset
representatives for $\til\PP(\L^\sm_{\oe})\backslash \tG_E \slash
\til\PP(\L^\sm_{\oe})$ intertwines $\tth_{\tP}$. By the Bruhat decomposition,
we may take these representatives to be in the $\tG_E$-normalizer
of the maximal torus $\tM_0\cap \tG_E$. 

\medskip

So suppose $y\in \tG_E$ normalizes $\tM_0$. Since $\tH^1_{\tP}$ has
Iwahori decompositions with respect to both $(\tM_0,\tP_0)$ and
$(\tM_0,\tP_0^y)$, the group $\tH^1_{\tP}\cap {}^y\tH^1_{\tP}$ has an
Iwahori decomposition with respect to $(\tM_0,\tP_0)$. Moreover,
$\tth_{\tP}$ is trivial on both $\tH^1_{\tP}\cap\tU_0$ and
$\tH^1_{\tP}\cap\tU_0^y$. Hence, to check that $y$ intertwines
$\tth_{\tP}$, we need only check that it intertwines
$\tth_{\tP}|_{\tH^1_{\tP}\cap\tM_0}$. But, since $\tM_0\subset\tM$
and $\tH_{\tP}^1\cap\tM=\tH^1\cap\tM$,
this restriction is just $\tth|_{\tH^1(\b,\L)\cap\tM_0}$, which is
surely intertwined by $y$, by Proposition~\ref{thetaint}.
\end{proof}

We remark that, in the same way as in~\cite[\S7.2]{BK}, we can
describe $\teta_{\tP}$ as the natural representation of
$\tJ^1_{\tP}$ on the space of $(\tJ^1(\b,\L)\cap\tU)$-fixed vectors in
$\teta$. It is also useful to note that the restriction of $\teta_{\tP}$
to $\tJ^1_{\tP}\cap\tM=\tJ^1(\b,\L)\cap\tM$ is given by
$$
\teta_{\tP}|_{\tJ^1(\b,\L)\cap\tM} \simeq \bigotimes_{j=1}^m \teta^{(j)},
$$
where $\teta^{(j)}$ is the unique irreducible representation of 
$\tJ^1(\b^{(j)},\L^{(j)})$ which contains the character
$\tth^{(j)}=\tau_{\L,\L^{(j)},\b}(\tth)$.


\subsection{Glauberman transfer}\label{S4.glaub}

In this section, we consider the case when we have a skew
semisimple stratum $[\L,n,0,\b]$. We say that a decomposition
$V=\bigoplus_{j=-m}^m W^{(j)}$ is \emph{self-dual} if, for 
$-m\le j\le m$, the orthogonal complement of $W^{(j)}$ is 
$\bigoplus_{k\ne -j}W^{(k)}$; we allow the possibility that $W^{(0)}=\{0\}$. 
Then, immediately from Proposition~\ref{IwHJ} we get:

\begin{Corollary} Let $[\L,n,0,\b]$ be a skew semisimple stratum and let 
$V=\bigoplus_{j=-m}^m W^{(j)}$ be a self-dual decomposition as
above. Let $\tM$ be
the Levi subgroup of $\tG$ which is the stabilizer of this
decomposition and let $\tP$ be any $\s$-stable parabolic subgroup with Levi
component $\tM$. Put $M^+=\tM\cap G^+$, $P^+=\tP\cap G^+$ and $M=\tM\cap G$,
$P=\tP\cap G$.
\begin{enumerate}
\item If the decomposition is subordinate to $[\L,n,0,\b]$, then the 
groups $\HH^1(\b,\L)$ and $\JJ^1(\b,\L)$ have Iwahori decompositions with 
respect to $(M,P)$.
\item If the decomposition is properly subordinate to $[\L,n,0,\b]$, 
then the group $\JJ^+(\b,\L)$ also has an Iwahori decomposition with 
respect to $(M^+,P^+)$ and the groups $\JJ(\b,\L),\JJ^\so(\b,\L)$ have
Iwahori decompositions with respect to $(M,P)$. 
\end{enumerate}
\end{Corollary}

\begin{proof} Note that, if $g\in G^+$ has Iwahori decomposition $g=lmu$,
with respect to some self-dual $(\tM,\tP)$ in $\tG$, then also $g=l^\s
m^\s u^\s$ so, by uniqueness of decomposition, $l,m,u\in G^+$. Now,
since $U\subset G$, if $g\in G$ then $l,m,u\in G$ also.
\end{proof}

Notice that, in this situation, the involution swaps the blocks
$\boe^{(j)}A\boe^{(k)}$ and $\boe^{(-k)}A\boe^{(-j)}$ so that 
the stratum $[\L^{(0)},n^{(0)},0,\b^{(0)}]$ is skew
semisimple. From Proposition~\ref{HintM} we get:

\begin{Corollary}\label{HintMforG}
 Let $[\L,n,0,\b]$ be a skew semisimple stratum and let 
$V=\bigoplus_{j=-m}^m W^{(j)}$ be a self-dual decomposition subordinate to the
stratum. Let $\tM$ be the Levi subgroup of $\tG$ which is the
stabilizer of this decomposition and put $M^+=\tM\cap G^+$, $M=\tM\cap
G$. Then
$$
\HH^1(\b,\L)\cap M \simeq \HH^1(\b^{(0)},\L^{(0)})\times 
\prod_{j=1}^m \tH^1(\b^{(j)},\L^{(j)}),
$$
and there is a similar decomposition for $\JJ^1(\b,\L)$. If the decomposition is properly
subordinate to the stratum then we also have
$$
\JJ^\so(\b,\L)\cap M \simeq \JJ^\so(\b^{(0)},\L^{(0)})\times 
\prod_{j=1}^m \tJ(\b^{(j)},\L^{(j)}),
$$
and similar decompositions for $\JJ^+(\b,\L)\cap M^+$ and
$\JJ(\b,\L)\cap M$.
\end{Corollary}

Note here that $\tH^1(\b^{(j)},\L^{(j)})$ has been identified with
$$
\left\{
(k,h)\in\tH^1(\b^{(-j)},\L^{(-j)})\times\tH^1(\b^{(j)},\L^{(j)})
\hbox{ such that }h\ov k=1
\right\}
$$
(Recall that the involution $\ov{\phantom a}$ swaps the blocks
$\boe^{(j)}A\boe^{(j)}$ and $\boe^{(-j)}A\boe^{(-j)}$.)

\medskip

We can also form the groups $\HH^1_P=\tH^1_{\tP}\cap
G=\HH^1(\b,\L)\left(\JJ^1(\b,\L)\cap U\right)$ and, likewise, $\JJ^1_P$ and
(when the decomposition is properly subordinate to the stratum)
$\JJ^+_P=\HH^1(\b,\L)\left(\JJ^+(\b,\L)\cap P\right)$, and
$\JJ_P$, $\JJ^\so_P$.

\medskip

Now, using the Glauberman correspondence, we transfer the
results of~\S\ref{S4.2} to the case of skew semisimple strata. We
recall briefly the nature of the correspondence, for details of which
(in our situation) we refer the reader to~\cite[\S2]{S2}.

\medskip

Let $\tK$ be a $\s$-stable pro-$p$ subgroup of $\tG$ and
$\KK=\tK^\s$ the group of $\s$-fixed points. Then $\s$ acts on the 
representations $\til\rho$ of $\tK$ by $\til\rho^\s(k)=\rho(k^\s)$, for 
$k\in\tK$. There is a bijection,
denoted $\til\rho\leftrightarrow\rho=\bs g(\til\rho)$, between (equivalence
classes of) irreducible representations of $\tK$ with
$\til\rho^\s\simeq\til\rho$ and (equivalence classes
of) irreducible representations of $\KK$. Further, this
correspondence commutes with irreducible restriction and irreducible
induction. Recall also that the representation
$\rho=\bs g(\til\rho)$ is characterized as the unique component of
$\til\rho|_\KK$ appearing with odd multiplicity; in particular, if 
$\til\rho$ is a character then $\rho=\til\rho|_\KK$. Moreover, for 
$g\in G^+$, we have that $\dim I_g(\til\rho)$ is odd if and only if 
$\dim I_g(\rho)$ is odd; in particular, if $\til\rho$ is a character, or 
if all the intertwining spaces are one-dimensional (as is often the case 
for us), then $I_{G^+}(\rho|\KK)=I_{\tG}(\til\rho|\tK)\cap G^+$.

\medskip

Let $[\L,n,0,\b]$ be a skew semisimple stratum with a subordinate self-dual 
decomposition $V=\bigoplus_{j=-m}^m W^{(j)}$. Let $\tM$ be the stabilizer 
in $\tG$ of this decomposition and let $\tP$ be a $\s$-stable parabolic 
subgroup with Levi component $\tM$. Put $M=\tM\cap G$ and $P=\tP\cap G$.

Let $\th\in\CC_-(\L,0,\b)$ be a skew semisimple character of 
$\HH^1(\b,\L)$, so that $\th=\tth|_{\HH^1(\b,\L)}$, for some 
$\s$-stable semisimple character $\tth\in\CC(\L,0,\b)$. From Corollary~\ref{HintMforG} and Proposition~\ref{Iwtheta}, we have
$$
\th|_{\HH^1(\b,\L)\cap M} = \th^{(0)}\otimes\bigotimes_{j=1}^m 
\left(\tth^{(j)}\right)^2,
$$
for some $\th^{(0)}\in\CC_-(\L^{(0)},0,\b^{(0)})$ and $\tth^{(j)}
\in \CC(\L^{(j)},0,\b^{(j)})$. Notice that we have
$\tH^1(\b^{(j)},\L^{(j)})=\tH^1(2\b^{(j)},\L^{(j)})$ and 
$\left(\tth^{(j)}\right)^2\in \CC(\L^{(j)},0,2\b^{(j)})$, as
in~\cite[\S4.3 Lemma 1]{Bl1}. 

We define the 
character $\th_P$ of $\HH^1_P=\HH^1(\b,\L)\left(\JJ^1(\b,\L)\cap U\right)$ by
$$
\th_P(hj)=\th(h),\quad\hbox{ for }h\in\HH^1(\b,\L),\ 
j\in\JJ^1(\b,\L)\cap U.
$$
Then $\th_P=\tth_{\tP}|_{\HH^1_P}$ and is the Glauberman transfer of 
$\tth_{\tP}$. 

Let $\eta$ be the unique irreducible representation of $\JJ^1(\b,\L)$ 
containing $\th$ and let $\teta$ be the unique irreducible representation 
of $\tJ^1(\b,\L)$ containing $\tth$. Then $\eta=\bs g(\teta)$. Let
$\teta_{\tP}$ be the unique irreducible representation of $\tJ^1_{\tP}$  
containing $\tth_{\tP}$, given by Corollary~\ref{etaP}, and define the 
representation $\eta_P$ of $\JJ^1_P=\HH^1(\b,\L)
\left(\JJ^1(\b,\L)\cap P\right)$ by $\eta_P=\bs g(\teta_{\tP})$. 
Then Corollary~\ref{etaP} and Lemma~\ref{intetaP}, together with the 
properties of the Glauberman correspondence, give us

\begin{Lemma}\label{etaPforG}
The representation $\eta_P$ is the unique irreducible representation of
$\JJ^1_P$ such that $\eta_P|_{\HH^1_P}$ contains
$\th_P$. We have $\eta\simeq \Indu{\JJ^1_P}{\JJ^1(\b,\L)}{\eta_P}$. 
Moreover, for each
$y\in G^+_E$, the coset 
$\JJ^1_Py\JJ^1_P$ is the unique $(\JJ^1_P,\JJ^1_P)$-double
coset in $\JJ^1(\b,\L)y\JJ^1(\b,\L)$ which intertwines $\eta_P$ and
$$
\dim I_g(\eta_P)=\begin{cases} 1&\hbox{ if }g\in \JJ_P^1 G^+_E
\JJ_P^1;\\
0&\hbox{ otherwise}. \end{cases}
$$
\end{Lemma}

\begin{proof}
The statement about the dimensions of the intertwining spaces comes
from the corresponding statement for $\eta$, which is Proposition~\ref{eta}.
\end{proof}

As above we can describe $\eta_P$ as being the natural representation of
$\JJ^1_P$ on the space of $(\JJ^1(\b,\L)\cap U)$-fixed vectors in
$\eta$ and the restriction of $\eta_P$
to $\JJ^1_P\cap M=\JJ^1(\b,\L)\cap M$ is given by
$$
\eta_P|_{\JJ^1(\b,\L)\cap M} \simeq \eta^{(0)} 
\otimes \bigotimes_{j=1}^m \teta^{(j)},
$$
where $\eta^{(0)}$ is the unique irreducible representation of 
$\JJ^1(\b^{(0)},\L^{(0)})$ which contains the character
$\th^{(0)}$ and, for $j$ non-zero, 
$\teta^{(j)}$ is the unique irreducible representation of 
$\tJ^1(\b^{(j)},\L^{(j)})=\tJ^1(2\b^{(j)},\L^{(j)})$ which contains
the character $\left(\tth^{(j)}\right)^2$. 

\medskip

We continue in the same situation and suppose from now on that 
the self-dual decomposition $V=\bigoplus_{j=-m}^m
W^{(j)}$ is \emph{properly} subordinate to the skew semisimple stratum
$[\L,n,0,\b]$. Let $\k$ be a $\b$-extension of $\eta$ to $\JJ^+(\b,\L)$ 
(relative to some $\L^{\sM}$). We can
form the natural representation $\k_P$ of $\JJ^+_P$ on the space of
$(\JJ^+(\b,\L)\cap U)=(\JJ^1(\b,\L)\cap U)$-fixed vectors in $\k$. Then
$\k_P|_{\JJ^1_P}=\eta_P$ so, in particular, $\k_P$ is irreducible. 
As in the proof of~\cite[Proposition 7.2.15]{BK}, the Mackey restriction
formula gives:

\begin{Proposition}[cf.\ {\cite[Proposition 7.2.15]{BK}}]\label{kappaP}
We have $\Indu{\JJ^+_P}{\JJ^+(\b,\L)}{\k_P}\simeq\k$.
\end{Proposition}

As for $\eta_P$, we can describe the restriction of $\k_P$ to $\JJ^+_P\cap
M=\JJ^+(\b,\L)\cap M$. Since $\k_P|_{\JJ^+(\b,\L)\cap M}$ restricts 
further to
$\eta_P|_{\JJ^1(\b,\L)\cap M}$, it is irreducible and
$$
\k_P|_{\JJ^+(\b,\L)\cap M}=\k^{(0)}\otimes\bigotimes_{j=1}^m \til\k^{(j)},
$$
where $\k^{(0)}$ is some irreducible representation of
$\JJ^+(\b^{(0)},\L^{(0)})$ which extends $\eta^{(0)}$ and
$\til\k^{(j)}$ is some irreducible representation of
$\tJ(\b^{(j)},\L^{(j)})$ which extends $\til\eta^{(j)}$.

We will see below in Proposition~\ref{k0beta} that, at least in certain
circumstances, $\k^{(0)}$ is a $\b^{(0)}$ extension and $\til\k^{(j)}$
is a $2\b^{(j)}$-extension for $1\le j\le m$. 

\medskip

We can also make the same constructions for the restriction
$\k|_{\JJ^\so(\b,\L)}$, and we clearly obtain the restriction
$\k_P|_{\JJ^\so_P}$. Moreover, we have
$$
\Indu{\JJ^\so_P}{\JJ^\so(\b,\L)}{\left(\k_P|_{\JJ^\so_P}\right)}
\simeq\k|_{\JJ^\so(\b,\L)}
\qquad\hbox{ and }
\qquad\Indu{\JJ_P}{\JJ(\b,\L)}{\left(\k_P|_{\JJ_P}\right)}
\simeq\k|_{\JJ(\b,\L)}.
$$
Indeed, these are special cases of Lemma~\ref{HeckeP} below.
We also have
$$
\k_P|_{\JJ^\so(\b,\L)\cap M}=
\k^{(0)}|_{\JJ^\so(\b^{(0)},\L^{(0)})}\otimes\bigotimes_{j=1}^m \til\k^{(j)}.
$$


\Section{Intertwining and supercuspidal representations}\label{Sinter}

The aim of this section is two-fold. First we seek to improve the
lower bound for the
intertwining of $\b$-extensions from Corollary~\ref{intkappa}, in
particular for standard 
$\b$-extensions. In order to do this we need first to define certain Weyl
group elements. Secondly, we seek to bound above the intertwining of
representations of the form $\k|_{\JJ^\so}\otimes\rho$, where $\k$ is
a standard $\b$-extension and $\rho$ is a \emph{cuspidal}
representation of the reductive quotient $\JJ^\so/\JJ^1$. This will
use both the lower bound for the intertwining of $\k$ and, critically,
the result of Morris, Lemma~\ref{Morlem}.


\subsection{A Hecke algebra isomorphism}\label{Sinter.1}

Let $[\L,n,0,\b]$ be a skew semisimple stratum and let
$\th\in\CC_-(\L,0,\b)$ be a skew semisimple character. 
We use our usual notation. Let $\eta$ be the unique irreducible
representation of $\JJ^1(\b,\L)$ containing $\th$ and let $\k$ be a
$\b$-extension of $\eta$ (relative to some $\L^{\sM}$).

Let $V=\bigoplus_{j=-m}^m W^{(j)}$ be a
self-dual decomposition which is properly subordinate to
$[\L,n,0,\b]$. As in~\S\ref{S4}, let $M^+$ be the Levi subgroup of $G^+$
which stabilizes the decomposition $V=\bigoplus_{j=-m}^m W^{(j)}$ and
let $P^+$ be a parabolic subgroup with Levi component $M^+$, and unipotent
radical $U$. We use all the other related notation from~\S\ref{S4}.

The following Lemma will be useful in several situations.

\begin{Lemma}\label{HeckeP}
Let $\KK$ be a compact open subgroup of $G$ with
$\JJ^1(\b,\L)\subset\KK \subset \JJ^+(\b,\L)$, such that $\KK$
has an Iwahori decomposition with respect to $(M^+,P^+)$. We write
$$
\KK\cap M^+ = \KK^{(0)}\times\prod_{j=1}^m \tK^{(j)}
$$
Let $\rho$ be
the inflation to $\KK$ of an irreducible representation of
$\KK/\JJ^1(\b,\L)$, and put $\l=\k|_\KK\otimes\rho$. Let $\l_P$
denote the natural representation of $\KK_P=\HH^1(\b,\L)\left(\KK\cap
P\right)$ on the space of $\JJ^1(\b,\L)\cap U$-fixed vectors in
$\l$. Then 
\begin{enumerate}
\item $\l_P$ is irreducible and $\l\simeq 
\Indu{\KK_P}{\KK}{\l_P}$.
\item We have $\l_P\simeq\k_P\otimes\rho$, considering $\rho$ as a
representation of $\KK_P/\JJ^1_P \cong \KK/\JJ^1(\b,\L)$.
\item $\l_P|_{\KK\cap M}=\l^{(0)}\otimes\bigotimes_{j=1}^m
\til\l^{(j)}$, where $\l^{(0)}=\k^{(0)}|_{\KK^{(0)}}\otimes\rho^{(0)}$ is a
representation of ${\KK^{(0)}}$ and
$\til\l^{(j)}=\til\k^{(j)}|_{\tK^{(j)}}\otimes\til\rho^{(j)}$ is a
representation of $\tK^{(j)}$, for $1\le j\le m$.
\item There is an algebra isomorphism
$$
\CH(G^+,\l_P) \cong \CH(G^+,\l)
$$
which preserves support: if $\phi\in\CH(G^+,\l)$ has
support $\KK y\KK$, for some $y\in G^+_E$, then the
corresponding function $\phi_P\in\CH(G^+,\l_P)$ has support
$\KK_P y\KK_P$. 
\end{enumerate}
\end{Lemma}

\begin{proof} (i), (ii) and (iii) are clear. We deduce from (ii) 
and~\cite[Proposition 4.1.3]{BK} that we have an algebra
isomorphism $\CH(G^+,\l_P) \cong \CH(G^+,\l)$. Moreover, if
$\phi\in\CH(G^+,\l)$ has support $\KK y\KK$, for some 
$y\in G^+_E$, then~\cite[Corollary 4.1.5]{BK} shows that 
the corresponding function $\phi_P\in\CH(G^+,\l_P)$ has
support on a union of double cosets $\KK_P x\KK_P$ contained in $\KK
y\KK$. Using the Iwahori decomposition for $\KK$, we may assume further
that $x\in(\KK\cap U_l)y(\KK\cap U_l) \subset
\JJ^1(\b,\L)y\JJ^1(\b,\L)$, where $U_l$ is the unipotent radical of
the parabolic subgroup opposite $P^+$.

Since $\l_P|_{\JJ^1_P}$ is a multiple of $\eta_P$, we also have that
$x$ intertwines $\eta_P$. Meanwhile, by
Lemma~\ref{etaPforG}, $\JJ^1_P y\JJ^1_P$ is the unique
$(\JJ^1_P,\JJ^1_P)$-double coset in $\JJ^1(\b,\L) y\JJ^1(\b,\L)$
which intertwines $\eta_P$. Hence $x\in\JJ^1_P y\JJ^1_P$ and, since
$\JJ^1_P\subset\KK_P$, we have $\KK_P x\KK_P=\KK_P y\KK_P$. Hence
$\phi_P$ has support $\KK_P y\KK_P$, as required.
Since
$$
\Supp\CH(G^+,\l)\subset I_{G^+}(\l|\JJ^1(\b,\L)) = I_{G^+}(\eta|\JJ^1(\b,\L))
=\JJ^1(\b,\L) G^+_E \JJ^1(\b,\L) = K G^+_E K,
$$
we are done.
\end{proof}

\begin{Corollary}\label{intkappaP}
Let $\k$ be a $\b$-extension of $\eta$ to $\JJ^+(\b,\L)$ (relative to $\L^{\sM}$). 
Then $\PP^+(\L^\sM_{\oe})\subset I_G(\k_P|\JJ^+_P)$.
\end{Corollary}

\begin{proof} 
Let $y\in\PP^+(\L^\sM_{\oe})$. By Corollary~\ref{intkappa}, $y$
intertwines $\k$ so the result follows by applying Lemma~\ref{HeckeP}
with $\KK=\JJ^+(\b,\L)$ and $\l=\k$. 
\end{proof}

We write 
$$
\k_P|_{\JJ^+(\b,\L)\cap M}=\k^{(0)}\otimes\bigotimes_{j=1}^m \til\k^{(j)},
$$
as before.

\begin{Proposition}[cf.\ {\cite[Corollary 7.2.16]{BK}}]\label{k0beta}
\begin{enumerate}
\item
Suppose $\fb_0(\L^{(0)})$ is a maximal self-dual $\oe$-order in
$B^{(0)}$. Then $\k^{(0)}$ is a $\b^{(0)}$-extension of
$\eta^{(0)}$ to $\JJ^+(\b^{(0)},\L^{(0)})$.
\item
If $1\le j\le m$, $W^{(j)}\subset V^i$, for some $i$, and
$\fb_0(\L^{(j)})$ is a maximal 
$\fo_{E_i}$-order in $B^{(j)}$, then 
$\til\k^{(j)}$ is a $2\b^{(j)}$-extension of
$\til\eta^{(j)}$ to $\tJ(\b^{(j)},\L^{(j)})$.
\end{enumerate}
\end{Proposition}

Note that this Proposition is clear when $|\kf|>3$, since then
\emph{any} extension of $\eta^{(0)}$ to $\JJ^+(\b^{(0)},\L^{(0)})$ is
a $\b^{(0)}$-extension in this maximal case and, for $j\ne 0$, any extension of
$\til\eta^{(j)}$ to $\tJ(\b^{(j)},\L^{(j)})$ is a
$2\b^{(j)}$-extension (the latter only requires $|\kf|>2$). 

\begin{proof} 
For $i=1,...,l$, let $\L^i_\sm$ be an $\fo_{E_i}$-lattice sequence in
$V^i$ such that $\fb_0(\L^i_\sm)$ is a minimal self-dual $\fo_{E_i}$-order
contained in $\fb_0(\L^i)$ and such that the self-dual decomposition
$$
V^i=\bigoplus_{j=-m}^m V^i\cap W^{(j)}
$$
splits $\L^i_\sm$. Put $\L^\sm=\bigoplus_{i=1}^l \L^i_\sm$, a
self-dual $\oe$-lattice sequence in $V$ such that $\PP(\L^\sm_{\oe})$
is an Iwahori subgroup contained in $\PP(\L_{\oe})$ and such that the
self-dual decomposition $V=\bigoplus_{j=-m}^m W^{(j)}$ is 
(properly) subordinate
to $[\L^\sm,n_\sm,0,\b]$.

For $j=-m,...,m$, set $E^{(j)}=F[\b^{(j)}]$. We also form the
$\fo_{E^{(j)}}$-lattice sequence $\L^{\sm^{(j)}}$ in $W^{(j)}$, given
by 
$$
\L^{\sm^{(j)}}(r)=\L^\sm(r)\cap W^{(j)}, \qquad\hbox{for }r\in\BZ.
$$
Then $\fb_0(\L^{\sm^{(j)}})$ is a minimal (self-dual, in the case
$j=0$) $\fo_{E^{(j)}}$-order contained in $\fb_0(\L^{(j)})$. In
particular, we can think of the pair $(\L^{\sm^{(j)}},\L^{(j)})$ in
the same way as the pair $(\L^\sm,\L)$, and make the constructions
analogous to those below in these cases too.

\medskip

We form the group $\JJ_{\sm,\L}^1=\PP_1(\L^\sm_{\oe})\JJ^1(\b,\L)$ and
let $\eta_{\sm,\L}$ be the unique extension of $\eta$ to
$\JJ^1_{\sm,\L}$ which is intertwined by $G_E^+$ 
(using Corollary~\ref{uniqueetatilde}). Recall that, by
Proposition~\ref{kapparesetatilde}, we have
$\k|_{\JJ^1_{\sm,\L}}=\eta_{\sm,\L}$. 

Notice that, since $\PP_1(\L^\sm_{\oe})\subset \PP(\L_{\oe})$ and
$\PP(\L_{\oe})\cap U=\PP_1(\L_{\oe})\cap U$, we have $\JJ^1_{\sm,\L}\cap U
= \JJ^1(\b,\L)\cap U$. Then we can form the natural representation
$\eta^{\sm,\L}_P$ of $\JJ^1_{\sm,\L,P}= 
H^1(\b,\L)\left(\JJ^1_{\sm,\L}\cap P\right)$ 
on the space of $\JJ^1(\b,\L)\cap U$-fixed vectors in
$\eta_{\sm,\L}$. By applying Lemma~\ref{HeckeP}, we have
$$
\eta_{\sm,\L}\simeq\Indu{\JJ^1_{\sm,\L,P}} 
{\JJ^1_{\sm,\L}}{\eta^{\sm,\L}_P}.
$$
It is irreducible and we can describe its restriction to
$$
\JJ^1_{\sm,\L}\cap M= \JJ^1_{\sm^{(0)},\L^{(0)}} \times \prod_{j=1}^m
\tJ^1_{\sm^{(j)},\L^{(j)}}, 
$$
where $\JJ^1_{\sm^{(0)},\L^{(0)}}=
\PP_1(\L^{\sm^{(0)}}_{\fo_{E^{(0)}}})\JJ^1(\b^{(0)},\L^{(0)})$ and 
$\tJ^1_{\sm^{(j)},\L^{(j)}}=
\tP_1(\L^{\sm^{(j)}}_{\fo_{E^{(j)}}})\tJ^1(\b^{(j)},\L^{(j)})$, for
$j=1,...,m$. We have
$$
\eta^{\sm,\L}_P |_{\JJ^1_{\sm,\L}\cap M}\ =\  
\hat\eta\ =\ \hat\eta^{(0)}\otimes\bigotimes_{j=1}^m \hat{\teta}^{(j)},
$$
where $\hat\eta^{(0)}$ is the restriction of $\k^{(0)}$, and
$\hat{\teta}^{(j)}$ the restriction of $\til\k^{(j)}$, for
$j=1,...,m$.

By Lemma~\ref{HeckeP}(iv), $G^+_E$ intertwines $\eta^{\sm,\L}_P$
so certainly $G_E^+\cap M$ intertwines $\hat\eta=
\eta^{\sm,\L}_P |_{\JJ^1_{\sm,\L}\cap M}$. Then the uniqueness
in Corollary~\ref{uniqueetatilde} shows that 
$\hat\eta^{(0)}$ is the unique extension $\eta_{\sm^{(0)},\L^{(0)}}$
of $\eta^{(0)}$ to $\JJ^1_{\sm^{(0)},\L^{(0)}}$ which is intertwined
by all of $G^+_{E^{(0)}}=G^+\cap\Aut_{E}(W^{(0)})$. But then,
if
$\PP^+(\L^{(0)}_{\fo_{E}})$ is a maximal compact
subgroup of $G^+_{E^{(0)}}$, the construction of
$\b$-extensions in Proposition~\ref{betamax} shows that $\k^{(0)}$ is a
$\b^{(0)}$-extension of $\eta^{(0)}$, since it restricts to
$\eta_{\sm^{(0)},\L^{(0)}}$. The same argument works for
$\til\k^{(j)}$, for $j=1,..,m$, using the uniqueness in 
Proposition~\ref{etatildetG}, whenever $W^{(j)}\subset V^i$ and 
$\tP(\L^{(j)}_{\fo_{E_i}})$ is a maximal compact
subgroup of $\tG_{E^{(j)}}=\Aut_{E_i}(W^{(j)})$.
\end{proof}

We also record the following for later use, which follows from the
application of Lemma~\ref{HeckeP} in the proof above:

\begin{Corollary}\label{etatildeP}
The intertwining of $\eta^{\sm,\L}_P$ is given by
$$
\dim I_g(\eta^{\sm,\L}_P)=
\begin{cases} 1&\hbox{ if }g\in \JJ^1_{\sm,\L,P} G^+_E
\JJ^1_{\sm,\L,P};\\
0&\hbox{ otherwise}. \end{cases}
$$
\end{Corollary}

Proposition~\ref{k0beta} motivates the following definition:

\begin{Definition}\label{exactdef}\rm
Let $[\L,n,0,\b]$ be a skew semisimple stratum. 
We say that a self-dual decomposition
$V=\bigoplus_{j=-m}^m W^{(j)}$ is \emph{exactly subordinate to
$[\L,n,0,\b]$} if it is minimal amongst all self-dual decompositions
which are properly subordinate to $[\L,n,0,\b]$ -- that is, no
refinement of the decomposition is properly subordinate to the stratum.
\end{Definition}

A properly subordinate decomposition is exactly
subordinate to $[\L,n,0,\b]$ 
if and only if $\fa_0(\L^{(0)})\cap B^{(0)}$ is a maximal self-dual 
$\oe$-order in $B^{(0)}$ and, for each $j\ne 0$,
there is an $i$ such that $W^{(j)}$ is contained in $V^i$ and   
$\fa_0(\L^{(j)})\cap B^{(j)}$ is a maximal $\fo_{E_i}$-order in $B^{(j)}$.
In particular, Proposition~\ref{k0beta} can
be applied to all blocks when the decomposition is exactly
subordinate.


\subsection{Some Weyl group elements}\label{S4.weyl}

Let $[\L,n,0,\b]$ be a skew semisimple stratum and let
$V=\bigoplus_{j=-m}^m W^{(j)}$ be a self-dual decomposition which is
subordinate (but not necessarily exactly subordinate) to the
stratum. We also assume, throughout this section, that each $W^{(j)}$
is contained in some $V^i$, for $j\ne 0$, and that 
$\fa_0(\L^{(j)})\cap B^{(j)}$ is a maximal $\fo_{E_i}$-order in $B^{(j)}$.
In particular, this hypothesis is satisfied whenever the decomposition
is exactly subordinate. 

It will be useful to
have a fixed ordering for the subspaces in the decomposition -- in
effect, this is giving a preferred choice of the parabolic subgroup
$P$ whose Levi component $M$ is the stabilizer of the decomposition.  
For each $j\ne 0$, write $i=i_j$ for the unique integer such
that $W^{(j)}\subset V^i$.
Since $\fa_0(\L^{(j)})\cap B^{(j)}$ is a maximal $\fo_{E_i}$-order,
for each $j\ne 0$ there is a unique integer $q_j=q_j(\L)$ such that
\begin{enumerate}
\item $-e_i/2\le q_j\le e_i/2$, where $e_i=e(\L^i|\fo_{E_i})$;
\item $\L(q_j)\cap W^{(j)}\supsetneq\L(q_j+1)\cap W^{(j)}$;
\item for $j>0$, $q_j>-e_i/2$; for $j<0$, $q_j<e_i/2$.
\end{enumerate}
By duality, we have $i_{-j}=i_j$ and $q_{-j}=-q_j$, for $j\ne 0$. 
In particular, by
swapping the numbering of $W^{(j)}$ and $W^{(-j)}$ if necessary, we
may, and do, assume that $q_j\ge 0$ for $j>0$. We put
$$
\nu(j) = \nu(j,\L) = ( i_j, q_j ) \in \BZ^2.
$$
Then, for $j>0$ we order the subspaces $W^{(j)}$ lexicographically
according to $\nu(j)$ -- that is, by reordering if necessary, we
assume that, for $0<j<k\le m$, either $i_j<i_k$ or $i_j=i_k$ and
$q_j\le q_k$. (If $\nu(j)=\nu(k)$ then we may order them either way.)
The order for $j<0$ is then determined by duality.

We remark also that this ordering depends on the ordering of the
pieces in the splitting $V=\bigoplus_{i=1}^l V^i$, which may be put in
any convenient order.

\medskip

We now need to define certain elements of $G^+_E$, which are essentially Weyl
group elements. We use the block decomposition for $A$ given by the
decomposition $V=\bigoplus_{j=-m}^m W^{(j)}$, ordering the blocks
according to $\nu(j)$. 
For each $j>0$ we choose an ordered $\fo_{E_i}$-basis $\CB^{(j)}=\{\bs
v_{j,1},...,\bs v_{j,d_j}\}$ for the lattice $\L^{(j)}(0)$ which
splits the lattice sequence $\L^{(j)}$, where $i=i_j$. (Indeed,
the second condition follows from the first since $\fb_0(\L^{(j)})$
is a maximal $\fo_{E_i}$-order.) Note that then $\CB^{(j)}\subset\L^{(j)}(q_j)\setminus\L^{(j)}(q_j+1)$.

Recall from~\S\ref{S2.1} that we have a nondegenerate $\ve$-hermitian
$E_i/E_{i,0}$-form $f_i$ on $V^i$ such that the duality of lattices in $V^i$
induced by $f_i$ coincides with that induced by the restriction of the
form $h$. On $W^{(j)}\oplus W^{(-j)}$, we have the nondegenerate
$\ve$-hermitian form which is the restriction of 
$f_i$, and we take the ordered basis
$\CB^{(-j)}=\{\bs v_{-j,1},...,\bs v_{-j,d_j}\}$ of $W^{(j)}$ for
which
$$
f_i(\bs v_{j,s},\bs v_{-j,t}) = \begin{cases} 
\varpi_i&\hbox{ if }s=t, \\ 0&\hbox{ otherwise,}\end{cases}
$$
where $\varpi_i$ is a (fixed) uniformizer for $E_i$ such that 
$\ov{\varpi_i}=(-1)^{e(E_i/E_{i,0})-1}\varpi_i$. By duality, this basis
is an $\fo_{E_i}$-basis for the lattice $\L^{(-j)}(1)$ and splits the
lattice sequence $\L^{(-j)}$; indeed $\CB^{(-j)}\subset\L^{(-j)}(e_i-q_j)\setminus\L^{(-j)}(e_i-q_j+1)$. Finally, for each $i=1,\ldots,l$, we choose a (possibly empty)
self-dual $E_i$-basis $\CB^{(i,0)}$ for $W^{(i,0)}=W^{(0)}\cap V^i$ which splits $\L^{(i,0)}=\boe^{(0)}\L^i$. 

Suppose $j,k\ne 0$, are such that $i_j=i_k=i$ and $\dim_{E_i} W^{(j)}
= \dim_{E_i} W^{(k)}$. Then we write $I_{j,k}$ for the element of
$B^i$ which sends the ordered basis $\CB^{(k)}$ to the ordered basis
$\CB^{(j)}$ and all other basis elements to $0$. Note that, in the case
$j=k$, we have $I_{j,j}=\boe^{(j)}$, 
the projection onto $W^{(j)}$. We also have the following properties,
which are immediate from the definitions:

\begin{Lemma}\label{Ijkproperties}
 For $j,k,s,t\ne 0$ such that $I_{j,k}$ and $I_{s,t}$ are defined, we have
\begin{enumerate}
\item $I_{j,k}I_{s,t}=\begin{cases} I_{j,t}&\hbox{ if }k=s; \\0&\hbox{
  otherwise}.\end{cases}\\$
\item if $j,k>0$ then $\ov{I_{j,k}} = I_{-k,-j}$ and $\ov{I_{j,-j}} =
  \ve (-1)^{e(E_i/E_{i,0})-1} I_{j,-j}$.
\item $I_{j,k}\in\begin{cases}\fb_{q_j-q_k}(\L)\setminus\fb_{q_j-q_k+1}(\L) &\hbox{ if }j,k>0,\\ \fb_{q_j-q_k-e_i}(\L)\setminus\fb_{q_j-q_k-e_i+1}(\L) &\hbox{ if }j>0,\ k<0.\end{cases}$ 
\end{enumerate}
\end{Lemma}

We remark that, if the decomposition is exactly subordinate then: for
$j> 0$, we have $0<q_j<e_i/2$: and, for $j\ne k$, if $i_j=i_k$ then $q_j\ne q_k$.

\medskip

We are now ready to define our Weyl group elements. If $j,k> 0$ are
distinct, such that $i_j=i_k=i$ and $\dim_{E_i} W^{(j)} = \dim_{E_i} W^{(k)}$,
then we put
$$
s_{j,k} = I_{j,k}+I_{k,j}+I_{-k,-j}+I_{-j,-k} +
\sum_{t\ne \pm j,\pm k} I_{t,t}.
$$
Then $s_{j,k}$ exchanges the block $\boe^{(j)}A\boe^{(j)}$ with 
$\boe^{(k)}A\boe^{(k)}$ and the block $\boe^{(-j)}A\boe^{(-j)}$ with 
$\boe^{(-k)}A\boe^{(-k)}$. Note also that Lemma~\ref{Ijkproperties}
implies that $\ov{s_{j,k}}=s_{j,k}=s_{j,k}^{-1}$, so that $s_{j,k}\in
G_E^+$; indeed, we have $s_{j,k}\in G_E$ since $\det_{B/E}(s_{j,k})=1$.

For $j>0$, we put
\begin{eqnarray*}
s_j &=& I_{-j,j} + \ve (-1)^{e(E_i/E_{i,0})-1} I_{j,-j} + \sum_{t\ne
  \pm j} I_{t,t},\\ 
s_j^\varpi &=& \varpi_i^{-1} I_{-j,j} + \ve \varpi_i I_{j,-j} + 
\sum_{t\ne \pm j} I_{t,t}.
\end{eqnarray*}
These elements exchange the blocks $\boe^{(j)}A\boe^{(j)}$ and
$\boe^{(-j)}A\boe^{(-j)}$. 
Again, Lemma~\ref{Ijkproperties} implies that $s_j,s_j^\varpi\in
G_E^+$. Indeed, except in the case when $E_i=F=F_0$, $\ve=1$ and
$\dim_F(W^{(j)})$ is odd, they lie in $G_E$.

\medskip

From Lemma~\ref{Ijkproperties}(iii), we also get:

\begin{Lemma}\label{sjkproperties}
If $j,k> 0$ are distinct, such that $i_j=i_k=i$ and $\dim_{E_i} W^{(j)} =
\dim_{E_i} W^{(k)}$, then:
\begin{enumerate}
\item $s_{j,k}\in \PP^+(\L_{\oe})$ if and only if $q_j(\L)=q_k(\L)$;
\item $\begin{cases}
s_j\in\PP^+(\L_{\oe})&\hbox{ if and only if }q_j(\L)=e_i/2;\\
s_j^\varpi\in\PP^+(\L_{\oe})&\hbox{ if and only if }q_j(\L)=0.
\end{cases}$
\end{enumerate}
\end{Lemma}

We continue in the same situation. Let $\L^\sM$ be an $\oe$-lattice sequence in $V$ such that
$\fb_0(\L^\sM)$ is a self-dual $\oe$-order containing
$\fb_0(\L)$. Then the decomposition $V=\bigoplus_{j=-m}^m W^{(j)}$ is subordinate to $\L^\sM$ and, for $j\ne 0$, $\fa_0(\L^{\sM})\cap B^{(j)}=\fa_0(\L)\cap B^{(j)}$ is a maximal $\fo_{E_i}$-order in $B^{(j)}$, where $i=i_j$; hence we can define the integers $q_j(\L^\sM)$, as above. In particular, this has the following consequences:
\begin{enumerate}
\item if $j>0$ then $\L^\sM(0)\cap W^{(j)}=\L(0)\cap W^{(j)}$;
\item if $j>0$ then $q_j(\L^\sM)\ge 0$;
\item if $j>k>0$ then $q_j(\L^\sM)\ge q_k(\L^\sM)$.
\end{enumerate}
Then (ii)--(iii) imply that the ordering of the terms in the decomposition $V=\bigoplus_{j=-m}^m W^{(j)}$ is also a good ordering for $\L^\sM$, and, having fixed such an ordering, (i) implies that the elements $s_{j,k}$, $s_j$ and $s_j^\varpi$ are the same for $\L$ and $\L^\sM$.

\begin{Definition}\rm
We say that two blocks $W^{(j)}$ and $W^{(k)}$, with $j,k\ne 0$, 
are \emph{companion with respect to $\L^\sM$} if both are contained
in the same subspace $V^i$, $\dim_{E_i}W^{(j)}=\dim_{E_i}W^{(k)}$ and,
for $r\in\BZ$,
$$
\L^\sM(r)\cap W^{(j)} = \L^\sM(r+1)\cap W^{(j)}
\quad\iff\quad
\L^\sM(r)\cap W^{(k)} = \L^\sM(r+1)\cap W^{(k)}.
$$
So, in the notation above, $W^{(j)}$ and $W^{(k)}$ are companion
with respect to $\L^\sM$ if and only if $i_j=i_k=i$,
$\dim_{E_i}W^{(j)}=\dim_{E_i}W^{(k)}$ and $q_j(\L^\sM)\equiv
q_k(\L^\sM) \pmod{e_i^\sM}$, where $e_i^\sM=e(\L^{\sM,i}|\fo_{E_i})$.
\end{Definition} 

Lemma~\ref{sjkproperties} then says that, for $j,k>0$ distinct, $W^{(j)}$ and
$W^{(k)}$ are companion with respect to $\L^\sM$ if and only if
$s_{j,k}\in\PP^+(\L^\sM_{\oe})$, while $W^{(j)}$ and
$W^{(-j)}$ are companion with respect to $\L^\sM$ if and only if
either $s_j$ or $s_j^\varpi$ lies in $\PP^+(\L^\sM_{\oe})$.

\medskip

Finally, suppose that $\L^\sM$ is such that $\fb_0(\L^{\sM,i})$ is a maximal self-dual $\fo_{E_i}$-order, for some (fixed) $i$. Then, with our numbering (in particular, $\L^\sM(r)^\#=\L^\sM(1-r)$), the lattices $\L^\sM(0)\cap V^i$ and $\L^\sM(1)\cap
V^i$ determine all the lattices in the image of $\L^{\sM,i}$ (they
are all $E_i$-multiples of these two lattices). Then we see that
$$
\L^{\sM,i}(0) \supset \L^{\sM,i}(1) = \cdots = 
\L^{\sM,i}(e^\sM_i/2) \supset \L^{\sM,i}(e^\sM_i/2+1) = \cdots = 
\L^{\sM,i}(e^\sM_i).
$$
[Note that the converse of this is not quite true: if precisely one of these containments is proper then $\fb_0(\L^{\sM,i})$ is necessarily maximal, but if both are proper then $\fb_0(\L^{\sM,i})$ can be the ``Siegel'' order.] 

In particular, for $j>0$ with $i_j=i$, we have $q_j(\L^\sM)\in\{0, e_i^\sM/2\}$ so $W^{(j)}$ and $W^{(-j)}$ are companion with respect to
$\L^\sM$, and 
$$
s_j^\varpi\in\PP^+(\L^\sM_{\oe})\ \iff\ q_j(\L^\sM)=0;\qquad
s_j\in\PP^+(\L^\sM_{\oe})\ \iff\ q_j(\L^\sM)=e_i^\sM/2.
$$


\subsection{Intertwining of beta extensions}\label{S4.beta}

We continue with the notation of the previous section and, throughout
this section and the next, we suppose that the self-dual decomposition
$V=\bigoplus_{j=-m}^m W^{(j)}$ is \emph{exactly} subordinate to the
skew semisimple stratum $[\L,n,0,\b]$, with the numbering as in 
the previous section. In particular, we are in a special case of the
situation of~\S\ref{S4.weyl}.

Using the Weyl group element $s_j$, we can define an involution $\s_j$ on
$\tG^{(j)}=\Aut_F(W^{(j)})$ as follows. Identifying $\tG^{(j)}$ with the
subgroup $\{(\ov g^{-1},g)\in\tG^{(-j)}\times\tG^{(j)}\}\subset M\cap
G$, we put, for $g\in\tG^{(j)}$,
$$
\s_j (g)\ =\ s_j g (s_j)^{-1}.
$$

\begin{Lemma}\label{matchingJ}
\begin{enumerate}
\item If $1\le j<k\le m$ and $W^{(j)}\cong W^{(k)}$ as $E_i$-vector spaces
(for some $i=i_j=i_k$) then conjugation by $s_{j,k}$ induces an isomorphism 
$\tJ(\b^{(j)},\L^{(j)})\cong\tJ(\b^{(k)},\L^{(k)})$.
\item For $1\le j\le m$, the group $\tJ(\b^{(j)},\L^{(j)})$ is stable
under the involution $\s_j$.
\end{enumerate}
\end{Lemma}

\begin{proof}
Suppose $1\le j\le m$. We can identify $\b^{(j)}$ with $\b_i$, for
$i=i_j$. 
Since the decomposition is exactly subordinate, the $\fo_{E_i}$-order
$\fb_0(\L^{(j)})$ is maximal in $\End_{E_i}(W^{(j)})$. Moreover,
$\L^{(j)}$ is the unique (upto translation of index) $\fo_{E_i}$-lattice
sequence of period $e(\L^i|\fo_{E_i})$ whose associated $\fo_{E_i}$-order is
$\fb_0(\L^{(j)})$. 

(i) If $W^{(j)}\cong W^{(k)}$ as $E_i$-vector spaces then conjugation
by $s_{j,k}$ identifies
an $\fo_{E_i}$-basis for $\L^{(j)}(0)$ with an $\fo_{E_i}$-basis for
$\L^{(k)}(0)$, so we get an isomorphism
$\fb_0(\L^{(j)})\cong\fb_0(\L^{(k)})$ of $\fo_{E_i}$-orders. By the
uniqueness above, this must then identify $\L^{(j)}$ with a translate
of $\L^{(k)}$. In particular
$\fa_r(\L^{(j)})$ is identified with $\fa_r(\L^{(k)})$, for all
$r\in\BZ$, and, since
$\b^{(j)}$ is identified with $\b^{(k)}$  via $\b_i$, we have an
isomorphism 
$\tJ(\b^{(j)},\L^{(j)})\cong\tJ(\b^{(k)},\L^{(k)})$.

(ii) As in (i), conjugation by $s_j$ 
stabilizes
$\tJ(\b^{(j)},\L^{(j)})\times\tJ(\b^{(-j)},\L^{(-j)})$, swapping the
two factors. Since $s_j\in G^+$, taking intersections with $G^+$ gives
the result.
\end{proof}

Let $\L^\sM$ be an $\oe$-lattice sequence in $V$ such that
$\fb_0(\L^\sM)$ is a maximal self-dual $\oe$-order containing
$\fb_0(\L)$ and let $\k$ be a $\b$-extension relative to $\L^\sM$. 
Let $M$ be the
Levi subgroup of $G$ which is the stabilizer of the decomposition
and let $P$ be a parabolic subgroup with Levi component $M$. We
form the representation $\k_P$ of $\JJ^\so_P$ and write
$$
\k_P|_{\JJ^\so_P\cap M}=\k^{(0)}\otimes_{j=1}^m \til\k^{(j)}.
$$
We recall that, since the decomposition is \emph{exactly} subordinate to the stratum, $\k^{(0)}$ is a $\b^{(0)}$-extension
and $\til\k^{(j)}$ is a $2\b^{(j)}$-extension, for $j=1,\ldots,m$, by Proposition~\ref{k0beta}.

\begin{Corollary}\label{matchingkappa}
\begin{enumerate}
\item If $1\le j<k\le m$ and $W^{(j)}$ is companion to
$W^{(k)}$ with respect to $\L^\sM$
then conjugation by $s_{j,k}$ induces an equivalence
$\til\k^{(j)}\simeq\til\k^{(k)}$. 
\item For $1\le j\le m$, conjugation by $s_j$ 
induces an equivalence
$\til\k^{(j)}\circ\s_j\simeq\til\k^{(j)}$. 
\end{enumerate}
\end{Corollary}

\begin{proof} (i) In this situation $s_{j,k}\in\PP^+(\L^\sM_{\oe})$ so, 
by Corollary~\ref{intkappaP},
$s_{j,k}$ intertwines $\k_P$. Moreover, it normalizes $\JJ^\so_P\cap
M$ so, since the restriction of $\k_P$ to this group is irreducible,
it normalizes $\k_P|_{\JJ^\so_P\cap M}=
\k_0\otimes\bigotimes_{j=1}^m\til\k^{(j)}$. The result is then clear,
from Lemma~\ref{matchingJ}.

(ii) If $s_j\in\PP^+(\L^\sM_{\oe})$ then the proof is the same as for
(i). Otherwise $s_j^\varpi\in\PP^+(\L^\sM_{\oe})$ so we get an
equivalence $\til\k^{(j)}\circ\s_j\simeq {}^{\varpi_i}\til\k^{(j)}$. But
$\til\k^{(j)}$ is normalized by $\varpi_i$ and the result follows.
\end{proof}

\begin{Corollary}
For $1\le j\le m$, there
exists a $2\b^{(j)}$-extension $\til\k^{(j)}$ which is $\s_j$-stable.
\end{Corollary}

\begin{Remark}\rm
Once we have such a $\s_j$-stable 
$\til\k^{(j)}$, it follows (using~\cite[Theorem~5.2.2]{BK}) that all $\s_j$-stable
$2\b^{(j)}$-extensions are of the form
$\til\k^{(j)}\otimes\til\chi\circ\det^{(j)}$, where:
\begin{itemize}
\item $\det^{(j)}:\End_{E_i}(W^{(j)})\to E_i$ is the determinant; and 
\item $\til\chi$ is the inflation of a $\s_i$-invariant character of
$\k_{E_i}^\times$, where $\s_i$ is the involution $x\mapsto
(\ov{x})^{-1}$ on $\k_{E_i}^\times$ and $\ov{\phantom a}$ is the Galois
involution on $k_{E_i}$ with fixed field $k_{E_{i,0}}$. 
\end{itemize}
If $E_i/E_{i,0}$ is ramified,
then there are only two such characters: the trivial character and the
quadratic character. If $E_i/E_{i,0}$ is unramified,
then there are $|k_{E_{i,0}}|+1$ such characters.
\end{Remark}

Since $\k_P|_{\JJ^\so_P\cap M}=\k^{(0)}\otimes_{j=1}^m \til\k^{(j)}$ and $\k_P$ is trivial on both $\JJ^\so_P\cap U$ and $\JJ^\so_P\cap U_l$, this observation gives us the following useful corollary:

\begin{Corollary}\label{betatwist}
Let $\L^\sM$, $\L^{\sM'}$ be $\oe$-lattice sequences in $V$ such that
$\fb_0(\L^\sM)$ and $\fb_0(\L^{\sM'})$ are maximal self-dual
$\oe$-orders containing $\fb_0(\L)$.
Let $\k$ be a $\b$-extension of $\eta$ relative to $\L^\sM$, and
let $\k'$ be a $\b$-extension of $\eta$ relative to
$\L^{\sM'}$. There are $\s_i$-invariant characters $\til\chi^{(j)}$ of
$\k_{E_i}^\times$ and a character $\chi^{(0)}$ of
$\JJ^\so(\b^{(0)},\L^{(0)})/\JJ^1(\b^{(0)},\L^{(0)})$ 
such that, writing $\chi=\chi^{(0)}\otimes\bigotimes_{j=1}^m
\til\chi^{(j)}\circ\det^{(j)}$, we have
$$
\k'=\Indu{\JJ^\so_P}{\JJ^\so(\b,\L)}{\left(\k_P\otimes\chi\right)}.
$$
Moreover, if $W^{(j)}$ is companion to $W^{(k)}$ with respect to
both $\L^\sM$ and $\L^{\sM'}$, then $\til\chi^{(j)}=\til\chi^{(k)}$.
\end{Corollary}

Corollary~\ref{matchingkappa} can also be applied to get more
information about the intertwining of $\k_P$ (and hence of
$\k$). 
Suppose that $\k$ is a standard $\b$-extension, so
$\L^\sM=\fM_\L$.  
If $j,k\ne 0$ and $W^{(j)}\simeq W^{(k)}$ as
$E_i$-spaces, for some $i$, then $W^{(j)}$ and $W^{(k)}$ are companion
with respect to $\fM_\L$. Moreover, by Lemma~\ref{sjkproperties}(ii), we
have $s_j\in\PP(\fM_{\L,\oe})$.

\smallskip

For each $i$, write $\CB^i=\CB^{(i,0)}\cup\bigcup_{j:i_j=i}\CB^{(j)}$, where $\CB^{(j)}$ is
the $E_i$-basis for $W^{(j)}$ chosen in the previous section and $\CB^{(i,0)}$ is the basis for $V^i\cap W^{(0)}$. Then
$\CB^i$ is an $E_i$-basis for $V^i$ which splits the lattice sequence
$\L^i$. We write $\til T_{E_i}$ for the maximal split torus in
$\tG_{E_i}$  which corresponds to the basis $\CB^i$ -- that is, 
$\til T_{E_i}$ is the stabilizer in $\tG_{E_i}$ of the decomposition 
$V^i=\bigoplus_{\bs v\in\CB^i} E_i\bs v$. 

We write $\til T_E=\prod_{i=1}^l \til T_{E_i}$ and put $T_E^+=\til
T_E\cap G_E^+$ and $T_E=\til T_E\cap G_E$. We also put
$N^+=N_{G^+_E}(T^+_E)$ and $N=N_{G_E}(T_E)=N^+\cap G$.

\begin{Proposition}\label{intkPstandard}
Let $w\in N^+$ be such that $w$ normalizes $\PP^\so(\L_{\oe})\cap
M$. Then $w$ normalizes $\k_P|_{\JJ^\so(\b,\L)\cap M}$ and intertwines
$\k_P|_{\JJ^\so_P}$.
\end{Proposition}

\begin{proof}
Suppose $w\in N$ normalizes $\PP^\so(\L_{\oe})\cap M$. In particular, it
then permutes the blocks of the decomposition $V=\bigoplus_{j=-m}^m
W^{(j)}$. Moreover, since $w$ is fixed by the involution $\s$, we have
that $wW^{(j)}=W^{(k)}$ if and only if $wW^{(-j)}=W^{(-k)}$. In
particular, this implies that $w$ stabilizes $W^{(0)}$ and its
orthogonal complement $\bigoplus_{j\ne 0} W^{(j)}$. But then the
$0$-component $w^{(0)}=\boe^{(0)}w\boe^{(0)}$ of $w$ normalizes
$\PP^\so(\L^{(0)}_{\oe})$ so, since $\fb_0(\L^{(0)})$ is a
maximal self-dual $\oe$-order in $B^{(0)}$, we have
$w^{(0)}\in\PP(\L^{(0)}_{\oe})$, and $w^{(0)}$ certainly normalizes
$\k^{(0)}$. 

The remaining part of $w$ can be written as the product of:
\begin{itemize}
\item a block diagonal element $z$ with $z^{(j)}=\boe^{(j)}z\boe^{(j)}$ a
non-zero scalar in $E_{i_j}$ for $1\le j\le m$; and 
\item a block permutation matrix $s$ which, in
turn, can be expressed as a product of the matrices $s_{j,k}$ and
$s_j$, for $1\le j<k\le m$. 
\end{itemize}
Since $\k$ is a standard $\b$-extension, Corollary~\ref{matchingkappa}
implies that $s$ normalizes $\k_P|_{\JJ^\so(\b,\L)\cap M}= 
\k^{(0)}\otimes\bigotimes_{j=1}^m\til\k^{(j)}$. On the other hand,
$z^{(j)}$ certainly normalizes $\til\k^{(j)}$, since it is a
$2\til\beta^{(j)}$-extension. Hence $w$ normalizes
$\k_P|_{\JJ^\so(\b,\L)\cap M}$, as required.

To show that $w$ also intertwines $\k_P|_{\JJ^\so_P}$, we
imitate~\cite[{Th\'eor\`eme 
2.19}]{Sec2}. Since $\JJ^\so(\b,\L)$ has an Iwahori decomposition with
respect to both $(M,P)$ and $(M,P^w)$, the group $\JJ^\so_P\cap
{}^w\JJ^\so_P$ also has an Iwahori decomposition with respect to
$(M,P)$. Now let $\Phi$ be a non-zero intertwining operator of
$\k_P|_{\JJ^\so(\b,\L)\cap M}$ with ${}^w\k_P|_{\JJ^\so(\b,\L)\cap
  M}$. For $g\in\JJ^\so_P\cap {}^w\JJ^\so_P$ we write $g=lmu$, with
$l\in U_l$, $m\in M$ and $u\in U$. Then
$$
\Phi\circ\k_P(g)=\Phi\circ\k_P(m)= {}^w\k_P(m)\circ\Phi.
$$
It remains only to check that $l,u\in \ker({}^w\k_P)$ or, equivalently,
that $l^w,u^w\in\ker(\k_P)$. 

By elementary row and column operations (cf.~\cite[Lemme 2.20]{Sec2}) 
we have $\tU^w=(\tU^w\cap \tU)(\tU^w\cap \tU_l)$ so, by uniqueness of
Iwahori decomposition, $U^w=(U^w\cap U)(U^w\cap U_l)$ and 
$u^w\subset U U_l$. Moreover, since $u^w\in\JJ^\so_P$, the uniqueness of
the Iwahori decomposition again implies $u^w\subset(\JJ^\so_P\cap U)
(\JJ^\so_P\cap U_l)\subset\ker(\k_P)$. The same argument applies for $l^w$.
\end{proof}

Write $N_\L=\{w\in N:w$~normalizes~$\PP^\so(\L_{\oe})\cap M\}$. 
It will be useful to understand this group in terms of the notation of
Morris in~\S\ref{S1}. If $\PP^\so(\L_{\oe})$ is a maximal parahoric subgroup
of $G_E$ (so the decomposition is just $V=W^{(0)}$ and $M=P=G$) then
we have $N_\L\subset\PP(\L_{\oe})$, since this is the 
normalizer of $\PP^\so(\L_{\oe})$ in $G_E$. Otherwise, the group
called $\fM_J$ in~\S\ref{S1} is just $M\cap G_E$ (because the
decomposition is exactly subordinate) and the group
$\fM_J\cap P_J$ is just $\PP^\so(\L_{\oe})\cap M$. In particular, our
group $N_\L$ is just the group written $N_N(\fM_J\cap P_J)$.


\subsection{Intertwining and supercuspidal types}\label{S3.2}

We continue with the notation of the previous section so that we have:
a self-dual decomposition $V=\bigoplus_{j=-m}^m W^{(j)}$ \emph{exactly}
subordinate to a skew semisimple stratum $[\L,n,0,\b]$, numbered as
usual; a semisimple character $\th\in\CC_-(\L,0,\b)$; a
\emph{standard} $\b$-extension $\k$ to $\JJ^\so(\b,\L)$ of the unique
irreducible representation $\eta$ of $\JJ^1(\b,\L)$ containing $\th$;
and a parabolic subgroup $P$ of $G$ whose Levi component $M$ is the
stabilizer of the decomposition.

The following proposition bounds the intertwining of a pair
$(\JJ^\so_P,\l_P)$, where $\l_P=\k_P\otimes\rho$ and $\rho$ is (the
inflation of) a
\emph{cuspidal} irreducible representation of $\JJ^\so_P/\JJ^1_P \cong
\PP^\so(\L_{\oe})/\PP^1(\L_{\oe})$.

\begin{Proposition}\label{boundintertwiningP}
We have $I_G(\l_P|\JJ^\so_P)\subset\JJ^\so_P N_\L \JJ^\so_P$.
\end{Proposition}

{\bf Remark} The proof, which is a variant of the proof 
of~\cite[Proposition 5.3.2]{BK}, is inspired by \cite[page 553]{Bl1}  
(using Proposition~\ref{Morlem} in place of \emph{op.\ cit.} Lemma 
4.2).

\begin{proof} 
Suppose $g\in G$ intertwines $\lambda_P=\kappa_P\otimes\rho$, so
that $g\in I_G(\eta_P|\JJ_P^1)=\JJ^\so_P G_E \JJ^\so_P$, as $\rho$ is
trivial on $\JJ^1_P$. Thus we may assume $g$ lies in $G_E$. 

Now let $\L^\sm$ be a self-dual $\oe$-lattice sequence such that
$\fb_0(\L^\sm_{\oe})$ is a minimal self-dual $\oe$-order contained in
$\fb_0(\L_{\oe})$ and put $\JJ^\so_{\sm,\L,P}=\PP^\so(\L^\sm_{\oe})\JJ^1_P$
and $\JJ^1_{\sm,\L,P}=\PP_1(\L^\sm_{\oe})\JJ^1_P$.
Then $G_E\cap \JJ^\so_P=\PP^\so(\L_{\oe})$ is a parahoric 
subgroup of $G_E$ containing the Iwahori subgroup
$\PP^\so(\L^\sm_{\oe})$. Thus we may further assume that $g$ is a
distinguished double coset representative for
$\PP^\so(\L_{\oe})\backslash G_E/\PP^\so(\L_{\oe})$ 
(see~\S\ref{S1}).

Since $\dim I_g(\eta_P, \JJ_P^1)=1$, we can imitate the proof 
of~\cite[Proposition 5.3.2]{BK} to show that any non-zero 
intertwining operator in $I_g(\k_P\otimes\rho,\JJ_P^\so)$ has the form
$S\otimes T$, with $S\in I_g(\eta_P, \JJ_P^1)$ and $T$ an
endomorphism of the space of $\rho$.

Now, by Lemma~\ref{etaPforG} and Corollary~\ref{etatildeP}, the
operator $S$ also intertwines the restriction 
$\kappa_P|_{\JJ_{\sm,\L,P}^1}$ so, again as in the proof 
of~\cite[Proposition 5.3.2]{BK}, it
follows that $T$ belongs to $I_g(\rho|\JJ_{\sm,\L,P}^1)$. In particular,
$g$ intertwines the restriction of $\rho$ to $\JJ_{\sm,\L,P}^1\cap
G_E=\PP_1(\L^\sm_{\oe})$. But $\PP_1(\L^\sm_{\oe})$ is the radical of
an Iwahori subgroup of $G_E$ so, by Proposition~\ref{Morlem} and the
remarks following it (see also the end of the previous section for a
translation of this into our notations here), $g\in N_\L$,
as required.
\end{proof}

\begin{Corollary}\label{tightboundinterP}
Put $N_{\L}(\rho)=\{g\in N_\L:{}^g\rho\simeq\rho\}$. Then
$I_G(\l_P|\JJ^\so_P)\subset\JJ^\so_P N_\L(\rho) \JJ^\so_P$.
\end{Corollary}

\begin{proof}
By Proposition~\ref{boundintertwiningP}, we need only show that 
$I_{N_\L}(\l_P|\JJ^\so_P)\subset\JJ_P^\so N_\L(\rho)
\JJ^\so_P$. 
So suppose $g\in N_\L$ intertwines $\l_P$. As in the proof of
Proposition~\ref{boundintertwiningP}, any non-zero 
intertwining operator in $I_g(\k_P\otimes\rho,\JJ_P^\so)$ has the form
$S\otimes T$, with $S\in I_g(\eta_P, \JJ_P^1)$ and $T$ an
endomorphism of the space of $\rho$.

Now, by Proposition~\ref{intkPstandard}, $g$ intertwines $\k_P$ so,
by Lemma~\ref{etaPforG}, the operator $S$ also intertwines  
$\kappa_P$ and it
follows that $T$ belongs to $I_g(\rho|\JJ^\so_P)$. In particular,
$g$ intertwines the restriction of $\rho$ to $\JJ^\so_P \cap
G_E=\PP^\so(\L_{\oe})$. Then, by Theorem~\ref{Morlem2}, $g\in 
\PP^\so(\L_{\oe}) N_\L(\rho) \PP^\so(\L_{\oe})$, as required.
\end{proof}

\medskip

In the case when $\PP(\L^\sM_{\oe})$ is a maximal compact subgroup of
$G_E$ we can now obtain supercuspidal representations from this
construction:

Let $[\L^\sM,n,0,\b]$ be a skew semisimple stratum with
$\fb_0(\L^\sM)$ a maximal self-dual $\oe$-order in $B$. Let $\th_\sM$
be a skew semisimple character in $\CC_-(\L^\sM,0,\b)$, let
$\eta_{\sM}$ be the unique irreducible representation of
$\JJ^1_\sM=\JJ^1(\b,\L^\sM)$ and let $\k_\sM$ be a $\b$-extension of
$\eta_\sM$ to $\JJ_\sM=\JJ(\b,\L^\sM)$.

Put $\l_\sM=\k_\sM\otimes\tau$, for $\tau$ an irreducible
cuspidal representation of $\JJ_\sM/\JJ^1_\sM\simeq
\PP(\L^\sM_{\oe})/\PP_1(\L^\sM_{\oe})$. Recall that
$\PP(\L^\sM_{\oe})/\PP_1(\L^\sM_{\oe})$ is not, in general, connected:
by an irreducible cuspidal representation we mean an irreducible
representation whose restriction to the
connected component contains an irreducible cuspidal representation
$\rho$. Write $\JJ_\sM^\so=\PP^\so(\L^\sM_{\oe})\JJ^1_\sM$.

\begin{Definition}\rm\label{maxsimtype}
We call a pair $(\JJ_\sM,\l_\sM)$ as above a \emph{maximal simple type
for $G$}.
\end{Definition}

\begin{Proposition}\label{boundintertwining}
Let $(\JJ_\sM,\l_\sM)$ be a maximal simple type for 
$G$. Then $I_G(\l_\sM|\JJ_\sM)=\JJ_\sM$.
\end{Proposition}

\begin{proof} 
Suppose $g\in G$ intertwines $\lambda_\sM=\kappa_\sM\otimes\tau$, so
that $g\in I_G(\eta_\sM|\JJ_\sM^1)=\JJ_\sM G_E \JJ_\sM$, as $\tau$ is
trivial on $\JJ^1_\sM$. Thus we may assume $g$ lies in $G_E$. Now
$\kappa_\sM|_{\JJ_\sM^\so}$ is irreducible and is normalized by
$\JJ_\sM$. If $\rho$ is an irreducible cuspidal component of
$\tau|_{\JJ_\sM^\so}$ then, by Clifford Theory, the restriction of
$\l_\sM$ to $\JJ_\sM^\so$ takes the form 
$$
\l_\sM|_{\JJ_\sM^\so} = m\sum_{p} \k_\sM\otimes\tau^p \simeq 
m\sum_p (\k_\sM\otimes\rho)^p,
$$
where $m\in\BN$ is a multiplicity and the sum is taken over a set of
representatives for
$\PP(\L^\sM_{\oe})/N_{\PP(\L^\sM_{\oe})}(\tau)$. Since $g$ certainly
intertwines $\l_\sM|_{\JJ_\sM^\so}$, we see that there exist $p_1,p_2\in
\PP(\L^\sM_{\oe})$ such that $p_1gp_2$ intertwines $\k_\sM\otimes\rho$
and, since $p_i\in \JJ_\sM$, we may assume that $g$ intertwines
$\k_\sM\otimes\rho$. But then, since $N_\L(\rho)\subset
N_{G_E}(\PP^\so(\L^\sM_{\oe})) = \PP(\L^\sM_{\oe})$,
Corollary~\ref{tightboundinterP} implies that $g\in\JJ_\sM$, as required.
\end{proof}

From Proposition~\ref{boundintertwining}, together
with~\cite[Proposition 1.5]{C} and~\cite[Proposition 5.4]{BK2}, we
immediately get: 

\begin{Corollary}\label{constructsc} 
Let $(\JJ_\sM,\l_\sM)$ be a maximal simple type in
$G$. Then $\pi=\cInd_{\JJ_\sM}^G \l_\sM$ is an irreducible supercuspidal
representation of $G$ and $(\JJ_\sM,\l_\sM)$ is a $[G,\pi]_G$-type.
\end{Corollary}

We will show in the following section that all irreducible supercuspidal
representations of $G$ arise from this construction.


\Section{Exhaustion}\label{S5}

In this final section, we show that any irreducible supercuspidal 
representation of $G$ contains a maximal simple type. In particular,
every irreducible supercuspidal 
representation of $G$ arises from our constructions. We do not, 
however, address any unicity issues -- for example, we do not have an 
intertwining implies conjugacy result for maximal simple types, so we
do not know when two apparently different maximal simple types might
induce to equivalent irreducible supercuspidal representations of $G$.


\subsection{A Hecke algebra injection}\label{5.0}

Let $[\L,n,0,\b]$ be a skew semisimple stratum and let
$\th\in\CC_-(\L,0,\b)$ be a skew semisimple character. 
We use our usual notation: $V=\bigoplus_{i=1}^l V^i$ is the
splitting associated to $[\L,n,0,\b]$, $E=F[\b]=\bigoplus_{i=1}^l E_i$,
$B=\End_E(V)=\bigoplus_{i=1}^l B^i$, $G_E=B\cap G=\prod_{i=1}^l G_{E_i}$,
etc.

Let $[\L',n',0,\b]$ be another skew semisimple stratum, with
$\fb_0(\L)\subset\fb_0(\L')$ and let $\th'=\tau_{\L,\L',\b}(\th)$. Let
$\k$ (respectively $\k'$) be a $\b$-extension of $\JJ^\so(\b,\L)$
(respectively $\JJ^\so(\b,\L')$) such that $\k,\k'$ are compatible. 
Let $\rho$ be the inflation of some irreducible representation of
$\JJ^\so(\b,\L)/\JJ^1(\b,\L)$ and put $\vth=\k\otimes\rho$. Set
$$
\CG^\so=\JJ^\so(\b,\L)/\JJ^1(\b,\L) \cong 
\prod_{i=1}^l \PP^\so(\L^i_{\fo_{E_i}})/\PP_1(\L^i_{\fo_{E_i}})
$$
and write $\rho=\otimes_{i=1}^l\rho_i$, for $\rho_i$ the inflation of an
irreducible representation of the finite reductive group 
$\CG^\so_i=\PP^\so(\L^i_{\fo_{E_i}})/\PP_1(\L^i_{\fo_{E_i}})$.

We put $\JJ^\so_{\L,\L'}=\PP^\so(\L_{\oe})\JJ^1(\b,\L')$; then we can
identify $\JJ^\so_{\L,\L'}/\JJ^1(\b,\L')$ with 
$\PP^\so(\L_{\oe})/\PP_1(\L'_{\oe})$, which has
$\PP^\so(\L_{\oe})/\PP_1(\L_{\oe})$ as a quotient. Hence
we can view $\rho$ as a representation of $\JJ^\so_{\L,\L'}$.
We define $\vth'$ to be the (irreducible) representation of
$\JJ^\so_{\L,\L'}$ given by
$$
\vth'=\left(\k'|_{\JJ^\so_{\L,\L'}}\right)\otimes\rho.
$$

\begin{Proposition}[cf.\ {\cite[Proposition 5.5.13]{BK}}]\label{vthprime}
There is a canonical algebra isomorphism
$$
\CH(G,\vth) \cong \CH(G,\vth')
$$
which preserves support: if $\phi\in\CH(G,\vth)$ is supported on
$\JJ^\so(\b,\L) y \JJ^\so(\b,\L)$, for some $y\in G_E$, then the corresponding
function $\phi'\in\CH(G,\vth')$ has support
$\JJ^\so_{\L,\L'} y \JJ^\so_{\L,\L'}$.
\end{Proposition}

\begin{proof}
Suppose first that $\fa_0(\L)\subset\fa_0(\L')$. Given the simple
intersection property from Lemma~\ref{inter}, the proof in this case
is identical to that of \cite[Proposition 5.5.13]{BK}.

\medskip

In general, by Lemma~\ref{lattcont}, there is an $\oe$-lattice
sequence $\L''$ with $\fb_0(\L'')=\fb_0(\L)$ and
$\fa_0(\L'')\subset\fa_0(\L')$. By Lemma~\ref{intermediate},
we have a sequence $\L_0=\L$, $\L_1,...,\L_t=\L''$ of
$\oe$-lattice sequences such that, for each $i$,
$\fb_0(\L_i)=\fb_0(\L)$ and either $\fa_0(\L_i)\subset\fa_0(\L_{i-1})$
or $\fa_0(\L_i)\supset\fa_0(\L_{i-1})$. For each $i$, let $\k_i$ be
the $\b$-extension of $\JJ^\so(\b,\L_i)$ compatible with $\k$ and put
$\vth_i=\k_i\otimes\rho$, where we use the fact that
$\JJ^\so(\b,\L_i)/\JJ^1(\b,\L_i)\cong \JJ^\so(\b,\L)/\JJ^1(\b,\L)$ to
view $\rho$ as a representation of $\JJ^\so(\b,\L_i)$. 

Applying the first case with $\{\L_{i-1},\L_i\}$ in place of
$\{\L,\L'\}$, we see that we have support-preserving Hecke algebra
isomorphisms 
$$
\CH(G,\vth_{i-1}) \cong \CH(G,\vth_i),
$$
for $1\le i\le t$. Again applying the first case, with $\L''=\L_t$ in place
of $\L$, we have a support-preserving Hecke algebra isomorphism
$$
\CH(G,\vth_t) \cong \CH(G,\vth').
$$
Putting these together gives the required isomorphism.
\end{proof}

\begin{Proposition}[cf.\ {\cite[(5.6.2), Lemma 5.6.3]{BK}}]
Thinking of $\rho$ as an irreducible representation of
$\PP^\so(\L_{\oe})$, there is an algebra isomorphism
$$
\CH(\JJ(\b,\L'),\vth') \cong \CH(\PP(\L'_{\oe}),\rho)
$$
which preserves support: if $\phi\in\CH(\JJ(\b,\L'),\vth')$ is
supported on $\JJ^\so_{\L,\L'} y \JJ^\so_{\L,\L'}$, for some $y\in 
\PP(\L'_{\oe})$, then the corresponding function $\phi_E$ has 
support $\PP^\so(\L_{\oe}) y \PP^\so(\L_{\oe})$.
\end{Proposition}

\begin{proof}As in \cite[Lemma 5.6.3]{BK}, the
irreducibility of $\k'|_{\PP(\L_{\oe})\JJ^1(\b,\L')}$ implies that we
have an algebra isomorphism $\CH(\JJ(\b,\L'),\rho)\cong
\CH(\JJ(\b,\L'),\vth')$, given as follows: map $\phi\in
\CH(\JJ(\b,\L'),\rho)$ to ${\k'}^\vee\otimes\phi\in
\CH(\JJ(\b,\L'),\vth')$, where ${}^\vee$ denotes the contragredient. 
The isomorphism with
$\CH(\PP(\L'_{\oe}),\rho)$ follows by reduction to
$\JJ(\b,\L')/\JJ^1(\b,\L')\cong \PP(\L'_{\oe})/\PP^1(\L'_{\oe})$.
\end{proof}

Let $V=\bigoplus_{j=-m}^m W^{(j)}$ be a
self-dual decomposition which is properly subordinate to
$[\L,n,0,\b]$. As in~\S\ref{S4}, let $M$ be the Levi subgroup of $G$
which stabilizes the decomposition $V=\bigoplus_{j=-m}^m W^{(j)}$ and
let $P$ be a parabolic subgroup with Levi component $M$, and unipotent
radical $U$. We use all the other related notation from~\S\ref{S4}.

By definition, the representation $\rho$ is trivial on $\JJ^1(\b,\L)$
so we can
apply Lemma~\ref{HeckeP} with $\KK=\JJ^\so(\b,\L)$ and $\l=\vth=\k\otimes\rho$ 
to obtain Iwahori factorization results for 
$\vth$. In particular, we get a support-preserving Hecke
algebra isomorphism
$$
\CH(G,\vth_P) \cong \CH(G,\vth).
$$

Putting all this together, we see that there is an injective algebra
map
\begin{equation}\label{injalgmap}
\CH(\PP(\L'_{\oe}),\rho)\cong\CH(\JJ(\b,\L'),\vth') \hookrightarrow
\CH(G,\vth') \cong \CH(G,\vth) \cong \CH(G,\vth_P)
\end{equation}
which preserves support: if $\phi_E\in\CH(\PP(\L'_{\oe}),\rho)$ has support
$\PP^\so(\L_{\oe}) y \PP^\so(\L_{\oe})$, for some
$y\in\PP(\L'_{\oe})$, then the 
corresponding function $\phi_P\in\CH(G,\vth_P)$ has support $\JJ^\so_P
y\JJ^\so_P$. 


\subsection{Covers}\label{5.1}

Let $\pi$ be an irreducible representation of $G$ and suppose that
there is a pair $([\L,n,0,\b],\th)$, consisting of a skew semisimple
stratum $[\L,n,0,\b]$ and a skew semisimple character 
$\th\in\CC_-(\L,0,\b)$, such that $\pi|_{\HH^1(\b,\L)}$ contains
$\th$.

Among all such pairs, for fixed $\b$, we choose
one for which the order $\fb_0(\L)=\fa_0(\L)\cap B$ is
minimal. Since there is a unique irreducible representation $\eta$ of
$\JJ^1=\JJ^1(\b,\L)$ containing $\th$, $\pi$ must also
contain $\eta$ and hence some irreducible representation $\vth$ of
$\JJ^\so=\JJ^\so(\b,\L)$ containing 
$\eta$. Since $\eta$ extends to $\JJ^\so$, we see that $\vth=\k\otimes\rho$,
where $\k$ is some standard $\b$-extension of $\eta$ and $\rho$ is the 
inflation of some irreducible representation of $\JJ^\so/\JJ^1$. As above,
we put 
$$
\CG^\so=\JJ^\so/\JJ^1 \cong 
\prod_{i=1}^l \PP^\so(\L^i_{\fo_{E_i}})/\PP_1(\L^i_{\fa_i})
$$
and write $\rho=\otimes_{i=1}^l\rho_i$, for $\rho_i$ the inflation of an
irreducible representation $\ov\rho_i$ of the finite reductive group 
$\CG^\so_i=\PP^\so(\L^i_{\fo_{E_i}})/\PP_1(\L^i_{\fo_{E_i}})$.

\begin{Lemma}\label{mincusp}
In the situation above, all the representations $\ov\rho_i$ are cuspidal. 
\end{Lemma}

\begin{proof} Suppose, for
contradiction, that $\ov\rho_1$ is not cuspidal. Then there exists a proper
parabolic subgroup $\CP^\so_1$ of $\CG^\so_1$, with unipotent radical $\CU_1$,
such that $\ov\rho_1|_{\CU_1}$ contains the trivial character. Let
${\L'}^1$ be a self-dual $\fo_1$-lattice sequence in $V^1$ such that
$\fb_0({\L'}^1)\subset\fb_0(\L^1)$ and the image of
$\PP^\so({\L'}^1_{\fo_1})$ under the quotient map
$\PP^\so(\L^1_{\fo_1})\to\CG_1$ is $\CP^\so_1$.
Such a lattice sequence ${\L'}^1$ is obtained as a refinement of $\L^1$:
$\CP^\so_1$ is defined by a self-dual $k_{E_1}$-flag in
$\L^1(0)/\L^1(1)$, $k_{E_1}$-flags in $\L^1(r)/\L^1(r+1)$, for $1\le
r<\frac e2$, where $e=e(\L^1|\fo_1)$, and a self-dual $k_{E_1}$-flag in
$\L^1(\lfloor\frac{e+1}2\rfloor)/\L^1(\lfloor\frac{e}2\rfloor+1)$;
then the image of ${\L'}^1$ consists of all the lattices which are
inverse images of these flags, their duals, and their $E_1$-scalar
multiples. 

We put $\L'={\L'}^1\oplus\bigoplus_{i=2}^l\L^i$ so that
$\fb_0(\L')\subset\fb_0(\L)$. By Lemma~\ref{lattcont}, we can choose
$\L''$ an $\oe$-lattice sequence so that
$\fb_0(\L'')=\fb_0(\L')$ and $\fa_0(\L'')\subset\fa_0(\L)$. 
Then $[\L'',n'',0,\b]$ is a skew
semisimple stratum in $A$, for a suitable integer $n''$. Let
$\th''=\tau_{\L,\L'',\b}(\th)$ be the transfer of $\th$ to $\HH^1(\b,\L'')$
and let $\eta''$ be the unique irreducible
representation of $\JJ^1(\b,\L'')$ which contains $\th''$. We show that
$\pi$ contains $\eta''$, and hence also $\th''$, contradicting the
minimality of $\fb_0(\L)$.

Since $\PP_1(\L''_{\oe})$ projects onto $\CU_1$ under the projection map
$\JJ^\so\mapsto\CG$, $\rho|_{\PP_1(\L''_{\oe})\JJ^1}$ contains the trivial
character so $\pi$ contains the restriction of $\k$ to
$\PP_1(\L''_{\oe})\JJ^1$. The lemma now follows from:

\begin{Lemma}[cf.\ {\cite[Lemma 8.1.6]{BK}}] \label{etakappelemma}
The representations of
$\PP_1(\L'')$ induced by $\eta''$ and by $\k|_{\PP_1(\L''_{\oe})\JJ^1}$ are
irreducible and equivalent to each other.
\end{Lemma}

The proof of Lemma~\ref{etakappelemma} is identical to that 
of~\cite[Lemma 8.1.6]{BK}. This also finishes the proof of 
Lemma~\ref{mincusp}.
\end{proof}

\medskip

We continue with the same situation, so that the semisimple character
$\th$ contained in $\pi$ has $\fb_0(\L)$ minimal for this
property and hence the representations $\rho_i$ are all cuspidal.
For the remainder of this section, we also make the following hypothesis:
\begin{enumerate}
\item[(H)] There is no skew semisimple stratum $[\L',n',0,\b]$ with 
$\fb_0(\L')$ a maximal self-dual $\oe$-order in $B$ such that 
$\pi$ contains a representation of the form $\k'\otimes\rho'$, where
  $\k'$ is a standard $\b$-extension and $\rho'$ is a cuspidal
  representation of $\PP^\so(\L'_{\oe})$.
\end{enumerate}
In particular, hypothesis (H) implies 
that $\fb_0(\L)$ is \emph{not} a maximal self-dual
$\oe$-order in $B$. Let $[\L',n',0,\b]$ be another skew semisimple stratum with
$\fM_{\L'}=\fM_{\L}$ and $\fb_0(\L')\supset \fb_0(\L)$. We let
$\th'=\tau_{\L,\L',\b}\th$ be the transfer of $\th$ to $\HH^1(\b,\L')$
and let $\eta'$ be the unique irreducible
representation of $\JJ^1(\b,\L')$ containing $\th'$. Let $\k'$ be the
unique standard $\beta$-extension of $\eta'$ compatible with $\k$ and let
$\rho'$ be the inflation to $\JJ^\so(\b,\L')$ of an irreducible component of 
$$
\Indu{\PP^\so(\L_{\oe})/\PP_1(\L'_{\oe})}{\PP^\so(\L'_{\oe})/\PP_1(\L'_{\oe})}
{\ov\rho} .
$$

\begin{Definition}\rm
With notation as above,
we say that a representation of the form $\vth'=\k'\otimes\rho'$
\emph{lies over\/} $\vth=\k\otimes\rho$.
\end{Definition}

Note that, thinking of $\rho$ as a representation of
$\JJ^\so_{\L,\L'}=\PP^\so(\L_{\oe})\JJ^1(\b,\L')$ trivial on
$\JJ^1(\b,\L')$, the representations lying over $\vth$ are precisely
the irreducible constituents of 
$$
\Indu{\JJ^\so_{\L,\L'}}{\JJ^\so(\b,\L')}
{\left({\k'|_{\JJ^\so_{\L,\L'}}\otimes\rho}\right)}
= \k'\otimes \Indu{\JJ^\so_{\L,\L'}}{\JJ^\so(\b,\L')}{\rho}.
$$

\begin{Lemma}[cf.\ {\cite[Proposition 8.3.5]{BK}}] \label{overlie}
Let $\vth$ and $\L'$ be as above and let $\pi$ be an irreducible smooth
representation of $G$. Then $\pi$ contains $\vth$ if and only if it
contains some $\vth'$ lying over $\vth$.
\end{Lemma}

\begin{proof} We can think of $\rho$ as a representation of
$\PP^\so(\L_{\oe})\PP_1(\L)$ trivial on $\PP_1(\L)$. Notice also that
$\Indu{\JJ^\so(\b,\L)}{\PP^\so(\L_{\oe})\PP_1(\L)}{\k}$ is irreducible and restricts
irreducibly to the representation $\Indu{\JJ^1(\b,\L)}{\PP_1(\L)}{\eta}$ of
$\PP_1(\L)$. Since $\rho$ is irreducible, but trivial on $\PP_1(\L)$,
we deduce that 
$\left(\Indu{\JJ^\so(\b,\L)}{\PP^\so(\L_{\oe})\PP_1(\L))}\k\right)\otimes\rho$ is also
irreducible. 

Suppose first that we also have $\fa_0(\L')\supset \fa_0(\L)$.
Then, by the compatibility of $\k$ and $\k'$, we have
\begin{eqnarray*}
\Indu{\JJ^\so(\b,\L)}{\PP^\so(\L_{\oe})\PP_1(\L)}{\vth} &\simeq&
\left(\Indu{\JJ^\so(\b,\L)}{\PP^\so(\L_{\oe})\PP_1(\L)}{\k}\right)\otimes\rho \\
&\simeq&
\left(\Indu{\JJ^\so_{\L,\L'}}{\PP^\so(\L_{\oe})\PP_1(\L)}
{\k'|_{\JJ^\so_{\L,\L'}}}\right) 
\otimes\rho \\
&\simeq& 
\Indu{\JJ^\so_{\L,\L'}}{\PP^\so(\L_{\oe})\PP_1(\L)}
{\left({\k'|_{\JJ^\so_{\L,\L'}}\otimes\rho}\right)}.
\end{eqnarray*}
Since these representations are irreducible, $\pi$ contains $\vth$ if
and only if it contains
$\k'|_{\JJ^\so_{\L,\L'}}\otimes\rho$. But this is the case if
and only if $\pi$ contains some representation lying over $\vth$.

\medskip

For the general case, by Lemma~\ref{lattcont}, there is an $\oe$-lattice
sequence $\L''$ with $\fb_0(\L'')=\fb_0(\L)$ and
$\fa_0(\L'')\subset\fa_0(\L')$. By Lemma~\ref{intermediate},
we have a sequence $\L_0=\L$, $\L_1,...,\L_t=\L''$ of
$\oe$-lattice sequences such that, for each $i$,
$\fb_0(\L_i)=\fb_0(\L)$ and either $\fa_0(\L_i)\subset\fa_0(\L_{i-1})$
or $\fa_0(\L_i)\supset\fa_0(\L_{i-1})$. For each $i$, let $\k_i$ be
the $\b$-extension of $\JJ^\so(\b,\L_i)$ compatible with $\k$ and put
$\vth_i=\k_i\otimes\rho$. 

Applying the first case with $\{\L_i,\L_{i-1}\}$ in place of
$\{\L,\L'\}$, we see that $\pi$ contains $\vth_i$ if and only if it
contains $\vth_{i-1}$. (Note that, when
$\fa_0(\L_i)\subset\fa_0(\L_{i-1})$, we have
$\JJ^\so_{\L_i,\L_{i-1}}=\JJ^\so(\b,\L_{i-1})$ so $\vth_{i-1}$ is the
only representation lying over $\vth_i$, and vice versa.)

Applying the first case again with $\L''=\L_t$ in place of $\L$, we see
that $\pi$ contains $\vth_t$ if and only if it contains some
representation $\vth'$ lying over $\vth_t$. Since
$\PP^\so(\L_{\oe})=\PP^\so(\L''_{\oe})$, these are precisely the
representations lying over $\vth=\vth_0$, and the result follows.
\end{proof}

We continue under the hypothesis (H), so that $\fb_0(\L)$ is not a maximal
self-dual $\oe$-order in $B$. Indeed, we can suppose moreover that
$\PP(\L_{\oe})$ is not a maximal compact subgroup of
$G_E$; for, if it were, then there would exist a self-dual $\oe$-lattice
sequence $\L'$ such that $\fb_0(\L')$ is a maximal self-dual
$\oe$-order in $B$ containing $\fb_0(\L)$ and
$\PP(\L'_{\oe})=\PP(\L_{\oe})$; then, by Lemma~\ref{overlie} and using
the notation there, $\pi$ would also contain the unique representation
$\vth'$ of $\JJ^\so(\b,\L')$ which lies over $\vth$, and we have
$\vth'=\k'\otimes\rho$, contradicting (H). 

We may also make another assumption from now on, namely that the
subspaces $V^i$ are numbered so that $\b_i\ne 0$, for $i>1$ -- that
is, if any block of the semisimple stratum $[\L,n,0,\b]$ is a null
stratum then it is the $1^{\rm st}$ block.

\medskip

Let $V=\bigoplus_{j=-m}^m W^{(j)}$ be a self-dual decomposition
exactly subordinate to the stratum $[\L,n,0,\b]$. Since the order
$\fb_0(\L)$ is not maximal self-dual, this decomposition is
non-trivial (that is $m\ge 1$). We use the usual numbering for exactly
subordinate decompositions (see~\S\ref{S4.weyl}). Then we have 
\begin{equation}\label{J/J1decomp}
\JJ^\so(\b,\L)/\JJ^1(\b,\L) \cong 
\PP^\so(\L^{(0)}_{\oe})/\PP_1(\L^{(0)}_{\oe}) \times
\prod_{j=1}^m \tP(\L^{(j)}_{\oe})/\tP_1(\L^{(j)}_{\oe})
\end{equation}
and we can write
$$
\ov\rho = \ov\rho^{(0)} \otimes\bigotimes_{j=1}^m \ov{\til\rho}^{(j)},
$$
where we can interpret $\til\rho^{(j)}$ as an irreducible representation of 
$\tP(\L^{(j)}_{\oe})$ or $\tJ(\b^{(j)},\L^{(j)})$, as appropriate.

Suppose that $0<j\le m-1$ is such that $i_j=i_{j+1}=i$ and
$\ov{\til\rho}^{(j)}\not\simeq\ov{\til\rho}^{(j+1)}$. 
Since the decomposition is
exactly subordinate, we have $0<q_j(\L)<q_{j+1}(\L)$ (see the remarks following Lemma~\ref{Ijkproperties}).
We write $e_i=e(\L^i|\fo_{E_i})$ and define another self-dual
$\oe$-lattice sequence $\L'$ in $V$, of the same $\of$-period as $\L$,
by
$$
\L'(ke_i+r)= \begin{cases} 
\L(ke_i+r) \qquad\qquad\quad\hbox{ if }-e_i/2<r\le -q_{j+1},\hbox{ or }
-q_j<r\le q_j,\hbox{ or }  q_{j+1}<r\le e_i/2;  \\ \\
\varpi_i^k \L^{(j)}(q_j)\oplus\varpi_i^k \L^{(j+1)}(q_{j+1}) 
\oplus\bigoplus_{t\ne j,j+1}\L^{(t)}(ke_i+r) 
\qquad\qquad\quad\hbox{ if }q_j<r\le q_{j+1};  \\ \\
\varpi_i^k \L^{(-j)}(1-q_j)\oplus\varpi_i^k \L^{(-j-1)}(1-q_{j+1}) 
\oplus\bigoplus_{t\ne -j,-j-1}\L^{(t)}(ke_i+r) \\
\hskip10cm\hbox{ if }-q_{j+1}<r\le -q_{j}. \end{cases}
$$
Then $[\L',n,0,\b]$ is a skew semisimple stratum, the decomposition is
exactly subordinate to this stratum, and
\begin{eqnarray*}
&q_t(\L')=q_t(\L)\quad\hbox{ for }0<t<j\hbox{ or } j+1<t<m,& \\
&q_j(\L')=q_{j+1}(\L),\qquad q_{j+1}(\L')=q_j(\L).
\end{eqnarray*}
In particular, for the standard numbering relative to $\L'$, the
$j^{\rm th}$ and $(j+1)^{\rm th}$ blocks have been exchanged (and
likewise for $-j$ and $-j-1$). Note also that we have
$\fM_\L=\fM_{\L'}$ so there is a unique standard $\b$-extension $\k'$
to $\JJ^\so(\b,\L')$ compatible with $\k$. 

We have 
$$
\JJ^\so(\b,\L')/\JJ^1(\b,\L') \cong \JJ^\so(\b,\L)/\JJ^1(\b,\L)
$$
(where the $j^{\rm th}$ and $(j+1)^{\rm th}$ terms in the decomposition~\eqref{J/J1decomp} have been interchanged) and we denote by $\rho'$ the inflation to $\JJ^\so(\b,\L')$ of $\ov\rho$.
We put $\vth'=\k'\otimes\rho'$.

\begin{Lemma}[cf.\ {\cite[Proposition 8.3.4]{BK}}]
The representation $\pi$ contains $\vth'$.
\end{Lemma}

\begin{proof}
We define a third self-dual $\oe$-lattice sequence $\L''$ in $V$, of
the same $\of$-period as $\L$, by
$$
\L''(ke_i+r)= \begin{cases} 
\L(ke_i+r) \qquad\quad\qquad\hbox{ if }-e_i/2<r\le -q_{j+1},\hbox{ or }
-q_j<r\le q_j,\hbox{ or }  q_{j+1}<r\le e_i/2;  \\ \\
\varpi_i^k \L^{(j)}(q_j)\oplus\varpi_i^k \L^{(j+1)}(q_{j}) 
\oplus\bigoplus_{t\ne j,j+1}\L^{(t)}(ke_i+r) 
\qquad\qquad\qquad\hbox{ if }q_j<r\le q_{j+1}; \\ \\
\varpi_i^k \L^{(-j)}(1-q_j)\oplus\varpi_i^k \L^{(-j-1)}(1-q_{j}) 
\oplus\bigoplus_{t\ne -j,-j-1}\L^{(t)}(ke_i+r) \\
\hskip10cm\hbox{ if }-q_{j+1}<r\le -q_{j}. \end{cases}
$$
Then $[\L'',n,0,\b]$ is a skew semisimple stratum but the decomposition
is no longer exactly (or even properly) subordinate. We have
\begin{eqnarray*}
&q_t(\L'')=q_t(\L)\quad\hbox{ for }0<t<j,\ j+1<t\le m,& \\
&q_j(\L'')=q_{j+1}(\L'')=q_{j+1}(\L).&
\end{eqnarray*}
Now the self-dual $\oe$-order $\fb_0(\L'')$ contains both $\fb_0(\L)$
and $\fb_0(\L')$. Moreover, $\PP^\so(\L_{\oe})/\PP_1(\L''_{\oe})$ is a
maximal parabolic subgroup in $\PP^\so(\L''_{\oe})/\PP_1(\L''_{\oe})$,
and $\PP^\so(\L'_{\oe})/\PP_1(\L''_{\oe})$ is the opposite parabolic.

Since $\ov{\til\rho}^{(j)}\not\simeq\ov{\til\rho}^{(j+1)}$, the induced
representation  
$$
\Indu{\PP^\so(\L_{\oe})/\PP_1(\L''_{\oe})}
{\PP^\so(\L''_{\oe})/\PP_1(\L''_{\oe})}
{\ov\rho}
$$
is irreducible and there is a unique irreducible representation
$\vth''$ of $\JJ^\so(\b,\L'')$ lying over $\vth$. By
Lemma~\ref{overlie}, $\pi$ therefore contains $\vth''$. But this 
representation also lies over $\vth'$ so, again by
Lemma~\ref{overlie}, $\pi$ contains $\vth'$.
\end{proof}

In particular, this means that we may repeat this process several
times and hence we may (and do) assume that:
\begin{itemize}
\item[(A)]\qquad if $k>j>0$ and $i_j=i_k$ (so $q_j<q_k$) then 
$\dim_{E_i}(W^{(j)}) \le \dim_{E_i}(W^{(k)})$.
\end{itemize}


\subsubsection{The case $\rho\not\simeq\rho\circ\s$}\label{5.2}

We suppose first that there is some index $k$ such that the
representation 
$\til\rho^{(k)}$ of $\tJ(\b^{(k)},\L^{(k)})$ is not fixed by the
involution $\s_k$ (see~\S\ref{S4.beta}). For $j>0$ we define $\til\rho^{(-j)}=\til\rho^{(j)}\circ\s_j$ (a representation of $\tJ(\b^{(j)},\L^{(j)})$) and set 
$$
J=\{-m\le j\le m: \til\rho^{(j)}\simeq\til\rho^{(k)}\}.
$$
Note that, since $\til\rho^{(k)}$ is not fixed by $\s_k$, if $j\in J$
then $-j\not\in J$. Put $-J=\{j:-j\in J\}$ and $J_0=\{j:\pm j\not\in
J\}$; then we set
$$
Y_1=\bigoplus_{j\in J} W^{(j)},\quad Y_0=\bigoplus_{j\in J_0}
W^{(j)},\quad Y_{-1}=\bigoplus_{j\in -J} W^{(j)}.
$$
The decomposition $V=Y_1\oplus Y_0\oplus Y_{-1}$ is then properly
subordinate to the stratum $[\L,n,0,\b]$. Let $M'$ be the Levi
subgroup of $G$ which is the stabilizer of this decomposition, and let
$P'=M'U'$ be a parabolic subgroup with Levi factor $M'$. Let $M$ be
the Levi subgroup of $G$ which stabilizes the decomposition
$V=\bigoplus_{j=-m}^m W^{(j)}$ and let $P=MU$ be a
parabolic subgroup with Levi factor $M$ such that $P\subset P'$. We
can then form the representation $\vth_P$ of $\JJ^\so_P$, as in Lemma~\ref{HeckeP}. 

\begin{Proposition}\label{nonsdcase}
The pair $(\JJ^\so_P,\vth_P)$ is a cover of 
$(\JJ^\so_P\cap M',\vth_P|_{\JJ^\so_P\cap M'})$. In particular, any
smooth irreducible representation $\pi$ of $G$ which contains $\vth$
is not supercuspidal.
\end{Proposition}

\begin{proof}
By the choice of the subspaces $Y_1,Y_0,Y_{-1}$, we have $N_\L(\rho)
\subset M'$, where we recall $N_\L(\rho)=\{g\in N_\L: g$~normalizes~$\rho\}$.
By Corollary~\ref{tightboundinterP}, we get $I_G(\vth_P)\subset \JJ^\so_P
M' \JJ^\so_P$ and it follows from~\cite[Theorem 7.2]{BK2} that any smooth irreducible representation $\pi$ of $G$ which contains $\vth_P$
is not supercuspidal. However, $\pi$ contains $\vth$ if and only if it contains $\vth_P$, by Lemma~\ref{HeckeP}. 
\end{proof}


\subsubsection{The case $\rho\simeq\rho\circ\s$}\label{5.3}

We suppose now that $\til\rho^{(j)}\circ\s_j\simeq\til\rho^{(j)}$, for
all $1\le j\le m$. Let $M$ be
the Levi subgroup of $G$ which stabilizes the decomposition
$V=\bigoplus_{j=-m}^m W^{(j)}$. 
We fix, throughout this section, a choice of parabolic subgroup $P$
with Levi component $M$, namely the stabilizer in $G$ of the self-dual flag
$$
\{0\} \subset W^{(-m)} \subset \cdots \subset \bigoplus_{j=-m}^k
W^{(j)} \subset \cdots \subset \bigoplus_{j=-m}^m W^{(j)} = V.
$$
Let $U$ bet he unipotent radical of $P$, and $U_l$ its opposite relative to $M$.

\medskip

According to the numbering of \S\ref{S4.weyl}, $i_m=l$ so that $W^{(m)}\subset V^l$. 
For $t=0,1$, we define a self-dual $\fo_{E_l}$-lattice sequence $\fM^l_t$
in $V^l$ as follows: if $t=1$, 
$$
\fM^l_1 (2 k+ r)=\begin{cases} 
\varpi_l^k \L^l(1-e_l/2)&\hbox{ for } r= 0,\\
\varpi_l^k \L^l(e_l/2)&\hbox{ for }r=1.\end{cases}
$$
if $t=0$ then recall that, again according to the numbering of \S\ref{S4.weyl}, $q_m\ge 1$ is the
greatest integer such that $\L^l(q_m)\supsetneq\L^l(e_l/2)$ and we put
$$
\fM^l_0 (2 k+ r)=\begin{cases} 
\varpi_l^k \L^l(1-q_m)&\hbox{ for } r= 0,\\
\varpi_l^k \L^l(q_m)&\hbox{ for }r=1;\end{cases}
$$
Then the $\fb_0(\fM^i_t)$ are maximal self-dual $\fo_{E_l}$-orders
containing $\fb_0(\L^i)$. For each $t$, we put
$\fM_t=\fM^i_t\oplus\bigoplus_{k\ne i}\L^k$. 

\medskip

We recall the elements $s_m$ and $s_m^\varpi$ of~\S\ref{S4.weyl}. Note
that $s_m\in\PP^+(\fM_0)$ and $s_m^\varpi\in\PP^+(\fM_1)$, by
Lemma~\ref{sjkproperties}.
We define $U^m=s_m U s_m^{-1}$ and $U^m_l=s_m U_l s_m^{-1}$. 

\begin{Lemma} \label{wjconj}
We have:
\begin{enumerate}
\item $\JJ^\so_P\cap M$, $\JJ^\so_P\cap U\cap U^m$ and
$\JJ^\so_P\cap U_l\cap U^m_l$ are stable under conjugation by $s_m$
and by $s_m^\varpi$; 
\item $s_m(\JJ^\so_P\cap U_l\cap U^m)s_m^{-1} \subset \JJ^\so_P\cap U\cap U_l^m$;
\item $\left(s_m^\varpi\right)^{-1}(\JJ^\so_P\cap U\cap U_l^m)s_m^\varpi 
\subset \JJ^\so_P\cap U_l\cap U^m$;
\end{enumerate}
In particular, $\JJ^\so_P s_m \JJ^\so_P s_m^\varpi \JJ^\so_P=\JJ^\so_P
s_m s_m^\varpi \JJ^\so_P$.
\end{Lemma}

\begin{proof} We prove the corresponding 
statements in $\tG$ and the Lemma follows by taking fixed points under $\s$.
Recall that, since the decomposition is properly subordinate,
$\fJ(\b,\L)$ can be decomposed as 
$$
\fJ(\b,\L)=\bigoplus_{-m\le j,k\le m} \boe^{(j)}\fJ(\b,\L)\boe^{(k)},
$$
and similarly for $\fH(\b,\L)$. We can also write, for example,
$$
\tU\cap s_m \tU s_m^{-1} = 1 + \sum_{-m<j<k<m} \boe^{(j)}A\boe^{(k)}.
$$
Then the properties from Lemma~\ref{Ijkproperties} imply
that the conjugations by $s_m$ and by $s_m^\varpi$ are trivial on 
$\tU\cap s_m \tU s_m^{-1}$ and (i) follows easily.

Parts (ii) and (iii) are rather similar so we will only address (iii).
We have 
$$
\tJ_{\tP}\cap \tU\cap s_m \tU_l s_m^{-1} = 1 + \sum_{-m<k\le m}
\boe^{(-m)}\fJ(\b,\L)\boe^{(k)}+ 
\sum_{-m<j< m} \boe^{(j)}\fJ(\b,\L)\boe^{(m)}.
$$
For $-m<k<m$, we have
\begin{eqnarray*}
(s_m^\varpi)^{-1} \boe^{(-m)}\fJ(\b,\L)\boe^{(k)} s_m^\varpi 
&=& \varpi_l I_{m,-m} \boe^{(-m)}\fJ(\b,\L)\boe^{(k)}I_{k,k} \\
&=& \boe^{(m)} \varpi_l I_{m,-m} \fJ(\b,\L) \boe^{(k)} \\
&\subset & \boe^{(m)} \fH(\b,\L) \boe^{(k)},
\end{eqnarray*}
since, by Lemma~\ref{Ijkproperties}, $\varpi_l I_{m,-m}\in\fb_1(\L)$
and, by~\cite[Lemma 3.11(ii)]{S4},
$\fb_1(\L)\fJ(\b,\L)\subset\fH(\b,\L)$. But $1+\boe^{(m)}
\fH(\b,\L) \boe^{(k)} \subset \tJ_{\tP}$, as required. The situation
for the other blocks is similar.

\medskip

The final assertion is immediate by the Iwahori decomposition of
$J_P^\so$: we have
$$
s_m \JJ^\so_P s_m^\varpi = (s_m (\JJ^\so_P\cap U_l)
s_m^{-1})  (s_m(\JJ^\so_P\cap M)s_m^{-1})
s_m s_m^\varpi 
((s_m^\varpi)^{-1} (\JJ^\so_P\cap U)s_m^\varpi) \subset 
\JJ^\so_P s_m s_m^\varpi \JJ^\so_P,
$$
by (i)--(iii).
\end{proof}

We put $\z=s_m s_m^\varpi $, so that $\JJ^\so_P s_m \JJ^\so_P
s_m^\varpi \JJ^\so_P=\JJ^\so_P \z \JJ^\so_P$.

\begin{Corollary}\label{zjconj}
For $k\ge 0$, we have:
\begin{enumerate}
\item $\z^{k}(\JJ^\so_P\cap U)\z^{-k} \subset \JJ^\so_P\cap U$;
\item $\z^{-k}(\JJ^\so_P\cap U_l)\z^{k} \subset \JJ^\so_P\cap U_l$;
\item $\z^k(\JJ^\so_P\cap M)\z^{-k}= \JJ^\so_P\cap M$.
\end{enumerate}
In particular, for $k_1,k_2\ge 0$, 
$$
\JJ^\so_P \z^{k_1} \JJ^\so_P \z^{k_2} \JJ^\so_P=\JJ^\so_P \z^{k_1+k_2}
\JJ^\so_P.
$$
\end{Corollary}

\begin{proof} This is immediate from Lemma~\ref{wjconj}. For example,
we have $\JJ^\so_P\cap U=(\JJ^\so_P\cap U\cap U^m)(\JJ^\so_P\cap U\cap
U_l^m)$ and certainly $\z^k$ normalizes the first factor here. But,
noting that $s_m,s_m^\varpi$ differ from their inverses only by an element
of $\JJ^\so_P\cap M$, by Lemma~\ref{wjconj}(iii),(ii) we have
$$
\z(\JJ^\so_P\cap U\cap U_l^m)\z^{-1} \subset  
s_m(\JJ^\so_P\cap U_l\cap U^m)s_m^{-1} \subset 
(\JJ^\so_P\cap U\cap U_l^m),
$$
as required.
\end{proof}

Now we will construct certain invertible elements of
$\CH(G,\vth_P)$. \emph{We assume for now that $s_m,s_m^\varpi$ are
elements of $G$.} We will treat the other case below.

For $t=0,1$, let $\k_t$ denote a $\b$-extension
of $\eta$ compatible with some standard $\b$-extension of $\JJ^\so(\b,\fM_t)$. By Corollary~\ref{betatwist}, we have
$$
\k_t \simeq \Indu{\JJ^\so_P}{\JJ^\so(\b,\L)}{\k_P\otimes\chi_t},
$$
for some self-dual character $\chi_t=\chi_t^{(0)}\otimes 
\bigotimes_{j=1}^m \chi_t^{(j)}\circ\det^{(j)}$. We write
$\rho_t=\rho\otimes\chi_t^{-1}$ so that $\vth=\k_t\otimes\rho_t$. 
Since the character $\chi_t^{(m)}$ is self-dual, we still have that
$\til\rho_t^{(m)}\simeq\til\rho_t^{(m)}\circ\s_m$ and, since it
factors through the determinant $\det^{(m)}$, the representation
$\til\rho_t^{(m)}$ is still cuspidal.

By \eqref{injalgmap}, we have 
support-preserving injective algebra maps 
$$
\CH(\PP(\fM_{t,\oe}),\rho_t) 
\into \CH(G,\vth_P).
$$
The former has structure
understood by the results of Morris~\cite[Theorem 7.12]{M1}. In particular: 
\begin{itemize}
\item when $t=0$, it has an invertible element with support
$\PP^\so(\L_{\oe}) s_m \PP^\so(\L_{\oe})$ and we denote by $T$ the image in
$\CH(G,\vth_P)$ of this, so that $T$ has support $\JJ^\so_P s_m \JJ^\so_P$.
\item when $t=1$, it has an invertible element with support
$\PP^\so(\L_{\oe}) s_m^\varpi \PP^\so(\L_{\oe})$ and we denote by
$T^\varpi$ the image in $\CH(G,\vth_P)$ of this, so that $T^\varpi$ has support $\JJ^\so_P s_m^\varpi\JJ^\so_P$;
\end{itemize}

Now we define $S=T*T^\varpi$. By
Lemma~\ref{wjconj}, it is supported on the single double coset
$\JJ^\so_P\z\JJ^\so_P$ and it is invertible, since both $T$ and
$T^\varpi$ are.

Finally, put $S'=S^{e(E_l/F)}$. 
Then $S'$ is invertible and, by Lemma~\ref{zjconj}, it has support
$\JJ^\so_P \z^{e(E_l/F)}\JJ^\so _P$. We also put
$Y_1=W^{(m)}$, $Y_0=\bigoplus_{j\ne\pm m}W^{(j)}$ and
$Y_{-1}=W^{(-m)}$.

\medskip

Now we turn to the case when \emph{$s_m$, $s_m^\varpi$ are not
elements of $G$}. Recall that this means that $E_l=F=F_0$, that
$\ve=1$, and that $\dim_{F}W^{(m)}$ is odd. In particular, this means
we have $\b_l=0$ so $i=i_m=1$ and $V^k\subset W^{(0)}$, for $k>1$, by
our ordering of the blocks $V^i$; in
particular, $\PP(\L^k_{\fo_{E_k}})$ is a maximal compact subgroup of
$G_{E_k}$, for $k>1$. 

Moreover, $\tJ(\b^{(m)},\L^{(m)})=\PP(\L^{(m)})$ and $\k^{(m)}$ is
just the trivial character (not the quadratic character, since that does not extend to $\PP^+(\L^{(m)})$). The involution $\s_m$ is transpose-inverse 
for the basis $\CB_m$ chosen to define $s_m$. In particular,
the conditions that $\rho^{(m)}$ be cuspidal and stable under $\s_m$
imply that $\dim_{F}W^{(m)}=1$ (see for example~\cite[Theorem
  7.1]{Adler}). But then, by our ordering of the subspaces $W^{(j)}$,
(see assumption (A) above) we have $\dim_F(W^{(j)})=1$, for all $j=1,...,m$.

We have two distinct cases here:
\begin{enumerate}
\item $V^1\cap W^{(0)}\ne \{0\}$. In this case there is an element
$p\in\PP^+(\L^{(0)}_{\oe})\setminus\PP(\L^{(0)}_{\oe})$ such that
$p^2=1$. (In fact, $p\in\PP^+(\L^{(l,0)}_{\oe})$, in the notation of~\S\ref{S4.1}--2.) Then we replace $s_m$ and $s_m^\varpi$ by $ps_m$ and
$ps_m^\varpi$ respectively, which lie in $G$. (Note that $p$ commutes with $s_m$ and $s_m^\varpi$.) The argument is then as
above and we obtain an invertible element $S'$ of $\CH(G,\vth_P)$ supported
on $J_P^\so \z J_P^\so$ (we have $e(E_1/F)=1$). We define the
subspaces $Y_1=W^{(m)}$, $Y_0=\bigoplus_{j\ne\pm m}W^{(j)}$ and
$Y_{-1}=W^{(-m)}$ as above also.

\item $V^1\cap W^{(0)}=\{0\}$. In this case, the first observation to
make is that $m>1$; for if $m=1$ then $G_{E_1}$ is a 2-dimensional
special orthogonal group, which is just $F^\times$, so
$\PP(\L^1_{\fo_1})$ is certainly a maximal compact subgroup of
$G_{E_1}$. But then, $\PP(\L_{\oe})$ is a maximal compact subgroup of
$G_{E}$, contradicting hypothesis (H).

Hence $m\ge 2$ and we may imitate the constructions above but with the
elements $s_m s_{m-1}$ and $s_m^\varpi s_{m-1}^\varpi$ in place of $s_m$
and $s_m^\varpi$ respectively. (Note that the element $\z$ does change
in this case.)
Then the analogue of Lemma~\ref{wjconj} holds in this situation (with
$U^m$ replaced by $s_ms_{m-1}Us_{m-1}^{-1}s_m^{-1}$) and the argument
is again as above. So we obtain an invertible element $S'$ of
$\CH(G,\vth_P)$ supported on $J_P^\so \z J_P^\so$. In this case we put
$Y_1=W^{(m)}\oplus W^{(m-1)}$, $Y_0=\bigoplus_{j\ne\pm m,\pm(m-1)}W^{(j)}$ and
$Y_{-1}=W^{(-m)}\oplus W^{(1-m)}$.
\end{enumerate}

\medskip

Finally, \emph{we suppose we are in any of the cases -- that is
  $s_m\in G$ or not} .
The decomposition $V=Y_1\oplus Y_0\oplus Y_{-1}$ is properly
subordinate to the stratum $[\L,n,0,\b]$. Let $M'$ be the Levi
subgroup of $G$ which is the stabilizer of this decomposition, and let
$P'=M'P$, a parabolic subgroup with Levi factor $M'$. We again form the representation $\vth_P$ of $\JJ^\so_P$.

\begin{Proposition}\label{sdcase}
The pair $(\JJ^\so_P,\vth_P)$ is a $G$-cover of $(\JJ^\so_P\cap
M',\vth_P|_{\JJ^\so_P\cap M'})$. In particular, any
smooth irreducible representation $\pi$ of $G$ which contains $\vth$
is not supercuspidal.
\end{Proposition}

\begin{proof}
Certainly $(\JJ^\so_P,\vth_P)$ has the Iwahori factorization
properties required, since $U'\subset U$. It remains only to observe
$S'\in\CH(G,\vth_P)$ is invertible with support $\JJ^\so_P
\z^{e(E_l/F)}\JJ^\so_P=\JJ^\so_P\z'\JJ^\so_P$, for $\z'$ 
a strongly $(P,\JJ^\so_P)$-positive element of the centre of $M'$. The result now follows as in Proposition~\ref{nonsdcase}, by appealing to~\cite[Theorem 7.2]{BK2} and Lemma~\ref{HeckeP}.
\end{proof}


\subsection{Conclusions}\label{5.4}

Let $\pi$ be any irreducible supercuspidal representation of
$G$. We
know from \cite[Theorem 5.1]{S4} that $\pi$ contains some semisimple character
$\th\in\CC_-(\L,0,\b)$, where $[\L,n,0,\b]$ is a skew semisimple
stratum. As in~\S\ref{5.1}, 
it contains some representation $\vth=\k\otimes\rho$ of
$\JJ^\so(\b,\L)$, with $\k$ a standard $\b$-extension and 
$\rho$ cuspidal. 

If hypothesis (H) were satisfied then either
Proposition~\ref{nonsdcase} or Proposition~\ref{sdcase} would give us
a contradiction. Hence $\pi$ contains some such $\vth=\k\otimes\rho$
with $\fb_0(\L)$ a maximal self-dual $\oe$-order in $B$. Then
$\pi$ contains some irreducible component $\l$ of
$\Indu{\JJ^\so(\b,\L)}{\JJ(\b,\L)}{\vth}$, and we have
$\l=\k\otimes\tau$, for $\k$ a $\b$-extension and $\tau$ some
irreducible component of
$\Indu{\JJ^\so(\b,\L)}{\JJ(\b,\L)}{\rho}$. Such a $\tau$ is cuspidal
so $(\JJ(\b,\L),\l)$ is a maximal simple type and, by
Corollary~\ref{constructsc}, $\cInd_\JJ^G\l$ is irreducible with
$\pi\simeq\cInd_\JJ^G\l$. We have proved our main result:

\begin{Theorem}
Let $\pi$ be an irreducible supercuspidal representation of $G$. Then
there exists a maximal simple type $(\JJ_\sM,\l_\sM)$, in the sense of
Definition~\ref{maxsimtype}, such that $\pi\simeq\cInd_{\JJ_\sM}^G\l_\sM$.
\end{Theorem}

In particular, this verifies, for the group $G$, the folklore
conjecture that all irreducible supercuspidal representations are
irreducibly compact-induced from compact (mod-centre) subgroups.


\bigskip

\small

Shaun Stevens \hbk\indent
School of Mathematics \hbk\indent
University of East Anglia \hbk\indent
Norwich NR4 7TJ \hbk\indent
United Kingdom 

\smallskip
email: Shaun.Stevens@uea.ac.uk

\end{document}